\documentclass[preprint]{elsarticle}

\usepackage[utf8]{inputenc}

\usepackage{amsmath}
\usepackage{amsfonts}
\usepackage{bm}
\usepackage{hyperref}
\usepackage{graphicx}
\usepackage{placeins}
\usepackage{multirow}
\usepackage{xcolor}
\usepackage{booktabs}
\usepackage{geometry}

\usepackage{tikz}

\usepackage{import}
\usepackage{pdfpages}
\usepackage{transparent}
\usepackage{xcolor}
\usepackage{amssymb}
\usepackage{empheq}
\usepackage{epstopdf} 
\usepackage{tocbasic}
\usepackage{caption}
\usepackage{subcaption}
\usepackage{siunitx}
\usepackage{todonotes}
\usepackage{nameref}
\usepackage[mathlines]{lineno}
\usepackage{etoolbox}
\usepackage[normalem]{ulem} % for strikeout

\usepackage{empheq}
\usepackage[most]{tcolorbox}
\newtcbox{\boxstyle}[1][]{%
    nobeforeafter, math upper, tcbox raise base, enhanced,
    colframe=white!25!black,
    colback=white!92!black,
    boxrule=1.2pt,
    #1}

\newcommand{\p}{\phantom{0}}
\newcommand{\e}[1]{\times 10^{#1}}

\newcommand{\figsize}{\small}
\newcommand{\sinc}{\operatorname{sinc}}

\newcommand{\R}[1]{}
\newcommand{\revone}[1]{#1}
\newcommand{\revtwo}[1]{#1}
\newcommand{\revboth}[1]{#1}
\newcommand{\souttoggle}[1]{}

\newcommand*\linenomathpatchAMS[1]{%
  \expandafter\pretocmd\csname #1\endcsname {\linenomathAMS}{}{}%
  \expandafter\pretocmd\csname #1*\endcsname{\linenomathAMS}{}{}%
  \expandafter\apptocmd\csname end#1\endcsname {\endlinenomath}{}{}%
  \expandafter\apptocmd\csname end#1*\endcsname{\endlinenomath}{}{}%
}
\expandafter\ifx\linenomath\linenomathWithnumbers
  \let\linenomathAMS\linenomathWithnumbers
  \patchcmd\linenomathAMS{\advance\postdisplaypenalty\linenopenalty}{}{}{}
\else
  \let\linenomathAMS\linenomathNonumbers
\fi
\linenomathpatchAMS{gather}
\linenomathpatchAMS{multline}
\linenomathpatchAMS{align}
\linenomathpatchAMS{alignat}
\linenomathpatchAMS{flalign}

\definecolor{gblue}{RGB}{51, 102, 204}
\definecolor{gred}{RGB}{204, 0, 0}
\definecolor{ggreen}{RGB}{102, 153, 0}
\definecolor{gyellow}{RGB}{255, 204, 0}
\definecolor{gorange}{RGB}{255, 153, 0}

\title{Discretize first, filter next: learning divergence-consistent closure models for large-eddy simulation}
\date{\today}
\journal{Journal of Computational Physics}

\begin{document}

% % Standalone graphical abstract
% \begin{center}
% \figsize
% \def\svgwidth{1\columnwidth}
% \input{./figures_graphabs.tex}
% \end{center}

% % Standalone novelty page
% \pagenumbering{gobble}
% \input{sections/novelty}

\begin{frontmatter}



% \title{}


% \author{}
% \author[cwi]{Syver Døving Agdestein}
% \author[cwi]{Benjamin Sanderse}
% \author[cwi,eindhoven]{Syver Døving Agdestein, Benjamin Sanderse}

\author[cwi,eindhoven]{Syver Døving Agdestein\corref{cor1}}
\ead{sda@cwi.nl}
\author[cwi,eindhoven]{Benjamin Sanderse}
\ead{b.sanderse@cwi.nl}

\cortext[cor1]{Corresponding author}

\affiliation[cwi]{
    organization={Scientific Computing Group, Centrum Wiskunde \& Informatica},%Department and Organization
    addressline={Science Park 123}, 
    city={Amsterdam},
    postcode={1098 XG}, 
    % state={Nord-Holland},
    country={The Netherlands}
}
\affiliation[eindhoven]{
    organization={Centre for Analysis, Scientific Computing and Applications, Eindhoven University of Technology},%Department and Organization
    addressline={PO Box 513}, 
    city={Eindhoven},
    postcode={5600 MB}, 
    country={The Netherlands}
}

\begin{abstract}
We propose a new neural network based large eddy simulation framework for the
incompressible Navier-Stokes equations based on the paradigm ``discretize first,
filter and close next''. This leads to full model-data consistency and allows
for employing neural closure models in the same environment as where they have
been trained. Since the LES discretization error is included in the learning
process, the closure models can learn to account for the discretization.

Furthermore, we
\revone{employ a}
divergence-consistent discrete filter defined through face-averaging
\revone{and provide novel theoretical and numerical filter analysis}.
\revone{ This}
filter preserves the discrete divergence-free constraint by construction, unlike
general discrete filters such as volume-averaging filters. We show that using a
divergence-consistent LES formulation coupled with a convolutional neural
closure model produces stable and accurate results for both a-priori and
a-posteriori training, while a general (divergence-inconsistent) LES model
requires a-posteriori training or other stability-enforcing measures.

\end{abstract}


\end{frontmatter}

% \linenumbers

% Sections
\section{Introduction} \label{sec:introduction}

The incompressible Navier-Stokes equations form a model for the movement of
fluids. They can be solved numerically on a grid using discretization techniques
such as finite differences~\cite{Harlow1965}, finite volumes, or pseudo-spectral
methods~\cite{Orszag1972,Rogallo1981}. The important dimensionless parameter in
the incompressible Navier-Stokes equations is the Reynolds number $\mathrm{Re} =
\frac{U L}{\nu}$, where $L$ is a characteristic length scale, $U$ is a
characteristic velocity, and $\nu$ the kinematic viscosity. For high Reynolds
numbers, the flow becomes turbulent. Resolving all the scales of motion of
turbulent flows requires highly refined computational grids. This is
computationally expensive~\cite{Pope2000,Sagaut2005,Rogallo1984,Verstappen1997}.

Large eddy simulation (LES) aims to resolve only the large scale features of the
flow, as opposed to direct numerical simulation (DNS), where all the scales are
resolved~\cite{Pope2000,Berselli2006}. The large scales of the flow, here
denoted $\bar{u}$, are extracted from the full solution $u$ using a spatial
filter. The equations for the large scales are then obtained by filtering the
Navier-Stokes equations. The large scale equations (filtered DNS equations) are
not closed, as they still contain terms depending on the small scales. It is
common to group the contributions of the unresolved scales into a single
commutator error term, that we denote $c(u)$. Large eddy simulation
requires modeling this term as a function of the large scales only. The common
approach is to introduce a closure model $m(\bar{u}, \theta) \approx c(u)$
to remove the small scale dependency
\cite{Pope2000,Sagaut2005,Berselli2006}, where $\theta$ are problem-specific
model parameters. The closure model accounts for the effect of the sub-filter
scales on the resolved scales. Traditional closure models account for the energy
lost to the sub-filter eddies. A simple approach to account for this energy
transfer is to add an additional diffusive term to the LES equations. These
closure models are known as \emph{eddy viscosity} models. This includes the well
known standard Smagorinsky model~\cite{Smagorinsky1963,Lilly1962,Lilly1967}.

The filter can be either be known explicitly or implicitly. In the latter case,
the (coarse) discretization itself is acting as a filter. Note that the
distinction explicit versus implicit is not always clear, as an explicitly
defined filter can still be linked to the discretization. Sometimes the filter
is only considered to be explicit when the filter width is fully independent
from the LES discretization size~\cite{Beck2023}. The type of filtering employed
results in a distinction between explicit and implicit LES. \emph{Explicit LES}
\cite{Lund2003,Gallagher2019a,Gallagher2019b} refers to a technique where the
filter is explicitly used inside the convection term of the LES equations (as
opposed to applying the filter in the exact filtered DNS evolution equations,
where the filter is applied by definition). The LES solution (that we denote by
$\bar{v}$) already represents a filtered quantity, namely the filtered DNS
solution $\bar{u}$. Applying the filter \emph{again} inside the LES equations
can thus be seen as computing a quantity that is filtered twice. The main
advantage is stability~\cite{Benjamin2023b} and to prevent the growth of
high-wavenumber components~\cite{Lund2003,Gallagher2019a}. Note that this
technique requires explicit access to the underlying filter. \emph{Implicit LES}
is used to describe the procedure where the DNS equations are used in their
original form, and the closure model is added as a correction (regardless of
whether the filter is known explicitly or not). Here, the filter itself is not
required to evaluate the LES equations, but knowledge of certain filter
properties such as the filter width is used to choose the closure model. In the
standard Smagorinsky model for example, the filter width is used as a model
parameter. If the exact filter is not known, this parameter is typically chosen
to be proportional to the grid spacing.

Recently, machine learning has been used to learn closure models, focusing
mostly on implicit LES
\cite{Beck2019,Beck2021,Beck2023,Duraisamy2019,Holl2020,Kurz2020,Kurz2021,List2022,Maulik2018,Sirignano2020,Sirignano2023a,Sirignano2023b}.
The idea is to represent the closure model $m(\bar{u},\theta)$ by an artificial
neural network (ANN). ANNs are in principle a good candidate as they are
universal function approximators~\cite{Barron1993,ChenTianping1995}. However,
using ANNs as closure models typically suffers from stability issues, which have
been attributed to a so-called model-data inconsistency: the environment in
which the neural network is trained is not the same in which it is being used
\cite{Beck2023}. Several approaches like backscatter clipping,
a-posteriori training and projection onto an eddy-viscosity basis have been used
to enforce stability
\revone{
    \R{inconsistency}
    \cite{Park2021,Kochkov2021,List2022,Beck2019}
}
--
for an overview, see~\cite{Sanderse2024}. Our view is that
\emph{one of the problems that lies at the root of the model-data inconsistency
is a discrepancy between the LES equations (obtained by first filtering the
continuous Navier-Stokes equations, then discretizing) and the training data
(obtained by discretizing the Navier-Stokes equations, and then applying a
discrete filter)}. 

\begin{figure}
    \centering
    \figsize
    \def\svgwidth{1\columnwidth}
    %% Creator: Inkscape 1.4 (e7c3feb100, 2024-10-09), www.inkscape.org
%% Accompanies image file 'figures_equations_same.pdf' (pdf, eps, ps)
\begingroup%
  \makeatletter%
  \providecommand\color[2][]{%
    \errmessage{(Inkscape) Color is used for the text in Inkscape, but the package 'color.sty' is not loaded}%
    \renewcommand\color[2][]{}%
  }%
  \providecommand\transparent[1]{%
    \errmessage{(Inkscape) Transparency is used (non-zero) for the text in Inkscape, but the package 'transparent.sty' is not loaded}%
    \renewcommand\transparent[1]{}%
  }%
  \providecommand\rotatebox[2]{#2}%
  \newcommand*\fsize{\dimexpr\f@size pt\relax}%
  \newcommand*\lineheight[1]{\fontsize{\fsize}{#1\fsize}\selectfont}%
  \ifx\svgwidth\undefined%
    \setlength{\unitlength}{680.31496063bp}%
    \ifx\svgscale\undefined%
      \relax%
    \else%
      \setlength{\unitlength}{\unitlength * \real{\svgscale}}%
    \fi%
  \else%
    \setlength{\unitlength}{\svgwidth}%
  \fi%
  \global\let\svgwidth\undefined%
  \global\let\svgscale\undefined%
  \makeatother%
  \begin{picture}(1,0.40042035)%
    \lineheight{1}%
    \setlength\tabcolsep{0pt}%
    \put(0,0){\includegraphics[width=\unitlength,page=1]{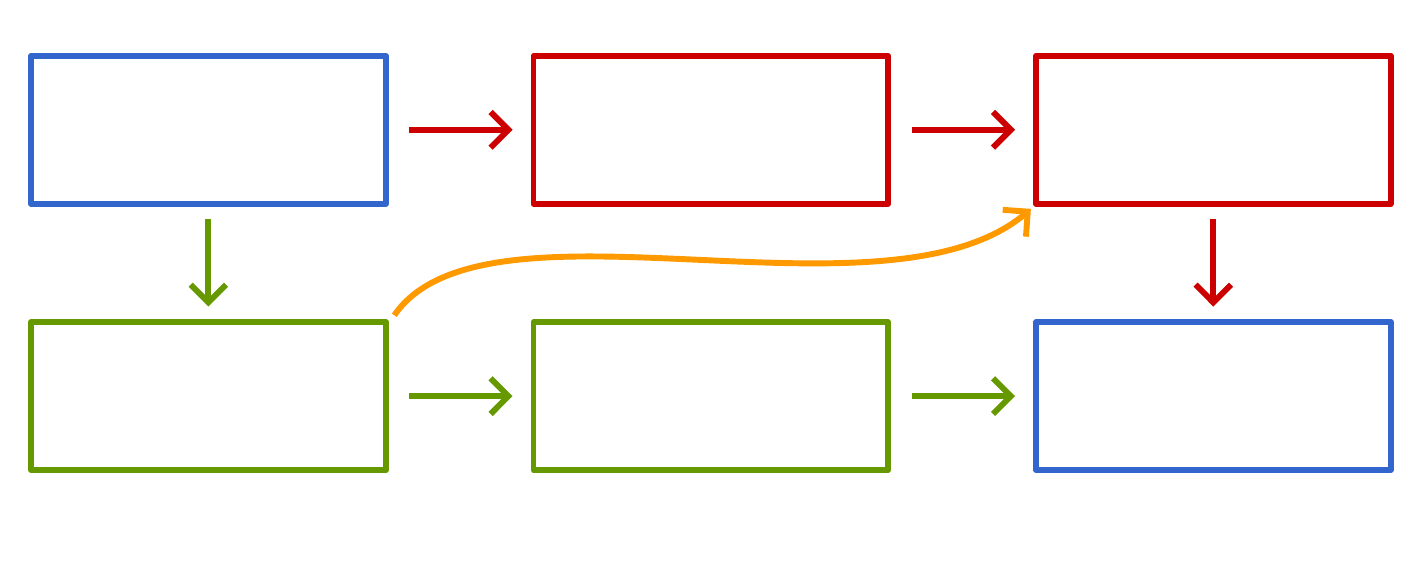}}%
    \put(0.09228117,0.32988149){\color[rgb]{0.2,0.4,0.8}\makebox(0,0)[lt]{\lineheight{0}\smash{\begin{tabular}[t]{l}$\nabla \cdot u = 0$\end{tabular}}}}%
    \put(0.07943198,0.37709932){\color[rgb]{0.2,0.4,0.8}\makebox(0,0)[lt]{\lineheight{0}\smash{\begin{tabular}[t]{l}Navier-Stokes\end{tabular}}}}%
    \put(0.06247483,0.28523745){\color[rgb]{0.2,0.4,0.8}\makebox(0,0)[lt]{\lineheight{0}\smash{\begin{tabular}[t]{l}$\dfrac{\partial u}{\partial t} = \mathcal{F}(u) - \nabla p$\end{tabular}}}}%
    \put(0.29539931,0.32991745){\color[rgb]{0.8,0,0}\makebox(0,0)[lt]{\lineheight{0}\smash{\begin{tabular}[t]{l}Filter\end{tabular}}}}%
    \put(0.39352338,0.37709932){\color[rgb]{0.8,0,0}\makebox(0,0)[lt]{\lineheight{0}\smash{\begin{tabular}[t]{l}Filtered Navier-Stokes\end{tabular}}}}%
    \put(0.78069907,0.37709932){\color[rgb]{0.8,0,0}\makebox(0,0)[lt]{\lineheight{0}\smash{\begin{tabular}[t]{l}Traditional LES\end{tabular}}}}%
    \put(0.55418399,0.19008078){\color[rgb]{1,0.6,0}\makebox(0,0)[lt]{\lineheight{0}\smash{\begin{tabular}[t]{l}Inconsistent training data\end{tabular}}}}%
    \put(0.04322525,0.2112612){\color[rgb]{0.4,0.6,0}\makebox(0,0)[lt]{\lineheight{0}\smash{\begin{tabular}[t]{l}Discretize\end{tabular}}}}%
    \put(0.87412674,0.2112612){\color[rgb]{0.8,0,0}\makebox(0,0)[lt]{\lineheight{0}\smash{\begin{tabular}[t]{l}Discretize\end{tabular}}}}%
    \put(0.12318577,0.04167682){\color[rgb]{0.4,0.6,0}\makebox(0,0)[lt]{\lineheight{0}\smash{\begin{tabular}[t]{l}DNS\end{tabular}}}}%
    \put(0.43687718,0.04135932){\color[rgb]{0.4,0.6,0}\makebox(0,0)[lt]{\lineheight{0}\smash{\begin{tabular}[t]{l}Filtered DNS\end{tabular}}}}%
    \put(0.29539931,0.14367524){\color[rgb]{0.4,0.6,0}\makebox(0,0)[lt]{\lineheight{0}\smash{\begin{tabular}[t]{l}Filter\end{tabular}}}}%
    \put(0.6530892,0.14206507){\color[rgb]{0.4,0.6,0}\makebox(0,0)[lt]{\lineheight{0}\smash{\begin{tabular}[t]{l}Close\end{tabular}}}}%
    \put(0.10883373,0.14316657){\color[rgb]{0.4,0.6,0}\makebox(0,0)[lt]{\lineheight{0}\smash{\begin{tabular}[t]{l}$D u = 0$\end{tabular}}}}%
    \put(0.05697876,0.09850932){\color[rgb]{0.4,0.6,0}\makebox(0,0)[lt]{\lineheight{0}\smash{\begin{tabular}[t]{l}$\dfrac{\mathrm{d} u}{\mathrm{d} t} = F(u) - G p$\end{tabular}}}}%
    \put(0,0){\includegraphics[width=\unitlength,page=2]{figures_equations_same.pdf}}%
    \put(0.63027657,0.00873086){\color[rgb]{1,0.6,0}\makebox(0,0)[lt]{\lineheight{0}\smash{\begin{tabular}[t]{l}Consistent training data\end{tabular}}}}%
    \put(0.79147357,0.04135932){\color[rgb]{0.2,0.4,0.8}\makebox(0,0)[lt]{\lineheight{0}\smash{\begin{tabular}[t]{l}Discrete LES\end{tabular}}}}%
    \put(0.44600225,0.32988145){\color[rgb]{0.8,0,0}\makebox(0,0)[lt]{\lineheight{0}\smash{\begin{tabular}[t]{l}$\nabla \cdot \bar{u} = 0$\end{tabular}}}}%
    \put(0.65183628,0.32991741){\color[rgb]{0.8,0,0}\makebox(0,0)[lt]{\lineheight{0}\smash{\begin{tabular}[t]{l}Close\end{tabular}}}}%
    \put(0.3853278,0.28523741){\color[rgb]{0.8,0,0}\makebox(0,0)[lt]{\lineheight{0}\smash{\begin{tabular}[t]{l}$\dfrac{\partial \bar{u}}{\partial t} = \mathcal{F}(\bar{u}) + \mathcal{C}(u) - \nabla \bar{p}$\end{tabular}}}}%
    \put(0.8001689,0.32988145){\color[rgb]{0.8,0,0}\makebox(0,0)[lt]{\lineheight{0}\smash{\begin{tabular}[t]{l}$\nabla \cdot \bar{v} = 0$\end{tabular}}}}%
    \put(0.7350847,0.28523741){\color[rgb]{0.8,0,0}\makebox(0,0)[lt]{\lineheight{0}\smash{\begin{tabular}[t]{l}$\dfrac{\partial \bar{v}}{\partial t} = \mathcal{F}(\bar{v}) + \mathcal{M}(\bar{v}) - \nabla \bar{q}$\end{tabular}}}}%
    \put(0.45452776,0.14315332){\color[rgb]{0.4,0.6,0}\makebox(0,0)[lt]{\lineheight{0}\smash{\begin{tabular}[t]{l}$\bar{D} \bar{u} = 0$\end{tabular}}}}%
    \put(0.39164846,0.09850928){\color[rgb]{0.4,0.6,0}\makebox(0,0)[lt]{\lineheight{0}\smash{\begin{tabular}[t]{l}$\dfrac{\mathrm{d} \bar{u}}{\mathrm{d} t} = \bar{F}(\bar{u}) + c(u) - \bar{G} \bar{p}$\end{tabular}}}}%
    \put(0.80682332,0.14315332){\color[rgb]{0.2,0.4,0.8}\makebox(0,0)[lt]{\lineheight{0}\smash{\begin{tabular}[t]{l}$\bar{D} \bar{v} = 0$\end{tabular}}}}%
    \put(0.73995899,0.09850928){\color[rgb]{0.2,0.4,0.8}\makebox(0,0)[lt]{\lineheight{0}\smash{\begin{tabular}[t]{l}$\dfrac{\mathrm{d} \bar{v}}{\mathrm{d} t} = \bar{F}(\bar{v}) + m(\bar{v}) - \bar{G} \bar{q}$\end{tabular}}}}%
  \end{picture}%
\endgroup%

    \caption{
        Proposed route (in green) to a discrete LES model, based on ``discretize
        first, then filter'' instead of ``filter first, then discretize'' (in
        red).
        \revboth{The term $\mathcal{F}(u)$ contains the convective and
        diffusive terms.}
    }
    \label{fig:equations}
\end{figure}

Our key insight is that the LES equations can also be obtained by ``discretizing
first" instead of ``filtering first'': following the green instead of the red
route in figure~\ref{fig:equations}. In other words, by discretizing the PDE
first, and then applying a \emph{discrete} filter, the model-data inconsistency
issue can be avoided and one can generate exact training data for the discrete
LES equations. Training data obtained by filtering discrete DNS solutions is
fully consistent with the environment where the discrete closure model is used.
The resulting LES equations do not have a coarse grid discretization error, but
an underlying fine grid discretization error and a commutator error from the
discrete filter, which can be learned using a neural network. In our recent work
\cite{Agdestein2022}, we showed the benefits of the ``discretize first''
approach on a 1D convection equation. With the discretize-first approach we
obtain stable results without the need for the stabilizing techniques mentioned
above (backscatter clipping, a-posteriori training, projection onto an
eddy-viscosity basis).

In the current paper, the goal is to extend the ``discretize first'' approach to
the full 3D incompressible Navier-Stokes equations. One major challenge that
appears in incompressible Navier-Stokes is the presence of the divergence-free
constraint. We show that discrete filters are in general not
divergence-consistent (meaning that divergence-free DNS solutions do not stay
divergence-free upon filtering).
\revone{
    \R{novel1}
    
    Kochkov et al.\ used a face-averaging filter to achieve
    divergence-consistency~\cite{Kochkov2021}. We employ this filter and show
    how it leads to a different set of equations than for
    non-divergence-consistent filters.
}
Overall, the main result is that our divergence-consistent neural
closure models lead to stable simulations.

In addition, we remark that divergence-consistent filtering is an important step
towards LES closure models that satisfy an energy inequality. Such models were
developed by us in~\cite{Vangastelen2023} for one-dimensional equations with
quadratic nonlinearity (Burgers, Korteweg - de Vries). When extending this
approach to 3D LES, the derivation of the energy inequality hinges on having a
divergence-free constraint on the filtered solution field. 

Our article is structured as follows. In section~\ref{sec:dns}, we present the
discrete DNS equations that serve as the ground truth in our problem, based on a
second order accurate finite volume discretization on a staggered grid. In
section~\ref{sec:filter}, we introduce discrete filtering. Unclosed equations
for the large scales are obtained. We show that the filtered velocity is not
automatically divergence-free. We then present
\revone{
    \R{novel2}
    
    two discrete filters on a staggered grid, of which one is
    divergence-consistent.
}
In section~\ref{sec:les}, we present our
discrete closure modeling framework, resulting in two discrete LES formulations.
We discuss the validity of the LES models in terms of divergence-consistency. A
discussion follows on the choice of closure model and how we can learn the
closure model parameters. In section~\ref{sec:results}, we present the results
of numerical experiments on a decaying turbulence test case.
Section~\ref{sec:conclusion} ends with concluding remarks.

\revboth{
    \R{appendices}
    Additional details are included in the appendices.
    In \ref{sec:discretization}, we include details about the discretization procedure.
    In \ref{sec:experiments}, we explain the numerical experiments.
    In \ref{sec:divfree}, we show that the face-averaging filter is divergence-consistent on both uniform and non-uniform grids.
    In \ref{sec:discretize_filter}, we show the problems that can occur when filtering before differentiating the divergence-free constraint.
    In \ref{sec:continuous_filtering}, we analyze the transfer functions of a continuous face-averaging filter and a continuous volume-averaging filter.
    % In \ref{sec:collocated}, we show how face-averaging can be performed on collocated grids.
    \R{commutator_ref}
    In \ref{sec:commutator}, we analyze the difference between ``discretize first'' and ``filter first'' for an analytical solution to the incompressible Navier-Stokes equations.
    In \ref{sec:other}, we provide results for two more test cases:
    a decaying turbulence test case in 2D, and a forced turbulence test case in 3D.
}

\section{Direct numerical simulation of all scales} \label{sec:dns}

In this section, we present the continuous Navier-Stokes equations and define a
discretization aimed at resolving all the scales of motion. The resulting
discrete equations will serve as the ground truth for learning an equation for
the large scales.

\subsection{The Navier-Stokes equations}

The incompressible Navier-Stokes equations describe conservation of mass and
conservation of momentum, which can be written as a divergence-free constraint
and an evolution equation:
\begin{align}
    \nabla \cdot u & = 0, \label{eq:mass_continuous} \\
    \frac{\partial u}{\partial t} +
    \nabla \cdot (u u^\mathsf{T}) & =
    -\nabla p +
    \nu \nabla^2 u +
    f, \label{eq:momentum_continuous}
\end{align}
where $\Omega \subset \mathbb{R}^d$ is the domain, $d \in \{2, 3\}$ is the
spatial dimension, $u = (u^1, \dots, u^d)$ is the velocity field, $p$ is the
pressure, $\nu$ is the kinematic viscosity, and $f = (f^1, \dots, f^d)$ is
the body force per unit of volume. The velocity, pressure, and body force are
functions of the spatial coordinate $x = (x^1, \dots, x^d)$ and time $t$.
For the remainder of this work, we assume that $\Omega$ is a rectangular domain
with periodic boundaries, and that $f$ is constant in time.

\subsection{Spatial discretization} \label{sec:spatial_discretization}

For the discretization scheme, we use a staggered Cartesian grid as proposed by
Harlow and Welch~\cite{Harlow1965}. Staggered grids have excellent conservation
properties~\cite{Lilly1965,Perot2011}, and in particular their exact discrete
divergence-free constraint is important for this work. Details about the
discretization can be found in \ref{sec:discretization}.

\begin{figure}
    \centering
    % \figsize
    \def\svgwidth{0.60\columnwidth}
    %% Creator: Inkscape 1.4 (e7c3feb100, 2024-10-09), www.inkscape.org
%% Accompanies image file 'figures_discretization.pdf' (pdf, eps, ps)
\begingroup%
  \makeatletter%
  \providecommand\color[2][]{%
    \errmessage{(Inkscape) Color is used for the text in Inkscape, but the package 'color.sty' is not loaded}%
    \renewcommand\color[2][]{}%
  }%
  \providecommand\transparent[1]{%
    \errmessage{(Inkscape) Transparency is used (non-zero) for the text in Inkscape, but the package 'transparent.sty' is not loaded}%
    \renewcommand\transparent[1]{}%
  }%
  \providecommand\rotatebox[2]{#2}%
  \newcommand*\fsize{\dimexpr\f@size pt\relax}%
  \newcommand*\lineheight[1]{\fontsize{\fsize}{#1\fsize}\selectfont}%
  \ifx\svgwidth\undefined%
    \setlength{\unitlength}{553.25643224bp}%
    \ifx\svgscale\undefined%
      \relax%
    \else%
      \setlength{\unitlength}{\unitlength * \real{\svgscale}}%
    \fi%
  \else%
    \setlength{\unitlength}{\svgwidth}%
  \fi%
  \global\let\svgwidth\undefined%
  \global\let\svgscale\undefined%
  \makeatother%
  \begin{picture}(1,0.71221847)%
    \lineheight{1}%
    \setlength\tabcolsep{0pt}%
    \put(0,0){\includegraphics[width=\unitlength,page=1]{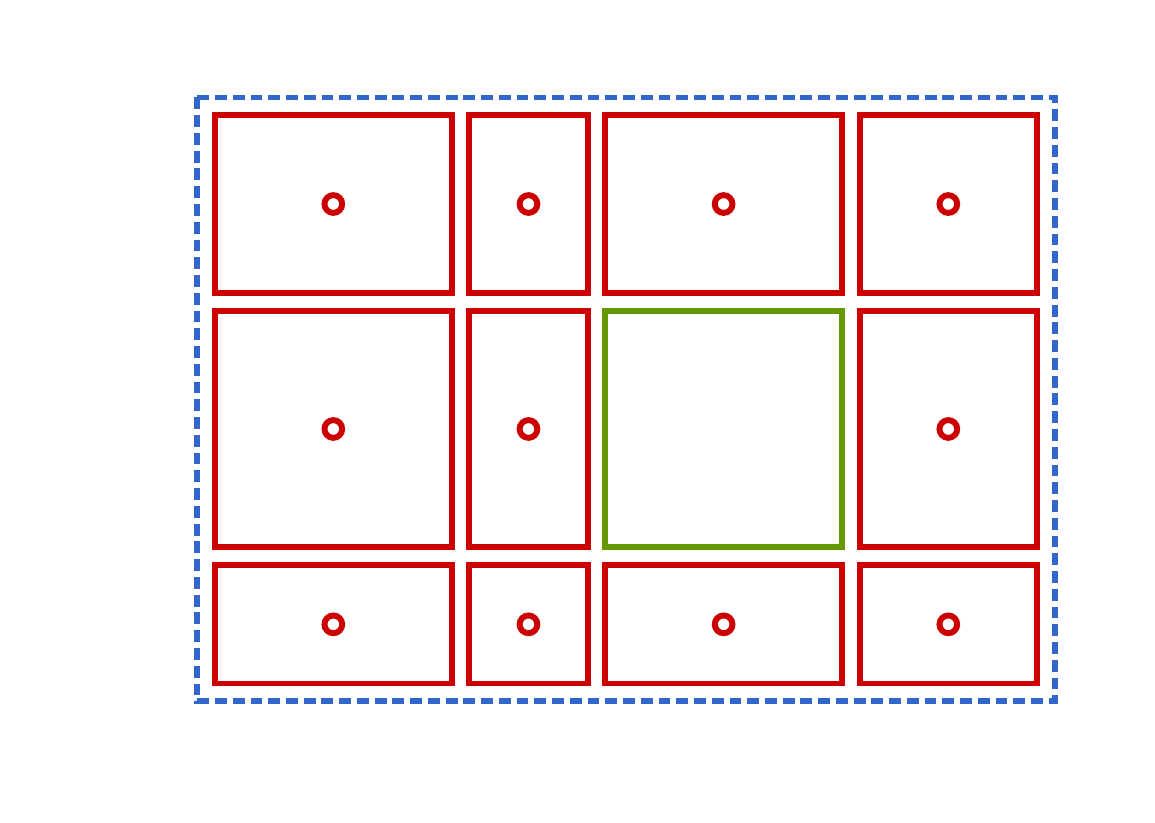}}%
    \put(0.17282437,0.6435522){\color[rgb]{0.2,0.4,0.8}\makebox(0,0)[lt]{\lineheight{1.25}\smash{\begin{tabular}[t]{l}$\Omega$\end{tabular}}}}%
    \put(0,0){\includegraphics[width=\unitlength,page=2]{figures_discretization.pdf}}%
    \put(0.23763415,0.01166483){\color[rgb]{0.2,0.4,0.8}\makebox(0,0)[lt]{\lineheight{1.25}\smash{\begin{tabular}[t]{l}$x^1$\end{tabular}}}}%
    \put(0.05299244,0.17012688){\color[rgb]{0.2,0.4,0.8}\makebox(0,0)[lt]{\lineheight{1.25}\smash{\begin{tabular}[t]{l}$x^2$\end{tabular}}}}%
    \put(0,0){\includegraphics[width=\unitlength,page=3]{figures_discretization.pdf}}%
    \put(0.58386744,0.36856025){\color[rgb]{0.4,0.6,0}\makebox(0,0)[lt]{\lineheight{1.25}\smash{\begin{tabular}[t]{l}$p$\end{tabular}}}}%
    \put(0.75402925,0.37263451){\color[rgb]{0.4,0.6,0}\makebox(0,0)[lt]{\lineheight{1.25}\smash{\begin{tabular}[t]{l}$u^1$\end{tabular}}}}%
    \put(0.64584411,0.47342431){\color[rgb]{0.4,0.6,0}\makebox(0,0)[lt]{\lineheight{1.25}\smash{\begin{tabular}[t]{l}$u^2$\end{tabular}}}}%
    \put(0,0){\includegraphics[width=\unitlength,page=4]{figures_discretization.pdf}}%
  \end{picture}%
\endgroup%

    \caption{
        Finite volume discretization on a staggered grid. The pressure is
        defined in the volume center, and the velocity components on the volume
        faces.
    }
    \label{fig:discretization}
\end{figure}

We partition the domain $\Omega$ into $N$ finite volumes.
Let $u(t) \in \mathbb{R}^{d N}$ and $p(t) \in \mathbb{R}^N$ be vectors
containing the unknown velocity and pressure components in their canonical
positions as shown in figure~\ref{fig:discretization}. They are not to be
confused with their space-continuous counterparts $u(x, t)$ and $p(x, t)$, which
will no longer be referred to in what follows. The discrete and continuous
versions of $u$ and $p$ have the same physical dimensions.

Equations \eqref{eq:mass_continuous} and \eqref{eq:momentum_continuous} are
discretized as
\begin{align}
    D u & = 0, \label{eq:mass} \\
    \frac{\mathrm{d} u}{\mathrm{d} t} & = F(u) - G p, \label{eq:momentum}
\end{align}
where
$D \in \mathbb{R}^{N \times d N}$ is the divergence operator, $G =
-\Omega_u^{-1} D^\mathsf{T} \Omega_p \in \mathbb{R}^{d N \times N}$ is the
gradient operator, $\Omega_u \in \mathbb{R}^{d N \times d N}$ and $\Omega_p
\in \mathbb{R}^{N \times N}$ are element-wise scaling operators containing the
velocity and pressure volume sizes, and $F(u) \in \mathbb{R}^{d N}$ contains
the convective, diffusive, and body force terms.

\subsection{Pressure projection}

The two vector equations \eqref{eq:mass} and \eqref{eq:momentum} form an index-2
differential-algebraic equation system~\cite{Hairer1991,Hairer2006}, consisting
of a divergence-free constraint and an evolution equation. Given $u$, the
pressure can be obtained by solving the discrete Poisson equation $0 = D F(u) -
D G p$, which is obtained by differentiating the divergence-free constraint in
time. This can also be written
\begin{equation} \label{eq:poisson}
    L p = \Omega_p D F(u),
\end{equation}
where the Laplace matrix $L = \Omega_p D G = -\Omega_p D \Omega_u^{-1} D^\mathsf{T}
\Omega_p$ is symmetric and negative semi-definite. We denote by $L^\dagger$ the
solver to the scaled pressure Poisson equation \eqref{eq:poisson}.
Since no boundary value for the pressure is prescribed,
$L$ is rank-$1$ deficient and the pressure is only
determined up to a constant. We set this constant to zero by choosing
\begin{equation}
    L^\dagger = 
    \begin{pmatrix} I & 0\end{pmatrix}
    \begin{pmatrix} L & e \\ e^\mathsf{T} & 0 \end{pmatrix}^{-1}
    \begin{pmatrix} I \\ 0 \end{pmatrix},
\end{equation}
where $e = (1, \dots, 1) \in \mathbb{R}^N$ is a vector of ones and the
additional degree of freedom enforces the constraint of an average pressure of
zero (i.e.\ $e^\mathsf{T} L^\dagger = 0$)~\cite{Luenberger1970}. In this case we
still have $L L^\dagger = I$, even though $L^\dagger L \neq I$.

We now introduce the projection operator $P$, which plays an important role in the
development of our new closure model strategy in section~\ref{sec:filter}: $P =
I - G L^\dagger \Omega_p D$. It is used to make velocity fields divergence-free
\cite{Chorin1968}, since $D P = D - \Omega_p^{-1} L L^\dagger \Omega_p D = 0$.
It naturally follows that $P$ is a projector, since $P^2 = P - G L^\dagger
\Omega_p D P = P$.

Having defined $P$, it is (at least formally) possible to eliminate the pressure
from  equations \eqref{eq:mass}-\eqref{eq:momentum} into a single
``pressure-free'' evolution equation for the velocity~\cite{Sanderse2013}, given
by
\begin{equation} \label{eq:dns}
    \frac{\mathrm{d} u}{\mathrm{d} t} = P F(u).
\end{equation}
This way, an initially divergence-free velocity field $u$ stays divergence-free
regardless of what forcing term $F$ is applied. We note that equation \eqref{eq:dns}
alone does not enforce the divergence-free constraint, as it also requires the
initial conditions to be divergence-free.

$L^\dagger$ and $P$ are non-local operators that are not explicitly
assembled, but their action on vector fields is computed on demand using an
appropriate linear solver. Formulation \eqref{eq:dns} is used as a starting
point for developing a new filtering technique in section~\ref{sec:filter}.

\subsection{Time discretization} \label{sec:RK}

\revboth{
    \R{Wray3}
    
    We use Wray's low storage third order Runge-Kutta method (Wray3) for the
    incompressible Navier-Stokes equations~\cite{Wray1990,Sanderse2013}.
}
While explicit
methods may require smaller time steps (depending on the problem-specific
trade-off of stability versus accuracy), they are easier to differentiate with
automatic differentiation tools.

Given the solution $u_n$ at a time $t_n$, the next solution at a time $t_{n +
1}$ is given by
\begin{equation}
    u_{n + 1} = u_{n} + \Delta t_{n} \sum_{i = 1}^{s} b_{i} k_i,
\end{equation}
where
\begin{equation} \label{eq:RKstep}
    k_i = P F \left(u_n + \Delta t_{n} \sum_{j = 1}^{i - 1} a_{i j} k_j \right),
\end{equation}
$\Delta t_n = t_{n + 1} - t_n$, $s$ is the number of stages, $a \in
\mathbb{R}^{s \times s}$, and $b \in \mathbb{R}^s$.
In practice, each of the RK steps \eqref{eq:RKstep} are performed by first
computing a tentative (non-divergence-free) velocity field, subsequently solving
a pressure Poisson equation, and then correcting the velocity field to be
divergence-free. As the method is explicit, this is equivalent to
\eqref{eq:RKstep} and does not introduce a splitting error~\cite{Sanderse2013}.
\revboth{
    
    For Wray3, we set
    $s = 3$,
    $b_1 = 1 / 4$,
    $b_2 = 0$,
    $b_3 = 3 / 4$,
    $a_{2 1} = 8 / 15$,
    $a_{3 1} = 1 / 4$,
    $a_{3 2} = 5 / 12$,
    and the other coefficients
    $a_{i j} = 0$.
}

\section{Discrete filtering and divergence-consistency} \label{sec:filter}

The DNS discretization presented in the previous section is in general too
expensive to simulate for problems of practical interest, and it is only used to
generate reference data for a limited number of test cases. In this section, we
present discrete filtering from the fine DNS grid to a coarse grid in order to
alleviate the computational burden. This means we filter the discretized
equations, in contrast to many existing approaches, which filter the
\textit{continuous} Navier-Stokes equations, and then apply a discretization.
The advantages of ``discretizing first'' were already mentioned in
section~\ref{sec:introduction}.
However, one disadvantage of ``discretizing first'' is
that the filtered velocity field is in general not divergence free.
\R{novel3}
\revone{
    
    Kochkov et al. pointed out that using a face-averaging filter preserves the
    divergence-free constraint for the filtered velocity~\cite{Kochkov2021}. We
    will employ this filter and compare it to a non-divergence-preserving
    volume-averaging filter.
}

\subsection{Filtering from fine to coarse grids}

We consider two computational grids: a fine grid of size $N \in \mathbb{N}^d$
and a coarse grid of size $\bar{N} \in \mathbb{N}^d$, with $\bar{N}_\alpha
\leq N_\alpha$ for all $\alpha \in \{1, \dots, d\}$. The operators $D$, $F$,
$G$, $P$, etc., are defined on the fine grid. On the coarse grid, similar
operators are denoted $\bar{D}$, $\bar{F}$, $\bar{G}$, $\bar{P}$, etc.

Consider a flow problem. We assume that the flow is fully resolved on the fine
grid, meaning that the grid spacing is at least twice as small as the smallest
significant spatial structure of the flow. The resulting fully resolved solution
$u \in \mathbb{R}^{d N}$ is referred to as the \emph{DNS solution}. In
addition, we assume that the flow is not fully resolved on the coarse grid,
meaning that the coarse grid spacing is larger than the smallest significant
spatial structure of the flow. The aim is to solve for the large scale features
of $u$ on the coarse grid. For this purpose, we construct a discrete spatial
filter $\Phi \in \mathbb{R}^{d \bar{N} \times d N}$. The resulting
filtered DNS velocity field is given by 
\begin{equation}
    \bar{u} = \Phi u \in \mathbb{R}^{d \bar{N}},
\end{equation}
and is a coarse-grid quantity. We stress that $\bar{u}$ is by definition a
consequence of the DNS. It is \textit{not} obtained by solving the Navier-Stokes
equations on the coarse grid. That is instead the goal in the next sections.

Since $\Phi$ is a coarse-graining filter, it does not generally
commute with discrete differential operators. In particular, the divergence-free
constraint is preserved for continuous convolutional filters, but this is not
automatically the case for discrete filters. We consider this property in
detail, and investigate its impact on the resulting large scale equations.

\subsection{Equation for large scales} \label{sec:filtered_dns}

When directly filtering the differential-algebraic system
\eqref{eq:mass}-\eqref{eq:momentum}, multiple challenges arise. These are detailed
in \ref{sec:discretize_filter} and summarized here:
\begin{itemize}
    \item The filter $\Phi$ works on inputs defined in the velocity points, and
        is targeted at filtering the momentum equation. It is not directly clear
        how to filter the divergence-free constraint (which is defined in the
        pressure points), and whether a second filter
        % $\Psi$
        needs to be defined for the pressure points.
    \item The momentum equation includes a pressure term. While its gradient
        can be filtered with $\Phi$, it is not clear what a filtered pressure
        $\bar{p}$ should be, or how it should appear in the filtered momentum
        equation.
\end{itemize}

\begin{figure}
    \centering
    \def\svgwidth{0.8\columnwidth}
    %% Creator: Inkscape 1.4 (e7c3feb100, 2024-10-09), www.inkscape.org
%% Accompanies image file 'figures_routes.pdf' (pdf, eps, ps)
\begingroup%
  \makeatletter%
  \providecommand\color[2][]{%
    \errmessage{(Inkscape) Color is used for the text in Inkscape, but the package 'color.sty' is not loaded}%
    \renewcommand\color[2][]{}%
  }%
  \providecommand\transparent[1]{%
    \errmessage{(Inkscape) Transparency is used (non-zero) for the text in Inkscape, but the package 'transparent.sty' is not loaded}%
    \renewcommand\transparent[1]{}%
  }%
  \providecommand\rotatebox[2]{#2}%
  \newcommand*\fsize{\dimexpr\f@size pt\relax}%
  \newcommand*\lineheight[1]{\fontsize{\fsize}{#1\fsize}\selectfont}%
  \ifx\svgwidth\undefined%
    \setlength{\unitlength}{594.20671598bp}%
    \ifx\svgscale\undefined%
      \relax%
    \else%
      \setlength{\unitlength}{\unitlength * \real{\svgscale}}%
    \fi%
  \else%
    \setlength{\unitlength}{\svgwidth}%
  \fi%
  \global\let\svgwidth\undefined%
  \global\let\svgscale\undefined%
  \makeatother%
  \begin{picture}(1,0.41673922)%
    \lineheight{1}%
    \setlength\tabcolsep{0pt}%
    \put(0,0){\includegraphics[width=\unitlength,page=1]{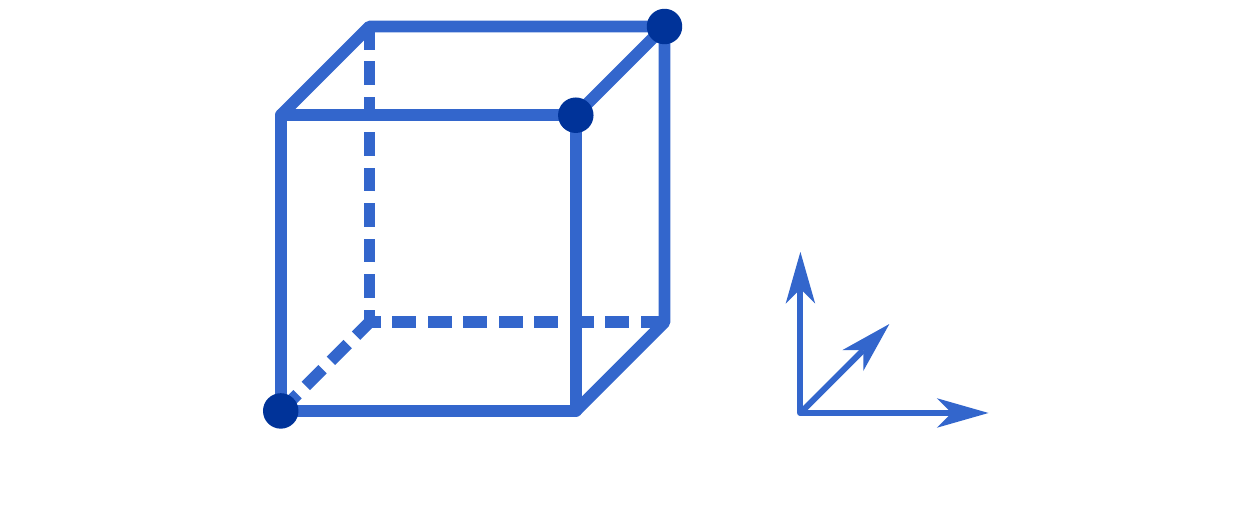}}%
    \put(0.66953191,0.05682643){\color[rgb]{0.2,0.4,0.8}\makebox(0,0)[lt]{\lineheight{0}\smash{\begin{tabular}[t]{l}Discretize\end{tabular}}}}%
    \put(0.17084815,0.04739231){\color[rgb]{0,0.2,0.6}\makebox(0,0)[lt]{\lineheight{0}\smash{\begin{tabular}[t]{l}NS\end{tabular}}}}%
    \put(0.56070415,0.38689135){\color[rgb]{0,0.2,0.6}\makebox(0,0)[lt]{\lineheight{0}\smash{\begin{tabular}[t]{l}$\mathcal{M}_\text{DIF}$\end{tabular}}}}%
    \put(0.36602449,0.34498612){\color[rgb]{0,0.2,0.6}\makebox(0,0)[lt]{\lineheight{0}\smash{\begin{tabular}[t]{l}$\mathcal{M}_\text{DCF}$\end{tabular}}}}%
    \put(0.69320688,0.17711584){\color[rgb]{0.2,0.4,0.8}\makebox(0,0)[lt]{\lineheight{0}\smash{\begin{tabular}[t]{l}Differentiate constraint\end{tabular}}}}%
    \put(0.62145655,0.21882569){\color[rgb]{0.2,0.4,0.8}\makebox(0,0)[lt]{\lineheight{0}\smash{\begin{tabular}[t]{l}Filter\end{tabular}}}}%
    \put(0,0){\includegraphics[width=\unitlength,page=2]{figures_routes.pdf}}%
    \put(0.00582049,0.20757775){\color[rgb]{0.8,0,0}\makebox(0,0)[lt]{\lineheight{0}\smash{\begin{tabular}[t]{l}Filter-discretize\end{tabular}}}}%
    \put(0,0){\includegraphics[width=\unitlength,page=3]{figures_routes.pdf}}%
    \put(0.29945066,0.12759467){\color[rgb]{0.4,0.6,0}\makebox(0,0)[lt]{\lineheight{0}\smash{\begin{tabular}[t]{l}Discretize-\\filter\end{tabular}}}}%
    \put(0,0){\includegraphics[width=\unitlength,page=4]{figures_routes.pdf}}%
    \put(0.57406344,0.31173739){\color[rgb]{0.4,0.6,0}\makebox(0,0)[lt]{\lineheight{0}\smash{\begin{tabular}[t]{l}Discretize-differentiate-filter\end{tabular}}}}%
  \end{picture}%
\endgroup%

    \caption{
        Alternative view of figure~\ref{fig:equations} to highlight the effect
        of differentiating the constraint. The red arrows show the traditional
        route of filtering first and then discretizing. The solid green arrows
        show our proposed route of discretizing first, then differentiating the
        constraint, then filtering, and finally reintroducing a pressure term
        (if the filter is divergence-consistent). This is done to circumvent the
        pressure problems of the dashed green route (discretizing first, then
        filtering).
    }
    \label{fig:routes}
\end{figure}

To circumvent these issues, we propose to differentiate the discrete divergence
constraint first (to remove the pressure), and then apply the filter to the
pressure-free DNS equation \eqref{eq:dns}. This results in the sequence
``discretize -- differentiate constraint -- filter'', as shown by solid green
arrows in figure~\ref{fig:routes}. The advantage over the route defined by
dashed green arrows, ``discretize -- filter'', is that we do not need to
consider the pressure or the divergence-free constraint, and a single filter for
the velocity field is sufficient. The ``implied'' filtered pressure will be
discussed in section~\ref{sec:discrete_les}. The resulting equation for the
filtered DNS-velocity $\bar{u}$ is $\frac{\mathrm{d} \bar{u}}{\mathrm{d} t} =
\Phi P F(u)$, which is rewritten as
\begin{equation} \label{eq:filtered_dns}
    \frac{\mathrm{d} \bar{u}}{\mathrm{d} t} = \bar{P} \bar{F}(\bar{u}) + c(u),
\end{equation}
with the unclosed commutator error defined by 
\begin{equation}\label{eq:commutator_error}
    c(u) = \Phi P F(u) - \bar{P} \bar{F}(\Phi u).
\end{equation}

A crucial point is that when filtering from a fine grid to a coarse grid, one
generally does \textit{not} get a divergence-free filtered velocity field
($\bar{D} \bar{u} \neq 0$), as discrete filtering and discrete differentiation
do not generally commute. 

\subsection{Divergence-consistent discrete filter} \label{sec:newfilter}

\begin{figure}
    \centering
    \def\svgwidth{0.48\columnwidth}
    %% Creator: Inkscape 1.4 (e7c3feb100, 2024-10-09), www.inkscape.org
%% Accompanies image file 'figures_filter.pdf' (pdf, eps, ps)
\begingroup%
  \makeatletter%
  \providecommand\color[2][]{%
    \errmessage{(Inkscape) Color is used for the text in Inkscape, but the package 'color.sty' is not loaded}%
    \renewcommand\color[2][]{}%
  }%
  \providecommand\transparent[1]{%
    \errmessage{(Inkscape) Transparency is used (non-zero) for the text in Inkscape, but the package 'transparent.sty' is not loaded}%
    \renewcommand\transparent[1]{}%
  }%
  \providecommand\rotatebox[2]{#2}%
  \newcommand*\fsize{\dimexpr\f@size pt\relax}%
  \newcommand*\lineheight[1]{\fontsize{\fsize}{#1\fsize}\selectfont}%
  \ifx\svgwidth\undefined%
    \setlength{\unitlength}{801.92029097bp}%
    \ifx\svgscale\undefined%
      \relax%
    \else%
      \setlength{\unitlength}{\unitlength * \real{\svgscale}}%
    \fi%
  \else%
    \setlength{\unitlength}{\svgwidth}%
  \fi%
  \global\let\svgwidth\undefined%
  \global\let\svgscale\undefined%
  \makeatother%
  \begin{picture}(1,0.76747034)%
    \lineheight{1}%
    \setlength\tabcolsep{0pt}%
    \put(0,0){\includegraphics[width=\unitlength,page=1]{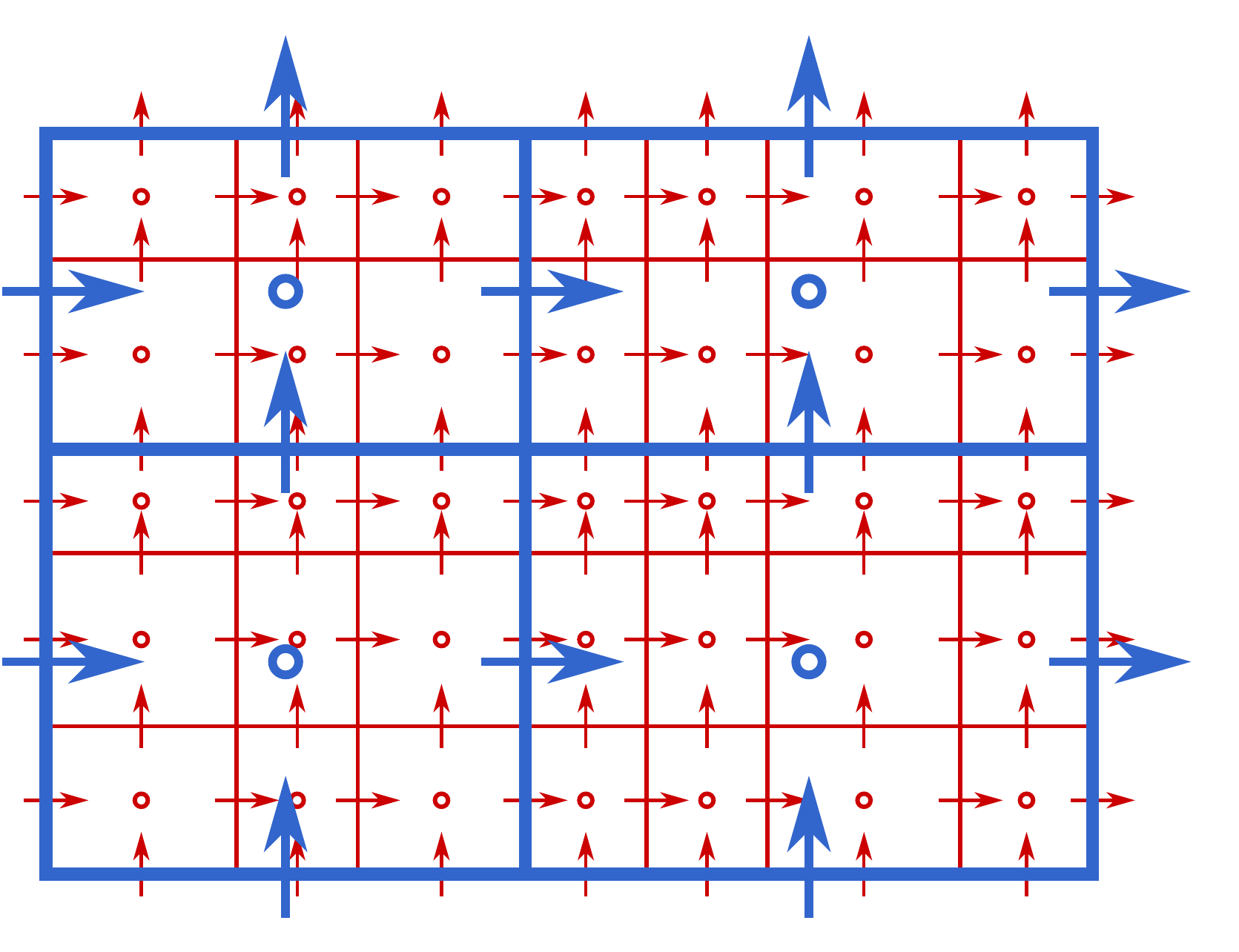}}%
  \end{picture}%
\endgroup%
    \def\svgwidth{0.48\columnwidth}
    %% Creator: Inkscape 1.4 (e7c3feb100, 2024-10-09), www.inkscape.org
%% Accompanies image file 'figures_filter_nonstruc.pdf' (pdf, eps, ps)
\begingroup%
  \makeatletter%
  \providecommand\color[2][]{%
    \errmessage{(Inkscape) Color is used for the text in Inkscape, but the package 'color.sty' is not loaded}%
    \renewcommand\color[2][]{}%
  }%
  \providecommand\transparent[1]{%
    \errmessage{(Inkscape) Transparency is used (non-zero) for the text in Inkscape, but the package 'transparent.sty' is not loaded}%
    \renewcommand\transparent[1]{}%
  }%
  \providecommand\rotatebox[2]{#2}%
  \newcommand*\fsize{\dimexpr\f@size pt\relax}%
  \newcommand*\lineheight[1]{\fontsize{\fsize}{#1\fsize}\selectfont}%
  \ifx\svgwidth\undefined%
    \setlength{\unitlength}{801.92029097bp}%
    \ifx\svgscale\undefined%
      \relax%
    \else%
      \setlength{\unitlength}{\unitlength * \real{\svgscale}}%
    \fi%
  \else%
    \setlength{\unitlength}{\svgwidth}%
  \fi%
  \global\let\svgwidth\undefined%
  \global\let\svgscale\undefined%
  \makeatother%
  \begin{picture}(1,0.76747034)%
    \lineheight{1}%
    \setlength\tabcolsep{0pt}%
    \put(0,0){\includegraphics[width=\unitlength,page=1]{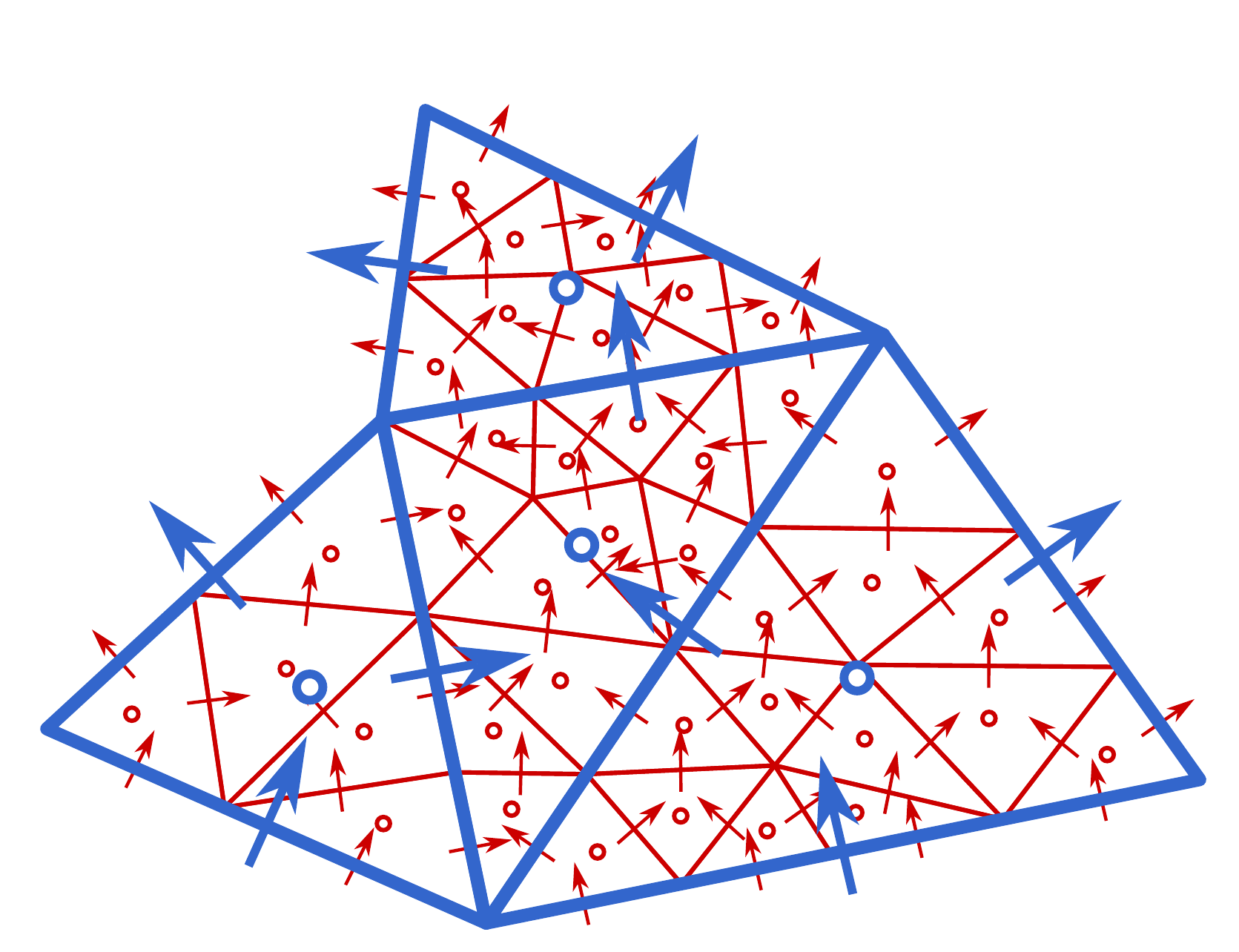}}%
  \end{picture}%
\endgroup%

    \caption{Four coarse volumes (blue) and their fine grid sub-grid volumes
        (red) in 2D. For each of the coarse volume faces, the discrete filter
        $\Phi^\text{FA}$ combines the DNS velocities $u$ into one LES velocity
        $\bar{u}$ using averaging. The interior sub-grid velocities are not
        present in $\bar{u}$. The coarse grid pressure $\bar{p}$ is defined in
        the coarse volume centers, but is not obtained by filtering $p$.
        Instead, it is computed from $\bar{u}$.
        \textbf{Left:} Structured grid, used in this work.
        \textbf{Right:} Unstructured grid.
    }
    \label{fig:coarsegrid}
\end{figure}

\R{novel4} 
Our approach is to design the filter and coarse grid such
that the divergence-free constraint is preserved. This is achieved by merging
fine DNS volumes to form coarse LES volumes, such that the faces of DNS and LES
volumes overlap (as shown in figure~\ref{fig:coarsegrid}). The extracted large
scale velocities $\bar{u} = \Phi u$ are then obtained by averaging the DNS
velocities that are found on the LES volume \emph{faces}
\revone{\cite{Kochkov2021}}. We denote this
face-averaging discrete filter by $\Phi^\text{FA}$. This approach to filtering
also naturally generalizes to unstructured grids, as long as the coarse volume
faces overlap with the fine ones. An example is shown to the right in
figure~\ref{fig:coarsegrid} for triangular volumes.

The face-averaging filter $\Phi^\text{FA}$ is different to more traditional
volume-averaging filters such as the volume-averaging top-hat filter
\cite{Sagaut2005} (that we denote $\Phi^\text{VA}$). The face-averaging filter
can be thought of as a top-hat filter acting on the dimensions orthogonal to the
velocity components only, while the volume-averaging filter is averaging over
all dimensions.
\revboth{The associated transfer function of these two types of filters are further analyzed in \ref{sec:continuous_filtering}.}
In this work, we define the two discrete filters
$\Phi^\text{FA}$ and $\Phi^\text{VA}$ with uniform weights and with filter width
equal to the coarse grid spacing $\bar{\Delta}$. Since both filters are top-hat
like, the filter width is defined as the diameter of the averaging domain in the
infinity norm. Note that due to the normalization, all discrete velocity
components $\bar{u}^\alpha_J$ and $u^\alpha_I$ share the same dimension as the
continuous velocity $u^\alpha(x, t)$ (not velocity times area or velocity times
volume).

\begin{figure}
    \centering
    \def\svgwidth{0.4\columnwidth}
    %% Creator: Inkscape 1.4 (e7c3feb100, 2024-10-09), www.inkscape.org
%% Accompanies image file 'figures_volumeaverage_arrows.pdf' (pdf, eps, ps)
\begingroup%
  \makeatletter%
  \providecommand\color[2][]{%
    \errmessage{(Inkscape) Color is used for the text in Inkscape, but the package 'color.sty' is not loaded}%
    \renewcommand\color[2][]{}%
  }%
  \providecommand\transparent[1]{%
    \errmessage{(Inkscape) Transparency is used (non-zero) for the text in Inkscape, but the package 'transparent.sty' is not loaded}%
    \renewcommand\transparent[1]{}%
  }%
  \providecommand\rotatebox[2]{#2}%
  \newcommand*\fsize{\dimexpr\f@size pt\relax}%
  \newcommand*\lineheight[1]{\fontsize{\fsize}{#1\fsize}\selectfont}%
  \ifx\svgwidth\undefined%
    \setlength{\unitlength}{476.22047244bp}%
    \ifx\svgscale\undefined%
      \relax%
    \else%
      \setlength{\unitlength}{\unitlength * \real{\svgscale}}%
    \fi%
  \else%
    \setlength{\unitlength}{\svgwidth}%
  \fi%
  \global\let\svgwidth\undefined%
  \global\let\svgscale\undefined%
  \makeatother%
  \begin{picture}(1,0.65476195)%
    \lineheight{1}%
    \setlength\tabcolsep{0pt}%
    \put(0,0){\includegraphics[width=\unitlength,page=1]{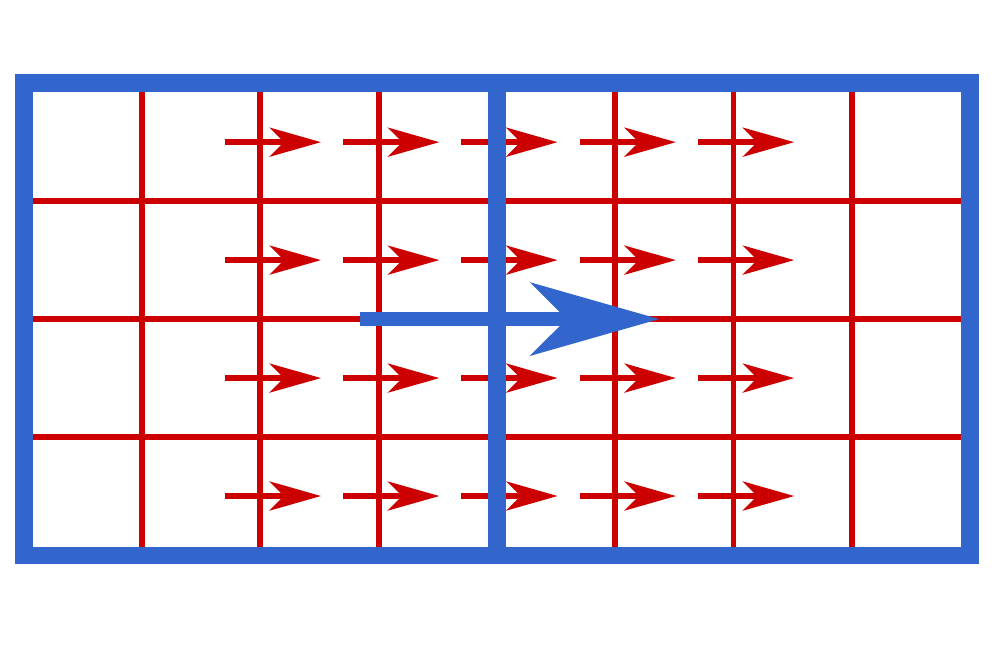}}%
    \put(0.52600472,0.02883271){\color[rgb]{0.2,0.4,0.8}\makebox(0,0)[lt]{\lineheight{0}\smash{\begin{tabular}[t]{l}$\bar{u}$\end{tabular}}}}%
    \put(0.27006762,0.60879581){\color[rgb]{0.2,0.4,0.8}\makebox(0,0)[lt]{\lineheight{0}\smash{\begin{tabular}[t]{l}Volume averaging\end{tabular}}}}%
    \put(0.4382977,0.02888992){\color[rgb]{0.8,0,0}\makebox(0,0)[lt]{\lineheight{0}\smash{\begin{tabular}[t]{l}$u$\end{tabular}}}}%
  \end{picture}%
\endgroup%
    \def\svgwidth{0.4\columnwidth}
    %% Creator: Inkscape 1.4 (e7c3feb100, 2024-10-09), www.inkscape.org
%% Accompanies image file 'figures_faceaverage_arrows.pdf' (pdf, eps, ps)
\begingroup%
  \makeatletter%
  \providecommand\color[2][]{%
    \errmessage{(Inkscape) Color is used for the text in Inkscape, but the package 'color.sty' is not loaded}%
    \renewcommand\color[2][]{}%
  }%
  \providecommand\transparent[1]{%
    \errmessage{(Inkscape) Transparency is used (non-zero) for the text in Inkscape, but the package 'transparent.sty' is not loaded}%
    \renewcommand\transparent[1]{}%
  }%
  \providecommand\rotatebox[2]{#2}%
  \newcommand*\fsize{\dimexpr\f@size pt\relax}%
  \newcommand*\lineheight[1]{\fontsize{\fsize}{#1\fsize}\selectfont}%
  \ifx\svgwidth\undefined%
    \setlength{\unitlength}{476.22047244bp}%
    \ifx\svgscale\undefined%
      \relax%
    \else%
      \setlength{\unitlength}{\unitlength * \real{\svgscale}}%
    \fi%
  \else%
    \setlength{\unitlength}{\svgwidth}%
  \fi%
  \global\let\svgwidth\undefined%
  \global\let\svgscale\undefined%
  \makeatother%
  \begin{picture}(1,0.65476195)%
    \lineheight{1}%
    \setlength\tabcolsep{0pt}%
    \put(0,0){\includegraphics[width=\unitlength,page=1]{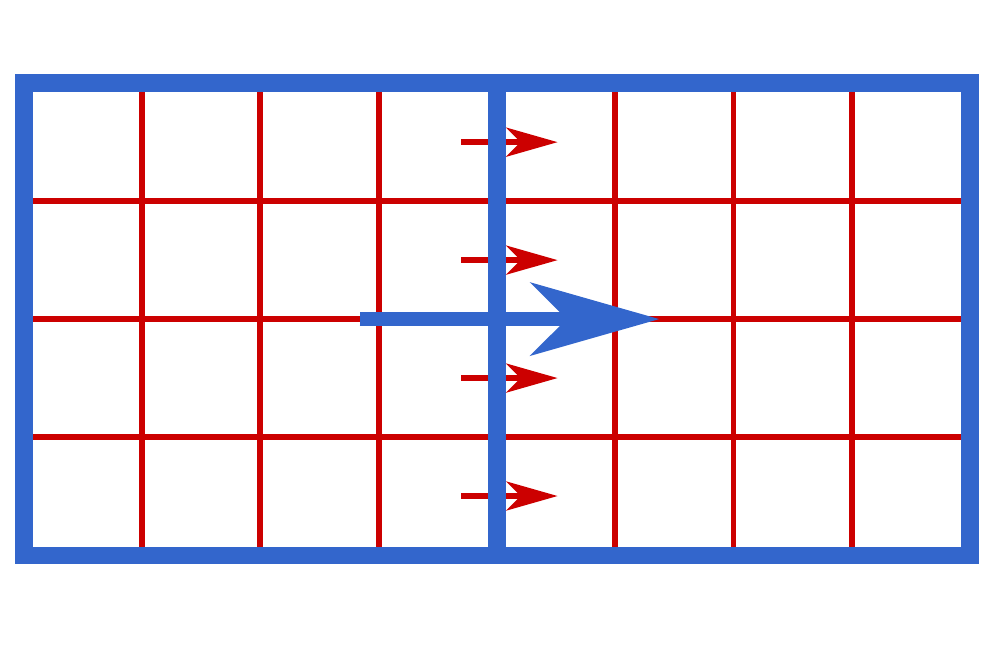}}%
    \put(0.52600472,0.02883271){\color[rgb]{0.2,0.4,0.8}\makebox(0,0)[lt]{\lineheight{0}\smash{\begin{tabular}[t]{l}$\bar{u}$\end{tabular}}}}%
    \put(0.31101506,0.60879581){\color[rgb]{0.2,0.4,0.8}\makebox(0,0)[lt]{\lineheight{0}\smash{\begin{tabular}[t]{l}Face averaging\end{tabular}}}}%
    \put(0.43273055,0.02888992){\color[rgb]{0.8,0,0}\makebox(0,0)[lt]{\lineheight{0}\smash{\begin{tabular}[t]{l}$u$\end{tabular}}}}%
  \end{picture}%
\endgroup%

    \caption{
        DNS velocity components $u$ contributing to a single filtered DNS
        component $\bar{u}$ for two filters. Both filters have a filter width
        equal to the grid size. Both $u(t)$ and $\bar{u}(t)$ share the same
        dimension as the continuous velocity $u(x, t)$.
        \textbf{Left:} Volume-averaging filter $\Phi^\text{VA}$.
        \textbf{Right:} Face-averaging filter $\Phi^\text{FA}$.
    }
    \label{fig:twofilters}
\end{figure}
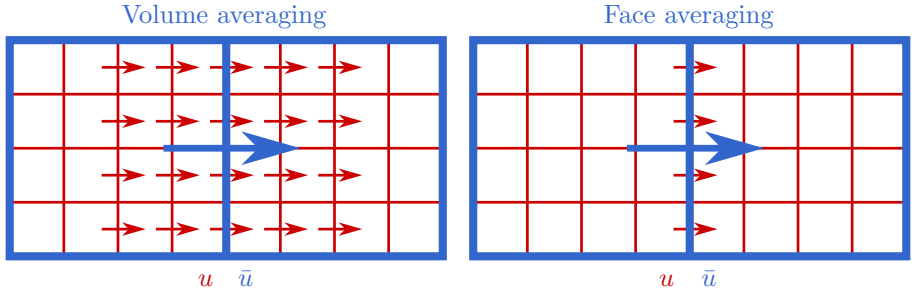

The two discrete filter supports are compared in figure~\ref{fig:twofilters}.
\R{FA-BC}
One advantage of the face-averaging filter is that it does not require 
modifications at non-periodic boundaries, while the volume-averaging filter does
(for example by using a volume of half the size and twice as large weights to
avoid averaging outside solid walls). But the main advantage lies in the
preservation of the divergence-free constraint, as we now show.

When $\bar{u}$ is obtained through face-averaging, the difference of fluxes
entering and leaving an LES volume is equal to a telescoping sum of all sub-grid
flux differences, which in turn is zero due to the fine grid divergence-free
constraint. The proof is shown in \ref{sec:divfree}. Note that this property no
longer holds if the \revboth{filter} weights are non-uniform or if the
\revboth{filter} width is different from the coarse grid spacing.
\R{novel5} \revone{ The} face-averaging filter
thus preserves \emph{by construction} the divergence-free constraint for the
filtered DNS velocity:
\begin{equation} \label{eq:ubar_divfree}
    \text{For} \ \Phi^\text{FA}: \ \boxed{\forall u, \ D u = 0 \implies \bar{D}
    \bar{u} = 0.}
\end{equation}
This property is the primary motivation for our filter choice (it can also be
thought of as a discrete equivalent of the divergence theorem $\nabla \cdot u =
0 \implies \int_{\partial \mathcal{O}} u \cdot n \, \mathrm{d} \Gamma = 0$ for
all $\mathcal{O}$). It can be used to show that the face-averaging filter has
the property
\begin{equation} \label{eq:DPhiP}
    \bar{D} \Phi P = 0
\end{equation}
since, for all $u$, if $w = P u$, then $D w = 0$ by definition
of $P$, and thus $\bar{D} \Phi P u = \bar{D} \Phi w = \bar{D} \bar{w} = 0$ from
\eqref{eq:ubar_divfree}.
In other words, divergence-free fine grid velocity fields stay divergence-free
upon filtering. For a general volume-averaging filter, we would not be able to
guarantee that all the sub-grid fluxes cancel out, and we would not get a
divergence-free constraint for $\bar{u}$. This constraint is often enforced
anyways, possibly introducing unforeseen errors, as pointed out by Sirignano et
al.~\cite{Sirignano2020}. With the face-averaging filter, we do not need to
worry about such errors.

With \revone{ the} face-averaging filter choice, the momentum
commutator error is divergence-free, since
\begin{equation}
    \bar{D} c = \bar{D} \left( \frac{\mathrm{d} \bar{u}}{\mathrm{d} t} - \bar{P}
    \bar{F}(\bar{u}) \right) = \frac{\mathrm{d} (\bar{D} \bar{u})}{\mathrm{d} t}
    - (\bar{D} \bar{P}) \bar{F}(\bar{u}) = 0
\end{equation}
since $\bar{D} \bar{u} = 0$ and $\bar{D} \bar{P} = 0$ on the coarse grid just
like $D P = 0$ on the fine grid.
As a result, $c = \bar{P} c$, and the right hand side of the large scale
equation is divergence-free. This allows us to rewrite equation
\eqref{eq:filtered_dns} as
\begin{equation} \label{eq:filtered_dns_divfree}
    \boxed{\frac{\mathrm{d} \bar{u}}{\mathrm{d} t} =
    \bar{P} \left( \bar{F}(\bar{u}) + c(u) \right).}
\end{equation}
The fact that the projection operator $\bar{P}$ now also acts on the commutator
error will play an important role in learning a new LES closure model in
section~\ref{sec:les}.

Since $D u(0) = 0$ and thus $\bar{D} \bar{u}(0) = 0$ (for the face-averaging
filter), we can  rewrite the filtered equation into an equivalent constrained
form, similar to the unfiltered equations \eqref{eq:mass}-\eqref{eq:momentum}:
\begin{align}
    \bar{D} \bar{u} & = 0, \\
    \frac{\mathrm{d} \bar{u}}{\mathrm{d} t}
        & = \bar{F}(\bar{u}) + c(u) - \bar{G} \bar{p},
\end{align}
where $\bar{p}$ is the ``implied'' pressure defined in the coarse volume centers.
It is obtained by solving the pressure Poisson equation with the additional
sub-grid forcing term $c$ in the right hand side:
\begin{equation}
    \bar{L} \bar{p} = \bar{\Omega}_p \bar{D} \left( \bar{F}(\bar{u}) + c(u) \right).
\end{equation}
Note that the coarse pressure $\bar{p}$ is not obtained defining a pressure
filter, but arises from enforcing the coarse grid divergence-free constraint. In
other words, by filtering the pressure-free momentum equation \eqref{eq:dns}
with a divergence-consistent filter, we discover what the (implicitly defined)
filtered pressure is. \revboth{An implicit volume-averaging pressure filter can
still be defined, see \ref{sec:divfree} for further details.}

\subsection{Other divergence-consistent filters}

We comment here on other approaches for constructing divergence-consistent
discrete filters.

\subsubsection{Discrete differential filters}

Continuous filters can be built using differential operators, for example
Germano's filter $\bar{u}(x, t) = (1 - \bar{\Delta}^2 / 24 \nabla^2)^{-1} u$
\cite{Germano1986,Vidal2016}. Differential filters can also be extended to the
discrete case. Trias and Verstappen propose using polynomials of the discrete
diffusion operator (which we will denote $D_2$) as a filter~\cite{Trias2011}:
\begin{equation}
    \Phi = I + \sum_{i = 1}^m \gamma_i D_2^i.
\end{equation}
where $m$ is the polynomial degree and $\gamma_i$ are filter coefficients. This
ensures that the filter has useful properties. However, divergence freeness is
only respected approximately, the convective operator may need to be modified to
preserve skew-symmetry, and there is no coarsening ($\bar{N} = N$).

\subsubsection{Spectral cut-off filters}

Spectral cut-off filters are divergence-consistent, but only so with respect to
the spectral divergence operator $\hat{u}_k \mapsto 2 \pi \mathrm{i}
k^\mathsf{T} \hat{u}_k$ (which acts element-wise in spectral space). For
pseudo-spectral discretizations, spectral cut-off filters are therefore natural
choices. On our staggered grid however, a spectral cut-off filter would not
automatically be such that $\bar{D} \bar{u} = 0$.

\subsubsection{Projected filters}

By including the projection operator into the filter definition, any filter can
be made divergence-consistent~\cite{Trias2011}. For example, the volume
averaging filter $\Phi^\text{VA}$ can be replaced with $\tilde{\Phi}^\text{VA} =
\bar{P} \Phi^\text{VA}$, which is a divergence-consistent filter. However, this
makes the filter non-local. The projection step is also more expensive, which
could be a problem when the filter is used to generate many filtered DNS
training data samples for a neural closure model.

\section{Learning a closure model for the large scale equation} \label{sec:les}

In this section, we present our new closure model formulation: a discrete LES
model based on the divergence-consistent \R{novel6} \revone{ face-averaging filter}.

\subsection{Discrete large eddy simulation} \label{sec:discrete_les}

Our ``discretize-differentiate-filter'' framework has led to equation
\eqref{eq:filtered_dns}, which describes the exact evolution of the large scale
components $\bar{u}$ for a general filter, but still contains the unclosed term
$c(u)$ from equation \eqref{eq:commutator_error}. 
% We define the goal of LES to predict $\bar{u}(t)$ from $\bar{u}(0)$, without computing $u(t)$.
Solving this equation would require access to the underlying DNS solution $u$.
We therefore replace $c(u)$ with a parameterized closure model
$m(\bar{u}, \theta) \approx c(u)$, which depends on $\bar{u}$ only
\cite{Pope2000,Sagaut2005,Berselli2006}. This produces a new \emph{approximate}
large scale velocity field $\bar{v} \approx \bar{u}$. It is defined as the
solution to the
% general (``$\text{DIF}$'')
discrete LES model
\begin{equation} \label{eq:les}
    \mathcal{M}_{\text{DIF}}: \quad
    \boxed{
    \frac{\mathrm{d} \bar{v}}{\mathrm{d} t} =
        \bar{P} \bar{F}(\bar{v}) + m(\bar{v}, \theta)
    }.
\end{equation}
Given the discretize-differentiate-filter framework, this constitutes a general LES
model formulation, which does not assume yet that the filter is
divergence-consistent. We therefore give it the label DIF (divergence-inconsistent formulation).

Since our face-averaging filter is divergence-consistent, we propose an
alternative LES model, by replacing $c$ with $m$ in equation
\eqref{eq:filtered_dns_divfree} instead of in equation \eqref{eq:filtered_dns}.
The result is a new divergence-consistent LES model: 
\begin{equation} \label{eq:les_divfree}
    \frac{\mathrm{d} \bar{v}}{\mathrm{d} t} =
        \bar{P} \left( \bar{F}(\bar{v}) + m(\bar{v}, \theta) \right).
\end{equation}
This equation is in ``pressure-free'' form, which was obtained by differentiating
the constraint. By reversing the process, i.e.\ integrating the constraint in
time, the model can now be written back into a constrained form in which the
pressure reappears:
\begin{equation} \label{eqn:LES_cons}
    \mathcal{M}_{\text{DCF}}: \quad
    \boxed{
        \begin{aligned}
        \bar{D} \bar{v} & = 0,\\
        \frac{\mathrm{d} \bar{v}}{\mathrm{d} t} & =  \bar{F}(\bar{v}) +
        m(\bar{v}, \theta)  - \bar{G} \bar{q}.
        \end{aligned}
    }
\end{equation}
This is our proposed divergence-consistent formulation that will be denoted by
``$\text{DCF}$''. The derivation of \eqref{eqn:LES_cons} from \eqref{eq:les_divfree}
hinges on the fact that $\bar{D} \bar{v}(0) = \bar{D} \bar{u}(0) = 0$. The
pressure field $\bar{q} \approx \bar{p}$ is the filtered pressure field ensuring
that the LES solution $\bar{v}$ stays divergence-free. We stress that in our
approach no pressure filter needs to be defined explicitly. With a
divergence-consistent formulation the filtered pressure can be seen as a
Lagrange multiplier.

We note that the system \eqref{eqn:LES_cons} seems to have the same
\textit{form} as the LES equations that are common in literature, being obtained
by the classic route ``filter first, then discretize'' (see red route in
figure~\ref{fig:routes}).
One might argue that the divergence-consistent face-averaging
filter is just a way to make sure that the operations of differentiation and
filtering commute. However, as we mentioned in section~\ref{sec:introduction},
there is an important difference: in contrast to the classic approach, in our
approach we have precisely defined what the filter is, and the training data
($\bar{u}$) is fully discretization-consistent with our learning target
($\bar{v}$). This is a key ingredient in obtaining model-data consistency and
hence stable closure models~\cite{Duraisamy2021,Kurz2021}.

\subsection{Divergence-consistency and LES}

We now compare in more detail the properties of the two models, the
``divergence-inconsistent'' formulation $\mathcal{M}_{\text{DIF}}$ and the
``divergence-consistent'' formulation $\mathcal{M}_{\text{DCF}}$,
see table \ref{tab:compatibility}.

$\mathcal{M}_{\text{DIF}}$ is valid for a general (non-divergence-consistent) filter
and leads to a non-divergence free $\bar{v}$. In case a divergence-consistent
filter is used, the model $\mathcal{M}_{\text{DIF}}$ is still different from
$\mathcal{M}_{\text{DCF}}$, unless $\bar{P} m = m$. This is because one can have an
exact commutator error $c$ with the property $\bar{P} c = c$ (meaning $\bar{D} c
= 0$) but still learn an approximate commutator error $m$ such that $\bar{P} m
\neq m$ (meaning $\bar{D} m \neq 0$). 

If the $\mathcal{M}_{\text{DCF}}$ model is used
with a volume-averaging filter, inconsistencies appear. For example, while the
filtered DNS data is not divergence-free, the $\text{DCF}$ model would enforce the
LES solution to be (incorrectly) divergence-free.

\begin{table}[htbp]
    \centering
    \begin{tabular}{l c c c}
    \toprule
    
    Model & Filter & $\bar{D} \bar{u} = 0$ & $\bar{D} \bar{v} = 0$ \\
    
    \midrule
    
    $\mathcal{M}_{\text{DIF}}$ & VA &
    False &
    False \\
    $\mathcal{M}_{\text{DIF}}$ & FA &
    True &
    True if $\bar{D} m = 0$ \\    
    $\mathcal{M}_{\text{DCF}}$ & VA &
    False &
    True \\
    $\mathcal{M}_{\text{DCF}}$ & FA &
    True &
    True \\
    
    \bottomrule
\end{tabular}

% vim: conceallevel=0 textwidth=0

    \caption{
        Divergence compatibility chart for LES models (rows)
        and filter properties (columns).
        The last row shows our preferred combination.
    }
    \label{tab:compatibility}
\end{table}

So far, we have only discussed the two LES formulations in terms of their
divergence properties, leaving out other properties that could also be
important. In section \eqref{sec:results}, we compare the four combinations from
table \ref{tab:compatibility} for a turbulent flow test case, and discuss
whether they are good choices or not.

\subsection{Choosing the objective function} \label{sec:model_parameters}

To learn the model parameters $\theta$, we exploit having access to exact
filtered DNS data samples $\bar{u}$ and exact commutator errors $c(u)$
obtained through (explicitly) filtered DNS solutions. A straightforward and
commonly used approach~\cite{King2016,Maulik2017b,Beck2019} is then to minimize
a loss function of a-priori type, that only depends on the DNS-solution $u$.
We use the commutator error loss 
\begin{equation} \label{eq:Lprior}
    L^\text{prior}(\mathcal{B}, \theta) =
    \frac{1}{\# \mathcal{B}}
    \sum_{u \in \mathcal{B}}
    \frac{\| m(\bar{u}, \theta) - c(u) \|^2}{\| c(u) \|^2}
\end{equation}
where $\mathcal{B}$ is a batch of 
$\# \mathcal{B}$ DNS snapshots.
Note that \eqref{eq:Lprior} does not involve $\bar{v}$, so the effect of the
closure model on the LES solution is not taken into account. We therefore call
this a-priori training~\cite{Sanderse2024}. This approach makes training fast,
since only gradients of the neural network itself are required for gradient
descent. 

Alternatively, one can minimize an a-posteriori loss function, that also depends
on the LES solution $\bar{v}$. We use the trajectory loss
\begin{equation} \label{eq:Lpost}
    L^\text{post}(u_0, \theta) =
    \frac{1}{n_\text{unroll}} 
    \sum_{i = 1}^{n_\text{unroll}}
    \frac{\|\bar{v}_i - \bar{u}_i\|^2}{\|\bar{u}_i\|^2},
\end{equation}
where $\bar{u}_i = \Phi u_i$ is obtained by filtering the DNS solution, $u_{i +
1} = \operatorname{RK}_{\Delta t}(u_i)$ is obtained using one RK4 time step from
section~\ref{sec:RK}, $u_0$ are random initial conditions, and the LES solution
$v_{i + 1} = \operatorname{RK}_{\Delta t, \mathcal{M}, m, \theta}(v_i)$ is
computed using the same time stepping scheme as $u_{i + 1}$ but with LES
formulation $\mathcal{M}$, closure $m$, and parameters $\theta$, starting from
the exact initial conditions $\bar{v}_0 = \bar{u}_0$.

The parameter $n_\text{unroll}$ determines how many time steps we unroll. If we
choose $n_\text{unroll} = 1$, $L^\text{post}$ will be very similar to
$L^\text{prior}$. If we choose a large $n_\text{unroll}$, we predict long
trajectories, and the loss may be more sensitive to small changes in $\theta$
(``exploding'' gradients). List et al.\ and Melchers et al.\ argue that the
number of unrolled time steps should depend on the characteristic time scale
(Lyapunov time scale) of the problem~\cite{List2022,Melchers2023}. For chaotic
systems (including turbulent flows), $L^\text{post}$ is expected to grow fast in
time, and the number of unrolled time steps should be small. List et al.\ found
good results with $n_\text{unroll} \in [30, 60]$ for the
incompressible Navier-Stokes equations in 2D~\cite{List2022}. For the chaotic
Kuramoto-Sivashinsky equation in 1D, Melchers et al.\ found $n_\text{unroll} =
30$ to give optimal long term results, with $n_\text{unroll} = 120$ performing
poorly~\cite{Melchers2023}.
\revone{
    \R{nunroll}
    Kochkov et al.\ found $32$ steps to be ideal for unrolling
    \cite{Kochkov2021}. For our test case, we choose $n_\text{unroll} = 50$, so
    that the trajectory becomes long enough in time while keeping the number of
    steps limited. See \ref{sec:training} for further details.
}

\begin{figure}
    \centering
    \figsize
    \def\svgwidth{0.45\columnwidth}
    %% Creator: Inkscape 1.4 (e7c3feb100, 2024-10-09), www.inkscape.org
%% Accompanies image file 'figures_chain_prior.pdf' (pdf, eps, ps)
\begingroup%
  \makeatletter%
  \providecommand\color[2][]{%
    \errmessage{(Inkscape) Color is used for the text in Inkscape, but the package 'color.sty' is not loaded}%
    \renewcommand\color[2][]{}%
  }%
  \providecommand\transparent[1]{%
    \errmessage{(Inkscape) Transparency is used (non-zero) for the text in Inkscape, but the package 'transparent.sty' is not loaded}%
    \renewcommand\transparent[1]{}%
  }%
  \providecommand\rotatebox[2]{#2}%
  \newcommand*\fsize{\dimexpr\f@size pt\relax}%
  \newcommand*\lineheight[1]{\fontsize{\fsize}{#1\fsize}\selectfont}%
  \ifx\svgwidth\undefined%
    \setlength{\unitlength}{651.96850394bp}%
    \ifx\svgscale\undefined%
      \relax%
    \else%
      \setlength{\unitlength}{\unitlength * \real{\svgscale}}%
    \fi%
  \else%
    \setlength{\unitlength}{\svgwidth}%
  \fi%
  \global\let\svgwidth\undefined%
  \global\let\svgscale\undefined%
  \makeatother%
  \begin{picture}(1,0.52173913)%
    \lineheight{1}%
    \setlength\tabcolsep{0pt}%
    \put(0,0){\includegraphics[width=\unitlength,page=1]{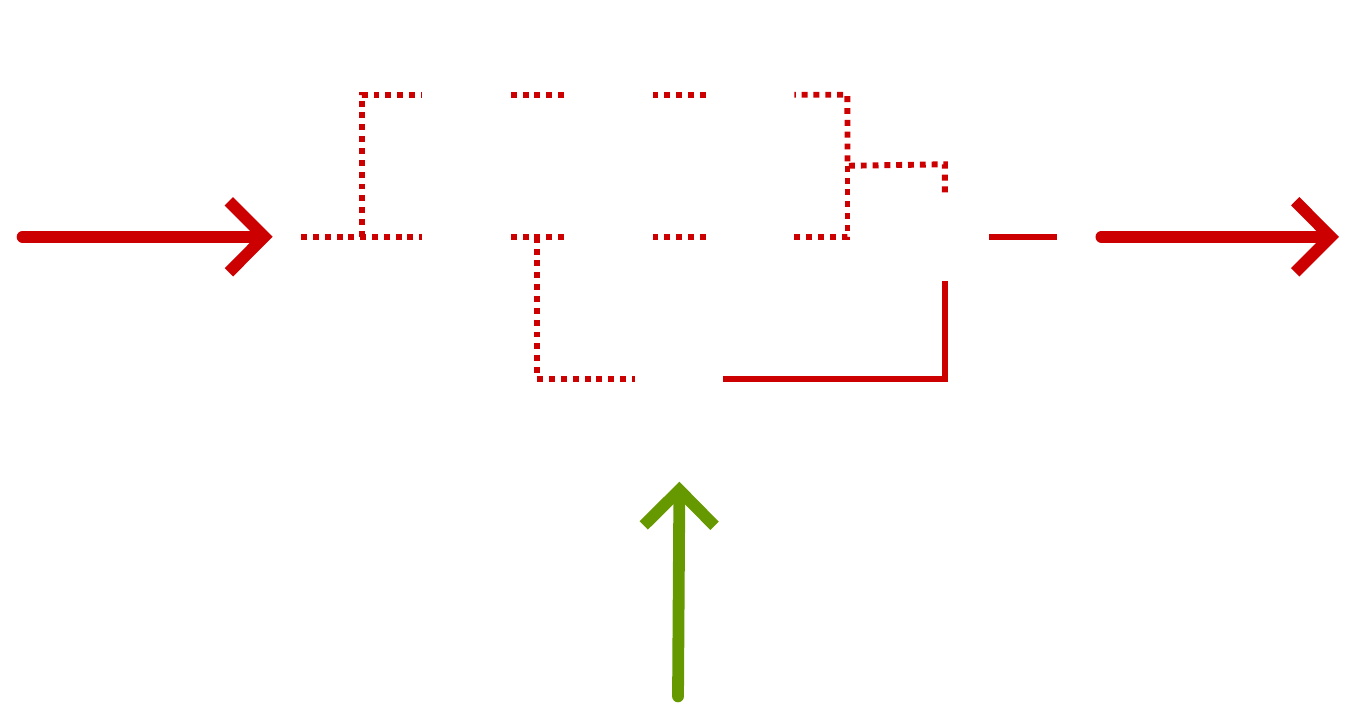}}%
    \put(0.08750607,0.37823289){\color[rgb]{0.8,0,0}\makebox(0,0)[lt]{\lineheight{1.25}\smash{\begin{tabular}[t]{l}$u$\end{tabular}}}}%
    \put(0.3439645,0.26596315){\color[rgb]{0.8,0,0}\makebox(0,0)[lt]{\lineheight{1.25}\smash{\begin{tabular}[t]{l}$\bar{u}$\end{tabular}}}}%
    \put(0.6496289,0.43802797){\color[rgb]{0.8,0,0}\makebox(0,0)[lt]{\lineheight{1.25}\smash{\begin{tabular}[t]{l}$c$\end{tabular}}}}%
    \put(0.32637669,0.43644202){\color[rgb]{0.2,0.4,0.8}\makebox(0,0)[lt]{\lineheight{1.25}\smash{\begin{tabular}[t]{l}$F$\end{tabular}}}}%
    \put(0.4776139,0.22674889){\color[rgb]{0.2,0.4,0.8}\makebox(0,0)[lt]{\lineheight{1.25}\smash{\begin{tabular}[t]{l}$m$\end{tabular}}}}%
    \put(0.67786751,0.33109671){\color[rgb]{0.2,0.4,0.8}\makebox(0,0)[lt]{\lineheight{1.25}\smash{\begin{tabular}[t]{l}$L$\end{tabular}}}}%
    \put(0.83363373,0.38370162){\color[rgb]{0.8,0,0}\makebox(0,0)[lt]{\lineheight{1.25}\smash{\begin{tabular}[t]{l}$L^\text{prior}$\end{tabular}}}}%
    \put(0.5133029,0.08065636){\color[rgb]{0.4,0.6,0}\makebox(0,0)[lt]{\lineheight{1.25}\smash{\begin{tabular}[t]{l}$\theta$\end{tabular}}}}%
    \put(0,0){\includegraphics[width=\unitlength,page=2]{figures_chain_prior.pdf}}%
    \put(0.4299065,0.43562397){\color[rgb]{0.2,0.4,0.8}\makebox(0,0)[lt]{\lineheight{1.25}\smash{\begin{tabular}[t]{l}$P$\end{tabular}}}}%
    \put(0,0){\includegraphics[width=\unitlength,page=3]{figures_chain_prior.pdf}}%
    \put(0.53425434,0.43562402){\color[rgb]{0.2,0.4,0.8}\makebox(0,0)[lt]{\lineheight{1.25}\smash{\begin{tabular}[t]{l}$\Phi$\end{tabular}}}}%
    \put(0,0){\includegraphics[width=\unitlength,page=4]{figures_chain_prior.pdf}}%
    \put(0.42990648,0.33146654){\color[rgb]{0.2,0.4,0.8}\makebox(0,0)[lt]{\lineheight{1.25}\smash{\begin{tabular}[t]{l}$\bar{F}$\end{tabular}}}}%
    \put(0,0){\includegraphics[width=\unitlength,page=5]{figures_chain_prior.pdf}}%
    \put(0.53425434,0.33127623){\color[rgb]{0.2,0.4,0.8}\makebox(0,0)[lt]{\lineheight{1.25}\smash{\begin{tabular}[t]{l}$\bar{P}$\end{tabular}}}}%
    \put(0,0){\includegraphics[width=\unitlength,page=6]{figures_chain_prior.pdf}}%
    \put(0.32555871,0.33127628){\color[rgb]{0.2,0.4,0.8}\makebox(0,0)[lt]{\lineheight{1.25}\smash{\begin{tabular}[t]{l}$\Phi$\end{tabular}}}}%
    \put(0,0){\includegraphics[width=\unitlength,page=7]{figures_chain_prior.pdf}}%
    \put(0.21580131,0.1313302){\color[rgb]{0.2,0.4,0.8}\makebox(0,0)[lt]{\lineheight{1.25}\smash{\begin{tabular}[t]{l}$L^\text{prior}$\end{tabular}}}}%
  \end{picture}%
\endgroup%
    \def\svgwidth{0.45\columnwidth}
    %% Creator: Inkscape 1.4 (e7c3feb100, 2024-10-09), www.inkscape.org
%% Accompanies image file 'figures_chain_post.pdf' (pdf, eps, ps)
\begingroup%
  \makeatletter%
  \providecommand\color[2][]{%
    \errmessage{(Inkscape) Color is used for the text in Inkscape, but the package 'color.sty' is not loaded}%
    \renewcommand\color[2][]{}%
  }%
  \providecommand\transparent[1]{%
    \errmessage{(Inkscape) Transparency is used (non-zero) for the text in Inkscape, but the package 'transparent.sty' is not loaded}%
    \renewcommand\transparent[1]{}%
  }%
  \providecommand\rotatebox[2]{#2}%
  \newcommand*\fsize{\dimexpr\f@size pt\relax}%
  \newcommand*\lineheight[1]{\fontsize{\fsize}{#1\fsize}\selectfont}%
  \ifx\svgwidth\undefined%
    \setlength{\unitlength}{651.96850394bp}%
    \ifx\svgscale\undefined%
      \relax%
    \else%
      \setlength{\unitlength}{\unitlength * \real{\svgscale}}%
    \fi%
  \else%
    \setlength{\unitlength}{\svgwidth}%
  \fi%
  \global\let\svgwidth\undefined%
  \global\let\svgscale\undefined%
  \makeatother%
  \begin{picture}(1,0.84782609)%
    \lineheight{1}%
    \setlength\tabcolsep{0pt}%
    \put(0,0){\includegraphics[width=\unitlength,page=1]{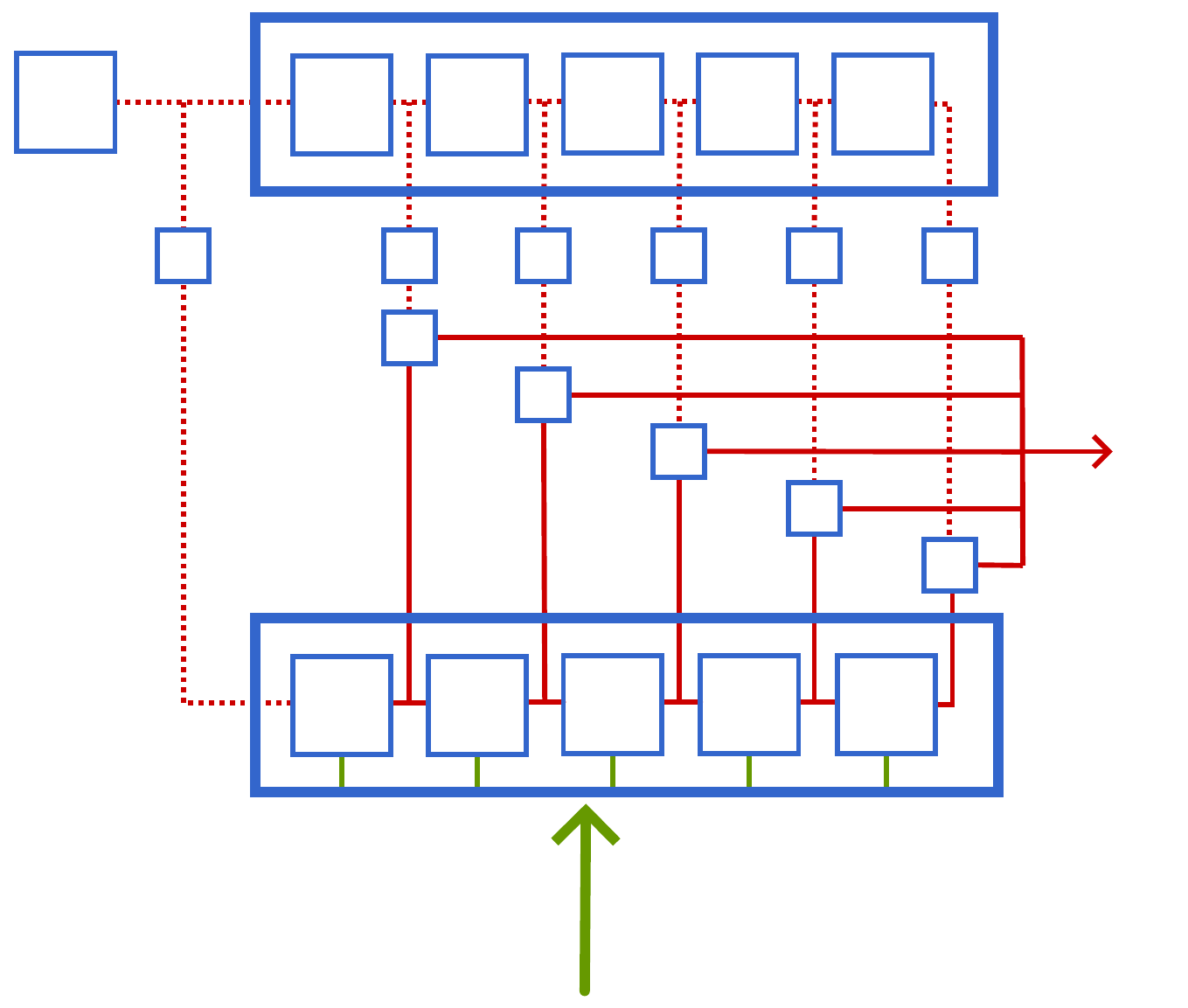}}%
    \put(0.09887656,0.22252989){\color[rgb]{0.8,0,0}\makebox(0,0)[lt]{\lineheight{1.25}\smash{\begin{tabular}[t]{l}$\bar{u}_0$\end{tabular}}}}%
    \put(0.87224823,0.48600867){\color[rgb]{0.8,0,0}\makebox(0,0)[lt]{\lineheight{1.25}\smash{\begin{tabular}[t]{l}$L^\text{post}$\end{tabular}}}}%
    \put(0.51581271,0.08278662){\color[rgb]{0.4,0.6,0}\makebox(0,0)[lt]{\lineheight{1.25}\smash{\begin{tabular}[t]{l}$\theta$\end{tabular}}}}%
    \put(0.21588562,0.12780751){\color[rgb]{0.2,0.4,0.8}\makebox(0,0)[lt]{\lineheight{1.25}\smash{\begin{tabular}[t]{l}$\operatorname{LES}$\end{tabular}}}}%
    \put(0.20898345,0.63687351){\color[rgb]{0.2,0.4,0.8}\makebox(0,0)[lt]{\lineheight{1.25}\smash{\begin{tabular}[t]{l}$\operatorname{DNS}$\end{tabular}}}}%
    \put(0.12665654,0.78000318){\color[rgb]{0.8,0,0}\makebox(0,0)[lt]{\lineheight{1.25}\smash{\begin{tabular}[t]{l}$u_0$\end{tabular}}}}%
    \put(0.13684272,0.61601612){\color[rgb]{0.2,0.4,0.8}\makebox(0,0)[lt]{\lineheight{1.25}\smash{\begin{tabular}[t]{l}$\Phi$\end{tabular}}}}%
    \put(0.32733001,0.61601607){\color[rgb]{0.2,0.4,0.8}\makebox(0,0)[lt]{\lineheight{1.25}\smash{\begin{tabular}[t]{l}$\Phi$\end{tabular}}}}%
    \put(0.32733002,0.54685255){\color[rgb]{0.2,0.4,0.8}\makebox(0,0)[lt]{\lineheight{1.25}\smash{\begin{tabular}[t]{l}$L$\end{tabular}}}}%
    \put(0.4415341,0.49903468){\color[rgb]{0.2,0.4,0.8}\makebox(0,0)[lt]{\lineheight{1.25}\smash{\begin{tabular}[t]{l}$L$\end{tabular}}}}%
    \put(0.44212729,0.61601612){\color[rgb]{0.2,0.4,0.8}\makebox(0,0)[lt]{\lineheight{1.25}\smash{\begin{tabular}[t]{l}$\Phi$\end{tabular}}}}%
    \put(0.03459953,0.74693446){\color[rgb]{0.2,0.4,0.8}\makebox(0,0)[lt]{\lineheight{1.25}\smash{\begin{tabular}[t]{l}$\mathcal{U}_0$\end{tabular}}}}%
    \put(0.55573817,0.45120603){\color[rgb]{0.2,0.4,0.8}\makebox(0,0)[lt]{\lineheight{1.25}\smash{\begin{tabular}[t]{l}$L$\end{tabular}}}}%
    \put(0.55462387,0.61601607){\color[rgb]{0.2,0.4,0.8}\makebox(0,0)[lt]{\lineheight{1.25}\smash{\begin{tabular}[t]{l}$\Phi$\end{tabular}}}}%
    \put(0.66764152,0.40337733){\color[rgb]{0.2,0.4,0.8}\makebox(0,0)[lt]{\lineheight{1.25}\smash{\begin{tabular}[t]{l}$L$\end{tabular}}}}%
    \put(0.66942113,0.61601607){\color[rgb]{0.2,0.4,0.8}\makebox(0,0)[lt]{\lineheight{1.25}\smash{\begin{tabular}[t]{l}$\Phi$\end{tabular}}}}%
    \put(0.78184561,0.35554873){\color[rgb]{0.2,0.4,0.8}\makebox(0,0)[lt]{\lineheight{1.25}\smash{\begin{tabular}[t]{l}$L$\end{tabular}}}}%
    \put(0.7819177,0.61601607){\color[rgb]{0.2,0.4,0.8}\makebox(0,0)[lt]{\lineheight{1.25}\smash{\begin{tabular}[t]{l}$\Phi$\end{tabular}}}}%
    \put(0.25194515,0.23880361){\color[rgb]{0.2,0.4,0.8}\makebox(0,0)[lt]{\lineheight{1.25}\smash{\begin{tabular}[t]{l}$\operatorname{RK}$\end{tabular}}}}%
    \put(0.36698138,0.23880361){\color[rgb]{0.2,0.4,0.8}\makebox(0,0)[lt]{\lineheight{1.25}\smash{\begin{tabular}[t]{l}$\operatorname{RK}$\end{tabular}}}}%
    \put(0.48201763,0.23880361){\color[rgb]{0.2,0.4,0.8}\makebox(0,0)[lt]{\lineheight{1.25}\smash{\begin{tabular}[t]{l}$\operatorname{RK}$\end{tabular}}}}%
    \put(0.36698139,0.74465522){\color[rgb]{0.2,0.4,0.8}\makebox(0,0)[lt]{\lineheight{1.25}\smash{\begin{tabular}[t]{l}$\operatorname{RK}$\end{tabular}}}}%
    \put(0.48103311,0.74465522){\color[rgb]{0.2,0.4,0.8}\makebox(0,0)[lt]{\lineheight{1.25}\smash{\begin{tabular}[t]{l}$\operatorname{RK}$\end{tabular}}}}%
    \put(0.59366402,0.74352514){\color[rgb]{0.2,0.4,0.8}\makebox(0,0)[lt]{\lineheight{1.25}\smash{\begin{tabular}[t]{l}$\operatorname{RK}$\end{tabular}}}}%
    \put(0.70764072,0.74352514){\color[rgb]{0.2,0.4,0.8}\makebox(0,0)[lt]{\lineheight{1.25}\smash{\begin{tabular}[t]{l}$\operatorname{RK}$\end{tabular}}}}%
    \put(0.25194516,0.74465522){\color[rgb]{0.2,0.4,0.8}\makebox(0,0)[lt]{\lineheight{1.25}\smash{\begin{tabular}[t]{l}$\operatorname{RK}$\end{tabular}}}}%
    \put(0.59741741,0.23880366){\color[rgb]{0.2,0.4,0.8}\makebox(0,0)[lt]{\lineheight{1.25}\smash{\begin{tabular}[t]{l}$\operatorname{RK}$\end{tabular}}}}%
    \put(0.7128172,0.23880366){\color[rgb]{0.2,0.4,0.8}\makebox(0,0)[lt]{\lineheight{1.25}\smash{\begin{tabular}[t]{l}$\operatorname{RK}$\end{tabular}}}}%
  \end{picture}%
\endgroup%

    \caption{Computational chain of loss functions. Solid lines are affected by
        changes in $\theta$. Dotted lines are not affected by $\theta$, and can
        be precomputed before training. \textbf{Left:} A-priori loss function
        \eqref{eq:Lprior}. The mean squared error $L$ is computed between the
        closure model $m$ and the commutator error $c$. \textbf{Right:}
        A-posteriori loss function \eqref{eq:Lpost} (here shown for five
        unrolled time steps). DNS initial conditions are sampled from the
        distribution $\mathcal{U}_0$, and filtered ($\Phi$) to produce LES
        initial conditions. After every time step, the LES solution is compared
        to the corresponding filtered DNS solution using the mean squared error
        $L$. The parameters $\theta$ are used in each LES RK time step, but not
        in the DNS time steps.}
    \label{fig:chain_loss}
\end{figure}

The a-priori loss function is easy to evaluate and easy to differentiate with
respect to $\theta$, as it does not involve solving the LES ODE given by the
model $\mathcal{M} \in \{
\mathcal{M}_{\text{DIF}},
% \mathcal{M}_{\dcf},
\mathcal{M}_{\text{DCF}}
\}$. However, minimizing $L^\text{prior}$ does not take into account the effect
of the prediction error on the LES solution error. The a-posteriori loss does
take into account this effect, but has a longer computational chain involving
the solution of the LES ODE~\cite{Um2020,MacArt2021,List2022,List2024,Kohl2023}.
This is illustrated in figure~\ref{fig:chain_loss}.

\subsection{Choosing the model architecture}

Traditionally, closure models are formulated in a continuous setting and they replace
the unclosed term $\nabla \cdot (\overline{u u} - \bar{u} \bar{u})$ by either
structural or functional models~\cite{Sagaut2005}. In recent machine learning
approaches, discrete data are inherently used for training the closure model,
and the loss function can take into account both structural and functional
elements~\cite{Guan2023}. In this work we use the common approach of
using a convolutional neural network for the closure model $m$
\cite{List2022}
(see section~\ref{sec:cnn}), and compare it to a traditional eddy-viscosity
model (see section~\ref{sec:eddyviscosity}).

\subsubsection{Eddy viscosity models} \label{sec:eddyviscosity}

Eddy viscosity models are functional models that consist of adding an additional
diffusive term
\begin{equation} \label{eq:smagorinsky_closure}
    m(\bar{u}, \theta) = \nabla \cdot (2 \nu_t \bar{S})
\end{equation}
to the continuously filtered Navier-Stokes equations, where $\bar{S} =
\frac{1}{2} \left( \nabla \bar{u} + \nabla \bar{u}^\mathsf{T} \right)$ is the
large scale strain rate tensor and $\nu_t$ is a turbulent viscosity
(parameterized by $\theta$). This term models transfer of energy from large to
unresolved scales. Note that $\bar{u}(x, t)$ and $m(\bar{u}(\cdot, t),
\theta)(x)$ are here continuous quantities that subsequently need to be
discretized to $\bar{u}(t)$ and $m(\bar{u}(t), \theta)$.

The Smagorinsky model~\cite{Smagorinsky1963,Lilly1967} predicts a local
viscosity of the form
\begin{equation}
    \nu_t = \theta^2 \bar{\Delta}^2 \sqrt{2 \operatorname{tr}(\bar{S} \bar{S})},
\end{equation}
where $\bar{\Delta}$ is the filter width and $\theta \in [0, 1]$ is the only
model parameter (the Smagorinsky coefficient). In our experiments, this
parameter is be fitted to filtered DNS data, similar to~\cite{Shankar2022,Guan2023}.
The model is discretized on the coarse grid, and $\bar{\Delta}$ is taken to be
the LES grid size.

\subsubsection{Convolutional neural networks (CNNs)} \label{sec:cnn}

Convolutional neural networks (CNNs) are commonly used in closure models when
dealing with structured data~\cite{Beck2019,Shankar2022,Guan2023,List2022}, and
since we are dealing with a structured Cartesian grid, this will be employed
here as well. A convolutional node $\operatorname{conv} : (u^n)_{n =
1}^{N_\text{chan}} \mapsto v$ in a CNN transforms $N_\text{chan}$ discrete input
channel functions into one discrete output function in the following non-linear
way:
\begin{equation}
    v_I = \sigma \left( b +
    \sum_{\substack{J \in \{-r, \dots, r\}^d \\ n \in \{1, \dots, N_\text{chan}\}}}
    K^n_J u^n_{I + J} \right),
\end{equation}
where
$\sigma$ is a non-linear activation function,
$b \in \mathbb{R}$ is a bias,
$K \in \mathbb{R}^{(2 r + 1)^d \times N_\text{chan}}$ is the kernel,
and $r$ is the kernel radius.

In a convolutional node it is typically assumed that all the input channels are
fields defined in the same grid points. Our closure model, on the other hand,
is defined on a staggered grid, with inputs and outputs in the velocity points,
whose locations differ in the different Cartesian directions, see
figure~\ref{fig:discretization}. Our CNN closure model is therefore defined as
follows:
\begin{equation}
        m_\text{CNN} = \operatorname{decollocate}
        \circ
        \begin{pmatrix} 
            \operatorname{conv} \\
            \vdots \\
            \operatorname{conv}
        \end{pmatrix}
        \circ \dots \circ
        \begin{pmatrix} 
            \operatorname{conv} \\
            \vdots \\
            \operatorname{conv}
        \end{pmatrix}
        \circ
        \operatorname{collocate}
\end{equation}
where the full vector of degrees of freedom $\bar{u}$ contains $d$
sub-vectors $\bar{u}^\alpha$ used as input channels to the CNN.
The degrees of freedom in these sub-vectors belong to their own canonical
velocity points. We therefore introduce a so-called collocation function as an
initialization layer to the CNN in order to produce quantities that are all
defined in the pressure points. The subsequent inner layers map from pressure
points to pressure points using kernels of odd diameters. Since the closure term
is required in the velocity points, a ``decollocation'' function is introduced in
the last layer to map back from pressure points to velocity points. Here, we use
a linear interpolation for both collocation and decollocation functions. It is
also possible to use ``divergence of a stress tensor'' as a decollocation
function in order to mimic the structure of the continuous commutator error.
However, our commutator error also includes discretization effects, where this
form may not be relevant. It would also require more (de)collocation functions
to produce the off-diagonal elements of the tensors (which should be the volume
corners in 2D and volume edges in 3D).

\section{Numerical experiment: \revone{forced} turbulence in a periodic box} \label{sec:results}

We consider a unit square domain $\Omega = [0, 1]^d$ with periodic boundaries.
The DNS initial conditions are sampled from a random velocity field defined
through its prescribed energy spectrum $\hat{E}_k$. Similar to
\cite{Orlandi2000,San2012,Maulik2017}, we create an initial energy profile by
multiplying a growing polynomial with a decaying exponential as
\begin{equation} \label{eq:energyprofile}
    \hat{E}_k = \frac{8\pi}{3 \kappa_p^5}
    \kappa^4
    \mathrm{e}^{-2\pi \left( \frac{\kappa}{\kappa_p} \right)^2},
\end{equation}
where $k \in \mathbb{Z}^d$ is the wavenumber, $\kappa = \| k \|$, and $\kappa_p$
is the peak wavenumber. The profile should grow for $\kappa < \kappa_p$, and the
decay should take over for $\kappa > \kappa_p$.
\revone{
    \R{bodyforce}
    To further distinguish between the different setups and prevent energy decay
    during long simulations, we add a 
    Kolmogorov-type body force to inject energy into the
    system, as in ~\cite{Chandler2013,Kochkov2021,List2024}.
    It is defined as a sinusoidal force at wavenumber $\kappa = 4$ as
    \begin{equation} \label{eq:bodyforce}
        f^\alpha(x^1, \dots, x^d) = \delta_{\alpha = 1} \sin(8 \pi x^2).
    \end{equation}
}
For more details about the
initialization procedure, see \ref{sec:ic}.

We perform all simulations in our open source package
\emph{IncompressibleNavierStokes.jl}, implemented in the Julia programming
language~\cite{Bezanson2017}. We use the \emph{KernelAbstractions.jl}
\cite{Churavy2023} framework for implementing back-end agnostic differential
operators, \emph{Lux.jl}~\cite{Pal2023} for neural networks components,
\emph{Zygote.jl}~\cite{Innes2018} for reverse mode automatic differentiation,
and \emph{Makie.jl}~\cite{Danisch2021} for visualization. All array operations
for DNS, LES, and training are performed on a CUDA-compatible GPU, using
\emph{CUDA.jl}~\cite{Besard2018,Besard2019}.

\subsection{Filtered DNS (2D and 3D)} \label{sec:prioranalysis}

Before showing results of our new divergence-consistent LES models, we perform
a-priori tests to investigate some characteristics of the DNS and filtered DNS
solutions. Note that ``a-priori'' here is used in relation to the analysis of the
results (namely before the LES model is employed), while ``a-priori'' in
section~\ref{sec:model_parameters} was related to the training procedure.

We generate two DNS trajectories $u(t)$ (one in 2D, one in 3D), starting from
the initial conditions defined above. The 2D simulation is performed with
resolution $N = (4096, 4096)$, and the 3D simulation with
\revtwo{ $N = (1024, 1024, 1024)$}.
\revtwo{
    
    For both simulations, we set $\kappa_p = 20$ and
    solve until $t_\text{end} = 1$ using adaptive time stepping.
    The Reynolds number is $\mathrm{Re} = 10^4$ in 2D
    and $\mathrm{Re} = 6000$ in 3D.
}
All array operations are performed on the GPU with double precision floating
point numbers
\revtwo{in 2D}
(to show divergence freeness), which works fine for this study even though GPUs
are not optimized for double precision, \revtwo{and single precision in 3D (to
fit all the arrays in the memory of a single H100 GPU}.

For the filter, we consider the face-averaging filter $\Phi^\text{FA}$ and the
volume-averaging filter $\Phi^\text{VA}$.
For the 2D setup, we use $\bar{N} = (\bar{n}, \bar{n})$
\revtwo{, and }
for the 3D setup, we use $\bar{N} = (\bar{n}, \bar{n}, \bar{n})$,
with \revtwo{$\bar{n} \in \{32, 64, 128, 256\}$}.

\subsubsection{Energy spectra}

\begin{figure}
    \centering
    \includegraphics[width=0.49\textwidth]{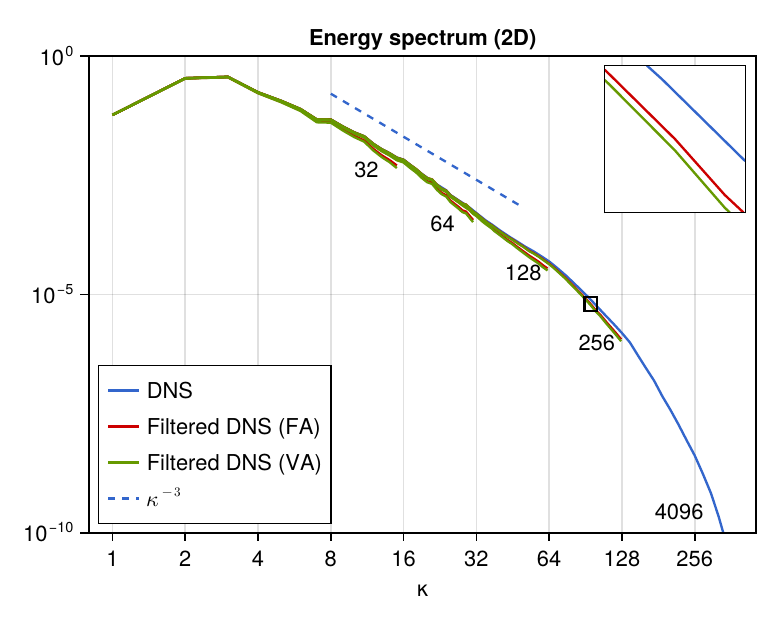}
    \includegraphics[width=0.49\textwidth]{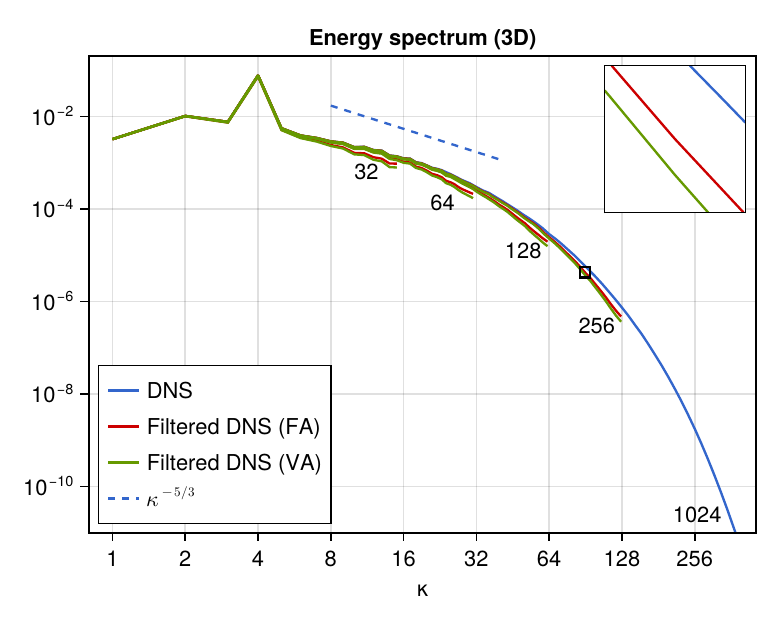}
    \caption{
        Kinetic energy spectra of DNS and filtered DNS at final time.
        The filters are both applied for \revtwo{ 4}
        different filter sizes $\bar{n}$,
        visible in the \revtwo{ four} sudden stops at the cut-off wavenumbers
        $\bar{n} / 2$.
        \revtwo{The numbers below the lines indicate the grid sizes $\bar{n}$ and $n$.}
        \textbf{Left:} 2D simulation, with theoretical scaling $\kappa^{-3}$.
        \textbf{Right:} 3D simulation, with theoretical scaling $\kappa^{-5 / 3}$.
    }
    \label{fig:priorspectra}
\end{figure}

Figure~\ref{fig:priorspectra} shows the kinetic energy spectra at the final time
for the 2D and 3D simulations. The initial velocity field is smooth (containing
only low wavenumbers), while the final DNS fields also contain higher
wavenumbers. The theoretical slopes of the inertial regions of $\kappa^{-3}$ in 2D
and $\kappa^{-5/3}$ in 3D~\cite{Pope2000} are also shown
(see section~\ref{sec:energy}). The inertial region is clearly visible in 2D,
but less so in 3D since the DNS-resolution is smaller.
The effect of diffusion is visible for
$\kappa > 128$ in the 2D plot and for $\kappa > 32$ in the 3D plot, with an
attenuation of the kinetic energy.
\revboth{The energy injection wavenumber is visible as a spike in the 3D plot,
but not in the 2D plot.}
The filtered DNS spectra are also shown. A grid of size $(n, \dots, n)$
can only fully resolve wavenumbers in the range $0 \leq \kappa \leq n / 2 - 1$,
which is visible in the sudden stops of the filtered spectra at
\revboth{
    $\bar{n} / 2 - 1$
    
}.
The face-averaging filter and the volume-averaging filter have very similar
energy profiles. Note that the filtered energy is slightly dampened even before
the filter cut-off wavelengths. Since both $\Phi^\text{FA}$ and $\Phi^\text{VA}$
are top-hat like, their transfer functions do not perform sharp cut-offs in
spectral space, but affect all wavenumbers~\cite{Pope2000,Berselli2006}.
The face-averaging filter is damping slightly less than the volume-averaging
filter. This is because it averages over one less dimension than the
volume-averaging filter, leaving the dimension normal to the face intact.
\revboth{}
Still, the coarse-graining of the discrete filters creates a spectral
cut-off effect that hides the damping of the top-hat transfer functions. This is
because the filter width is very close to the coarse-graining spectral cut-off
filter width.
\revtwo{
    \R{transferfunctions}
    For further details about the transfer functions,
    see \ref{sec:continuous_filtering}.
}

\subsubsection{Filtered fields}

\begin{figure}
    \centering
    \includegraphics[width=0.75\textwidth]{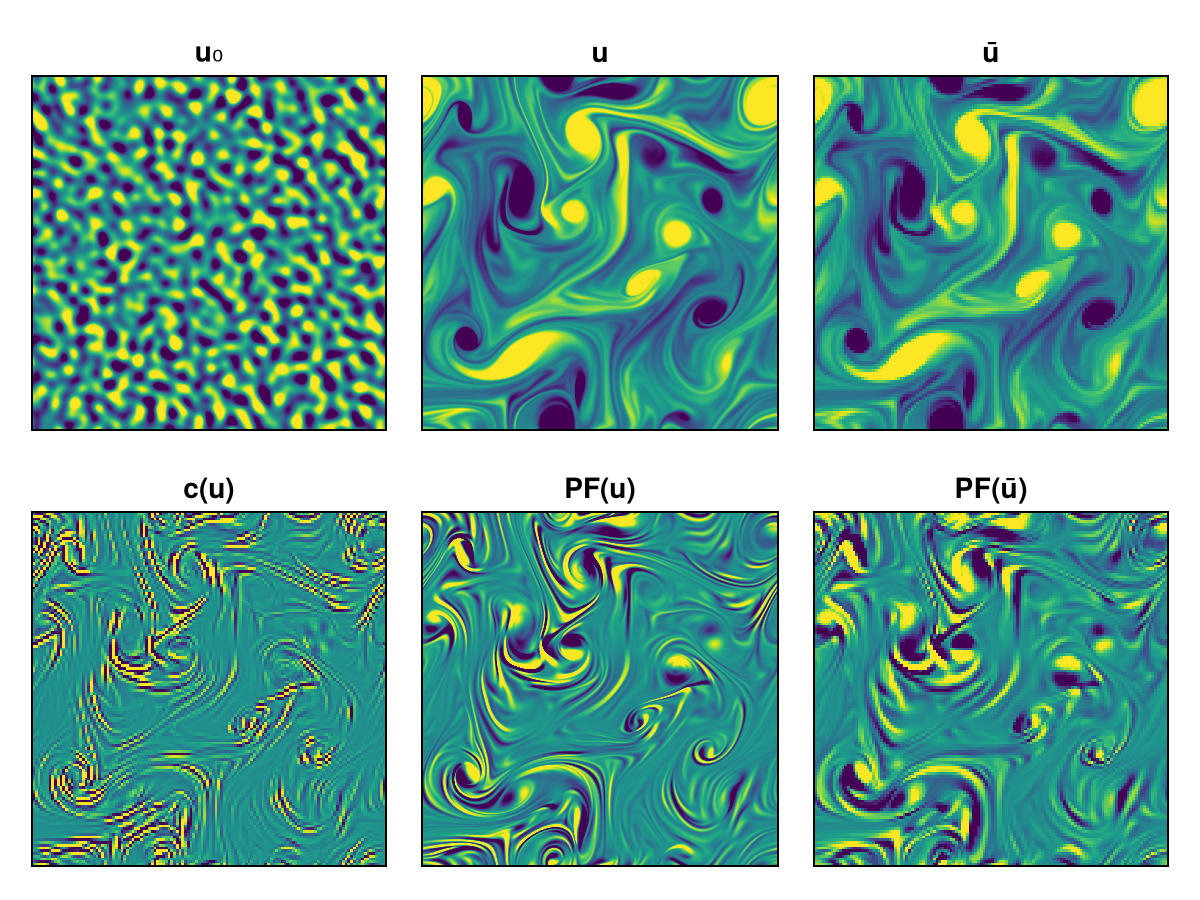}
    \caption{
        A-priori results: Discrete curl $-\delta_2 \varphi^1 + \delta_1 \varphi^2$
        of various 2D fields $\varphi \in \{ u(0), u, \bar{u}, P F(u), \bar{P}
        \bar{F}(\bar{u}), c(u)\}$.
        The filter is face-averaging, $\bar{N} = 128^2$.
    }
    \label{fig:priorfields_2D}
\end{figure}

Figure~\ref{fig:priorfields_2D} shows the discrete curl
$\nabla \times \varphi$
of various 2D fields
$\varphi \in \{ u(0), u, \bar{u}, P F(u), \bar{P} \bar{F}(\bar{u}), c(u)\}$
for the face-averaging filter with $\bar{N} = 128^2$.
We plot the curl since $\varphi$ is a vector field. Each
pixel corresponds to a pressure volume, in which the curl is interpolated for
visualization. The filtered field $\bar{u}$ in the top-right corner is clearly
unable to represent all the sub-grid fluctuations seen in the DNS field $u$, but
the large eddies of $u$ are still recognizable in $\bar{u}$. We stress again
that the aim of our neural closure models is to reproduce $\bar{u}$, without
having knowledge of the DNS field $u$. This will be shown in section
\ref{sec:analysis_les}.

Note in particular that the coarse grid right hand side $\bar{P}
\bar{F}(\bar{u})$ contains small oscillations which make the discrete curl look
grainy. These are due to the under-resolved central-difference discretization on
the coarse grid, and are subsequently also present in the discrete commutator error
$c$. The closure model $m$ thus has to predict the oscillations in $c$ in order to
correct for those in $\bar{P} \bar{F}(\bar{u})$, using information from the
smooth field $\bar{u}$ only. If $c$ was defined by filtering first and then
discretizing, as is commonly done in LES, these oscillations would not be part
of $c$ and other means of stabilization would be needed to correct for the
oscillations in $\bar{P} \bar{F}(\bar{u})$, such as
explicit LES~\cite{Lund2003,Gallagher2019a,Benjamin2023b}
(where the filter is applied to
the non-linear convective term in the LES right hand side).
This problem has also been addressed in literature; for example,
Geurts et al.~\cite{Geurts2006} and Beck and Kurz~\cite{Beck2023} argue that all commutator
errors should be modeled, and Stoffer et al.~\cite{Stoffer2021} show that
instabilities in the high wavenumbers can occur if the discretization is not
included in the commutator error.

\begin{figure}
    \centering
    \includegraphics[width=0.32\textwidth]{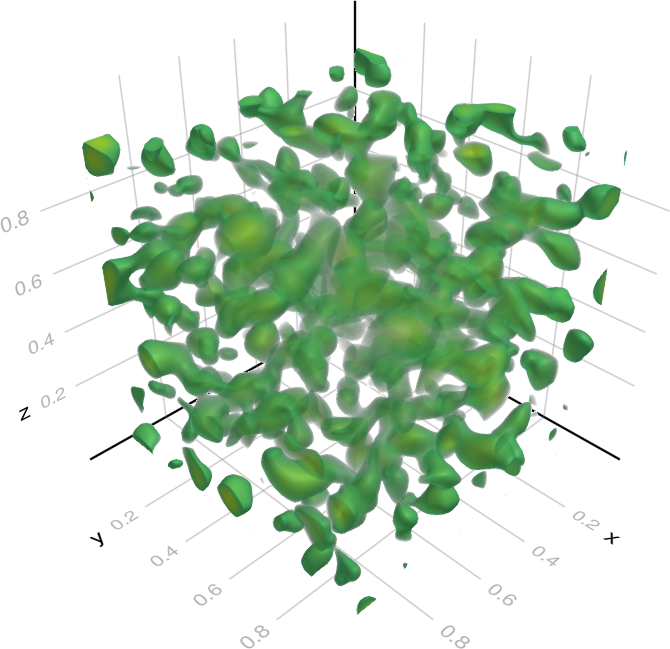}
    \includegraphics[width=0.32\textwidth]{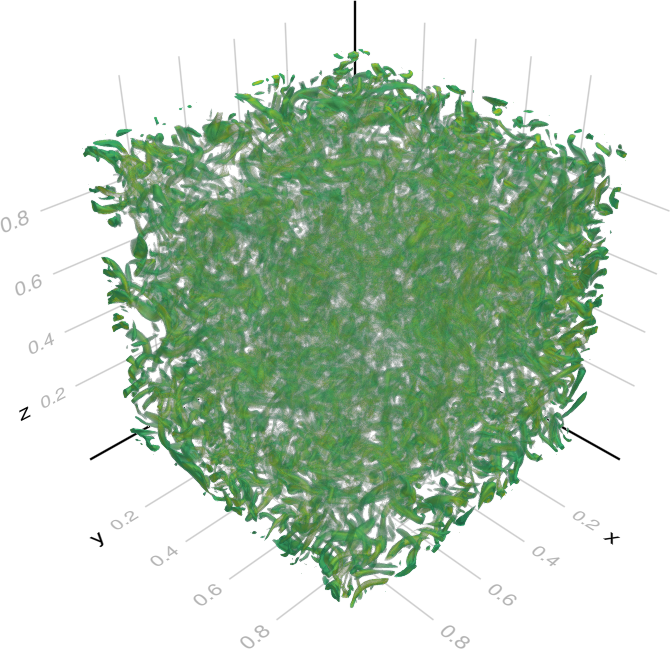}
    \includegraphics[width=0.32\textwidth]{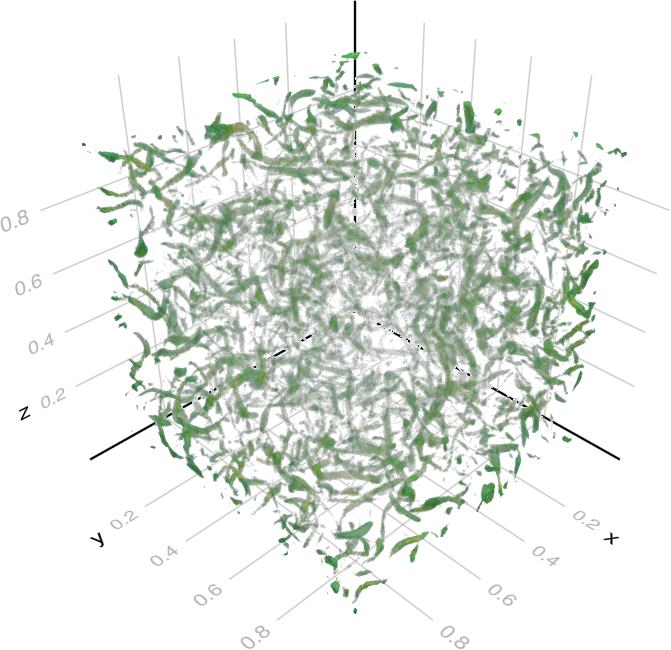}
    \caption{
        Vortex cores visualized as $10$ isocontours of negative regions of
        $\lambda_2(S^2 + T^2)$, where $\lambda_2$ denotes the second largest
        eigenvalue, and $S = \frac{1}{2} (\nabla u + \nabla u^\mathsf{T})$ and
        $T = \frac{1}{2} (\nabla u - \nabla u^\mathsf{T})$ are the symmetric and
        anti-symmetric parts of the velocity gradient. \textbf{Left:} DNS,
        initial time. \textbf{Middle:} DNS, final time. \textbf{Right:} Filtered DNS
        % (face-averaging, $\bar{N} = 64^3$),
        (face-averaging, $\bar{N} = 128^3$),
        final time.
    }
    \label{fig:priorfields_3D}
\end{figure}

Figure~\ref{fig:priorfields_3D} shows the vortex cores of the 3D simulation at
initial and final time. The vortex cores are visualized as isocontours of
$\lambda_2$-criterion~\cite{Jeong1995}. It is defined as negative regions of
$\lambda_2(S^2 + T^2)$, where $\lambda_2$ denotes the second largest eigenvalue
in absolute value of the $3 \times 3$-tensor, and
$S = \frac{1}{2} (\nabla u + \nabla u^\mathsf{T})$ and
$T = \frac{1}{2} (\nabla u - \nabla u^\mathsf{T})$ are
the symmetric and anti-symmetric parts of the velocity gradient tensor.
With the prescribed initial low wavenumber energy spectrum, only large DNS
vortex structures are present at the initial time. They are clearly visible in
the left plot. As energy gets transferred to higher wavenumbers, smaller
turbulent vortex structures are formed (middle plot). The same DNS field still
contains larger vortex structures, which become visible after filtering (right
plot).

\subsubsection{Divergence, commutator errors, and kinetic energy}

\begin{table}
    \centering
    % \revboth{
    \begin{tabular}{*{9}{c}}
    \toprule

    $N$
    & Bits
    & Filter
    & $\bar{N}$
    & $\frac{\| \bar{D} \bar{u} \|}{\| \bar{u} \|}$
    & $\frac{\| \bar{u} - \bar{P} \bar{u} \|}{\| \bar{u} \|}$
    & $\frac{\| c - \bar{P} c \|}{\| c \|}$
    & $\frac{\| c \|}{\| \bar{P} \bar{F} + c \|}$
    & $\frac{\| \bar{u} \|^2_{\bar{\Omega}}}{\| u \|^2_\Omega}$ \\

    \midrule

    % Re = 10^4, t = 1.0, double precision
    \multirow{8}{*}{$4096^2$} &
    \multirow{8}{*}{$64$} &
        FA & $\p32^2$ & $1.5\e{-14}$ & $2.5\e{-16}$ & $2.3\e{-13}$ & $0.56\p$ & $0.92$  \\
    & & FA & $\p64^2$ & $2.1\e{-14}$ & $1.9\e{-16}$ & $3.4\e{-13}$ & $0.35\p$ & $0.98$  \\
    & & FA &  $128^2$ & $3.4\e{-14}$ & $1.6\e{-16}$ & $6.1\e{-13}$ & $0.18\p$ & $0.99$  \\
    & & FA &  $256^2$ & $5.3\e{-14}$ & $1.3\e{-16}$ & $1.3\e{-12}$ & $0.077 $ & $1.0$\p \\
    & & VA & $\p32^2$ & $1.1\p$      & $0.017\p\p$  & $0.11\p$     & $0.54$   & $0.89$  \\
    & & VA & $\p64^2$ & $0.67$       & $0.0058\p$   & $0.096$      & $0.33$   & $0.97$  \\
    & & VA &  $128^2$ & $0.39$       & $0.0019\p$   & $0.088$      & $0.18$   & $0.99$  \\
    & & VA &  $256^2$ & $0.19$       & $0.00059$    & $0.12\p$     & $0.08$   & $1.0\p$ \\
    \midrule

    % Re = 6000, t = 1.0, single precision
    \multirow{8}{*}{$1024^3$} &
    \multirow{8}{*}{$32$} &
        FA & $\p32^3$ & $2.2\e{-5}$ & $2.8\e{-7}$ & $3.9\e{-6}$ & $0.85$ & $0.63$ \\
    & & FA & $\p64^3$ & $2.2\e{-5}$ & $1.4\e{-7}$ & $3.5\e{-6}$ & $0.70$ & $0.80$ \\
    & & FA & $128^3$  & $2.2\e{-5}$ & $7.6\e{-8}$ & $4.0\e{-6}$ & $0.51$ & $0.91$ \\
    & & FA & $256^3$  & $2.4\e{-5}$ & $4.5\e{-8}$ & $7.6\e{-6}$ & $0.31$ & $0.96$ \\
    & & VA & $\p32^3$ & $3.3$       & $0.043\p$   & $0.19$      & $0.80$ & $0.60$ \\
    & & VA & $\p64^3$ & $4.2$       & $0.028\p$   & $0.15$      & $0.66$ & $0.77$ \\
    & & VA & $128^3$  & $4.8$       & $0.017\p$   & $0.13$      & $0.49$ & $0.89$ \\
    & & VA & $256^3$  & $5.0$       & $0.0094$    & $0.13$      & $0.30$ & $0.95$ \\

    \bottomrule
\end{tabular}

% vim: conceallevel=0 textwidth=0 wrap=0

    % }
    \caption{Magnitude of various quantities derived from two DNS trajectories
        $u(t)$ (one in 2D, one in 3D). All quantities are averaged over time.
        \revboth{
            The machine precision for 64-bit numbers is
            $\epsilon \approx 2.22 \times 10^{-16}$,
            and for 32-bit numbers
            $\epsilon \approx 1.19 \times 10^{-7}$.
        }
        }
    \label{tab:prioranalysis}
\end{table}

Table \ref{tab:prioranalysis} shows the magnitude of various quantities derived
from the two DNS trajectories. All quantities $q(u(t))$ are averaged over
time with a frequency of $s = 20$ time steps as follows:
\begin{equation}
    \langle q \rangle_s = \frac{1}{n_t / s + 1} \sum_{i = 0}^{n_t / s}
    q(u(t_{s i})).
\end{equation}
The considered quantities $q(u)$ are:
Normalized divergence $\frac{\| \bar{D} \bar{u} \|}{\| \bar{u} \|}$,
magnitude of non-divergence-free part of filtered velocity field
$\frac{\| \bar{u} - \bar{P} \bar{u} \|}{\| \bar{u} \|}$,
magnitude of non-divergence-free part of commutator error
$\frac{\| c - \bar{P} c \|}{\| c \|}$, 
magnitude of commutator error in the total filtered right hand side
$\frac{\| c \|}{\| \bar{P} \bar{F} + \bar{c} \|}$,
and resolved kinetic energy
$\frac{\| \bar{u} \|^2_{\bar{\Omega}}}{\| u \|^2_\Omega}$.
The norms are defined as $\| u \| = \sqrt{\sum_{\alpha = 1}^d \| u^\alpha \|^2}$
for vector fields such as $u$ and
$\| u^\alpha \| = \sqrt{\sum_{I} | u^\alpha_I |^2}$
for scalar fields such as $u^\alpha$. Additionally, the norm
$\| \cdot \|_\Omega$ is weighted by the volume sizes.

It is clear that both $\bar{u}$ and $c(u)$ are divergence-free for the
face-averaging filter, in both 2D and 3D.
For the volume-averaging filter on the other hand, $\bar{u}$ and $c(u)$ are
not divergence-free.
At \revtwo{ $\bar{N} = 32^2$},
the orthogonal projected part $\bar{u} - \bar{P} \bar{u}$ comprises about
$1.7\%$ of $\bar{u}$. This is more visible
for the commutator error, for which the orthogonal part is
\revtwo{ $11\%$}
of the total commutator error.
When the grid is refined to $\bar{N} = 256^2$ however, the
non-divergence-free parts of $\bar{u}$
\revtwo{
    
    shrink to $0.059\%$ since the flow is more resolved.
    However, the non-divergence free part of commutator error is
    $12\%$, since the commutator error itself is smaller.
}
For \revtwo{volume-averaging} with $\bar{N} = 32^3$,
the non-divergence free parts of $\bar{u}$ and $c(u)$ are
\revtwo{ $4.3\%$ and $19\%$} respectively.
\revtwo{
    
    They shrink to $0.94\%$ and $13\%$ for $\bar{N} = 256^3$.
}
Note that the volume-averaging filter
width was chosen to be equal to the grid spacing, but these divergence errors
would be larger if we increased the filter width.

For both 2D and 3D, and both filters types, the commutator error becomes smaller
when the grid is refined, which is expected since more scales are resolved. The
commutator error magnitude does not seem to depend much on whether we use
face-averaging or volume-averaging, since they both have the same characteristic
filter width. For
\revtwo{
    
    face-averagining with $\bar{N} = 32^2$ and $\bar{N} = 128^3$
}, more than half of the
total right hand side is due to the commutator error, even though 
\revtwo{ $92\%$ and $91\%$}
of the kinetic energy is resolved by $\bar{u}$. For these resolutions, it
is very important to have a good closure model. For $\bar{N} = 256^2$, the
filtered DNS right hand side is much closer to the corresponding coarse
unfiltered DNS right hand side, with $c(u)$ comprising
\revtwo{ $7.7\%$}
of the total right hand side $\bar{P} \bar{F}(\bar{u}) + c(u)$.
A closure model is still clearly needed, even though
\revtwo{ $100\%$}
of the energy being resolved by $\bar{u}$.
Bae et al.~\cite{Bae2022} also found the discretization part of the commutator
errors to be quite significant, in particular near walls (which we do not
consider in this study).

In summary, the DNS and filtered DNS results confirm the theoretical analysis in
section~\ref{sec:newfilter}.
The volume-averaging filter lacks
divergence-consistency of the solution and the closure term. The magnitude of
the closure term is similar for both filters. The benefits of
divergence-consistency will be demonstrated in the a-posteriori analysis in the
subsequent section.

\subsection{LES (2D)} \label{sec:analysis_les}

Next, we turn to the results for the key challenge set out in this paper:
testing our neural closure models in an LES setting, with divergence-consistent
filters, aiming to approximate the trajectory $\bar{u}(t)$ given $\bar{u}(0)$.
We now only consider the 2D problem to reduce the computational time, since
the same network is trained repeatedly in multiple configurations. For 3D
results, see \ref{sec:les3D}.  The DNS resolution is set to $N = (4096, 4096)$.
The Reynolds number is
\revboth{
    
    $\mathrm{Re} = 6000$
}. The initial peak wavenumber is $\kappa_p = 20$. We use
single precision floating point numbers for all computations, including DNS
trajectory generation, to reduce memory usage and increase speed.
Details about the datasets are found in \ref{sec:datasets}.

For the closure model $m$, we consider three options:
\begin{enumerate}
    \item No closure model, $m_0 = 0$. This is the baseline model.
    \item Smagorinsky closure model, $m_\text{S}$.
    \item Convolutional neural closure model, $m_\text{CNN}$.
        The architecture is shown in \ref{sec:architecture}.
\end{enumerate}
\revtwo{
    The no-closure model can be thought of as ``coarse DNS'', and it is included
    to show the necessity of a closure model when compared to the filtered DNS
    reference data.
    \R{smagorinsky-FA}
    The Smagorinsky model is a traditional closure model. While it is
    conventionally used in combination with volume-averaging filters, it can
    also be used with a face-averaging filter.
    The idea behind the Smagorinsky eddy viscosity
    $\nu_t = (\theta \bar{\Delta})^2 \sqrt{2 \operatorname{tr}(\bar{S} \bar{S})}$
    is that it is equal to the large-scale strain weighted by the
    ``average sub-filter eddy size'' $\theta \bar{\Delta}$,
    which is parameterized by the fractional constant $\theta \in [0, 1]$.
    The largest unresolved eddy then has the size $\bar{\Delta}$.
    This argument also holds for the face-averaging filter, where the largest
    unresolved eddy has a characteristic size $\bar{\Delta}$, and the averaged
    unresolved eddy has a size $\theta \bar{\Delta}$ for some coefficient
    $\theta$ (potentially with a different value than for volume-averaging).
}
See \ref{sec:training} for details about the training.

\subsubsection{A-priori errors}

\begin{figure}
    \centering
    \includegraphics[width=\textwidth]{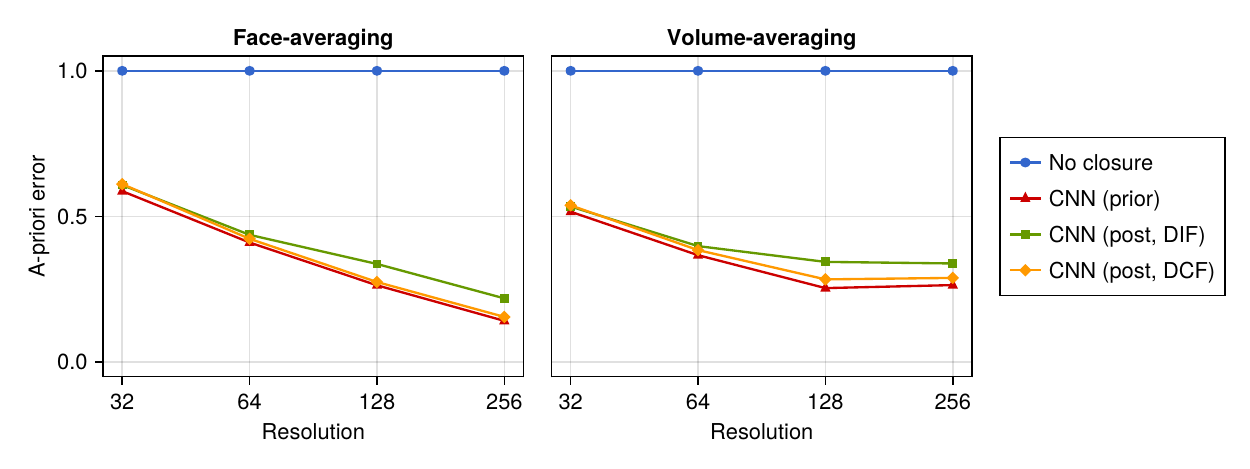}
    \caption{
        \revboth{Relative} a-priori errors
        $\frac{1}{n_\text{u}} \sum_u \frac{\| m - c \|}{\| c \|}$
        for the testing dataset.
        \textbf{Left:} Face-averaging filter.
        \textbf{Right:} Volume-averaging filter.
    }
    \label{fig:convergence_prior}
\end{figure}

Figure~\ref{fig:convergence_prior} shows the average relative a-priori errors
$\frac{1}{n_\text{u}} \sum_u \frac{\| m - c \|}{\| c \|}$
for the testing dataset and the three CNN
parameters $\theta^\text{prior}$, $\theta^\text{post}_\text{DIF}$,
and $\theta^\text{post}_\text{DCF}$. For the no-closure model $m_0$, the a-priori
error is always $\| 0 - c \| / \| c \| = 1$. For both FA and VA and for all grid
sizes, the a-priori trained parameters $\theta^\text{prior}$ perform the best.
This is expected, since the a-priori loss function used to obtain
$\theta^\text{prior}$ is similar to the a-priori error.

\subsubsection{A-posteriori errors}

\begin{figure}
    \centering
    \includegraphics[width=\textwidth]{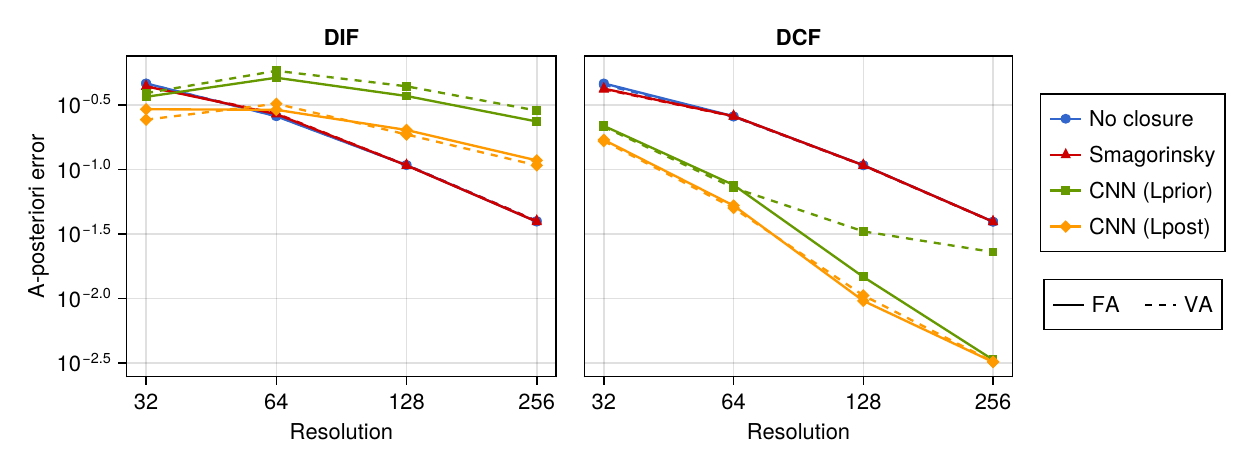}
    \caption{
        \revboth{Relative} a-posteriori errors
        $\frac{1}{n_t} \sum_{i = 1}^{n_t} \| \bar{v}_i - \bar{u}_i \| / \| \bar{u}_i \|$
        \revboth{at time $t = 0.27$}
        for the testing dataset. The CNN is
        trained separately for each resolution and each filter type with
        $L^\text{prior}$ (green squares) and
        $L^\text{prior}$-then-$L^\text{post}$ (yellow diamonds). \textbf{Solid
        lines:} Face-averaging ($\Phi^\text{FA}$). \textbf{Dashed lines:}
        Volume-averaging ($\Phi^\text{VA}$). \textbf{Left:} General model
        $\mathcal{M}_{\text{DIF}}$. \textbf{Right:} Divergence-consistent model
        $\mathcal{M}_{\text{DCF}}$.
    }
    \label{fig:convergence}
\end{figure}

The relative a-posteriori error
$\frac{1}{n_t} \sum_{i = 1}^{n_t} \| \bar{v}_i - \bar{u}_i \| / \| \bar{u}_i \|$
is computed for the trajectory in the testing dataset.
A \revboth{ closure model} parameter set $\theta$ is only used
for testing on the same coarse grid and same filter type that it was trained
for.
\revboth{
    Since the turbulent flow is a chaotic system, it does not make sense to plot
    the point-wise error over long time periods. We therefore plot the errors at
    time $t = 0.27$. For later time instances, one needs to consider statistical quantities.
}
Figure~\ref{fig:convergence} shows the average error over time.

For $\mathcal{M}_{\text{DCF}}$ (right), the no-closure and the Smagorinsky closure
have similar errors
\revboth{}.
The CNN is clearly outperforming the other two closures.
\revboth{
    
    The a-priori trained CNN is performing significantly better for
    face-averaging than for volume-averaging. Since the commutator errors
    are fully consistent with $\mathcal{M}_\text{DCF}$ in the face-averaging case,
    high accuracy can be achieved with a-priori training alone.
    In the volume-averaging case, the commutator errors are inconsistent
    with $\mathcal{M}_\text{DCF}$, leading to higher errors.
    When switching to the a-posteriori loss function however, the
    volume-averaging CNN drastically improves and catches up with the
    face-averaging one, and the errors become indistinguishable for the two
    filters.
    In the face-averaging case, the a-priori trained CNN is already performing
    well, and training with $L^\text{post}$ only leads to minor improvements.
}
Note that \revboth{ the a-posteriori} training is done starting from
the best performing a-priori trained parameters.

For $\mathcal{M}_{\text{DIF}}$ (left), the no-closure and the Smagorinsky closure have
similar profiles as for $\mathcal{M}_{\text{DCF}}$. For the no-closure, the two LES
formulations are actually identical.
\revboth{
    
    For all resolutions except for $\bar{n} = 32$,
    $m_\text{CNN}(\cdot, \theta^\text{prior})$ is showing signs
    of instability and is performing worse than $m_0$.
    The a-posteriori training does manage to detect and reduce this error, but
    it is still unstable and higher than $m_0$.
}

\subsubsection{Stability}

\begin{figure}
    \centering
    \includegraphics[width=\textwidth]{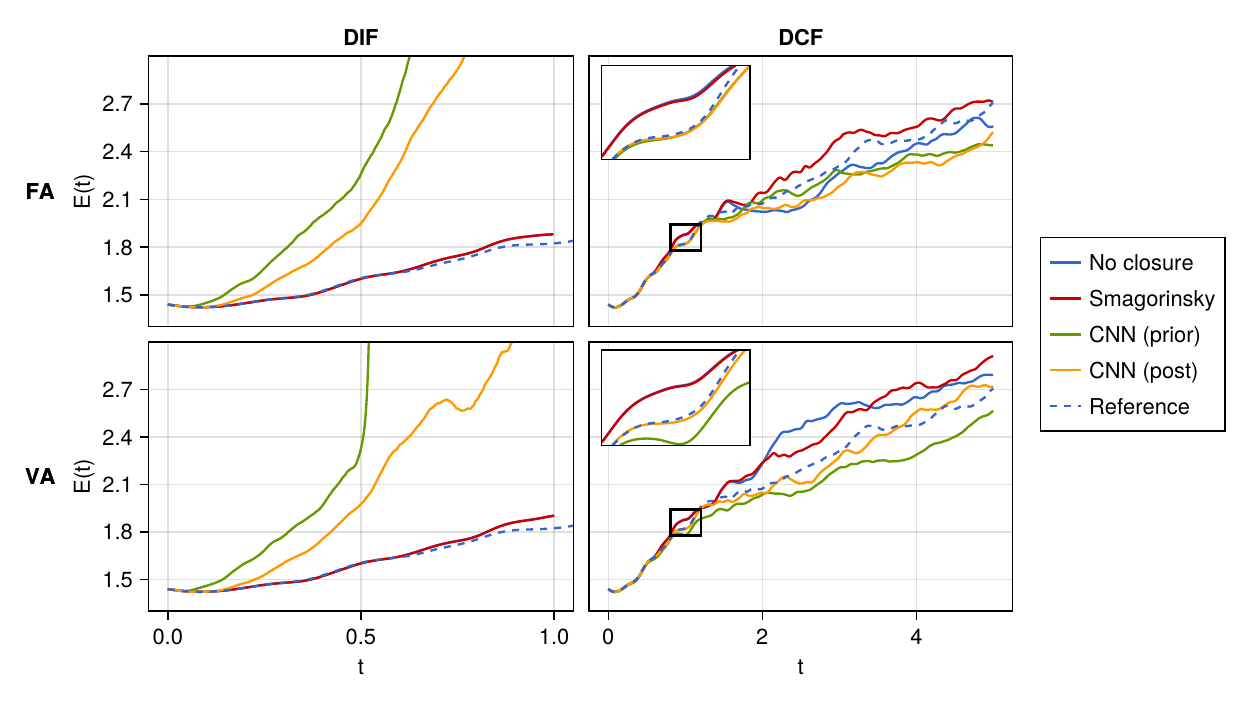}
    \caption{
        Total kinetic energy evolution for
        \revboth{ $\bar{n} = 256$}.
        \textbf{Left:} Unprojected closure model $\mathcal{M}_{\text{DIF}}$.
        \textbf{Right:} Constrained model $\mathcal{M}_{\text{DCF}}$.
        \textbf{Top:} Face-averaging filter.
        \textbf{Bottom:} Volume-averaging filter.
    }
    \label{fig:energy_evolution}
\end{figure}

To further investigate the stability, we compute the evolution of the total
kinetic energy as a function of time. This is shown in
figure~\ref{fig:energy_evolution}.
\revboth{
    The model DIF becomes unstable in the long term, so we only compute the
    energy until the time $t = 1$.
}
\revboth{
    
    The no-closure solution and Smagorinsky solution both stay close to the target energy
    during the first time unit, after which they become completely decorellated
    from the reference solution due to the chaotic nature of turbulence.
    However, their energy profile does follow the same trend as the reference,
    staying slightly above.
    For $\mathcal{M}_{\text{DCF}}$, the
    the a-posteriori trained CNN stays on the reference energy level for much longer
    than the two other models. As in figure~\ref{fig:convergence},
    $m_\text{CNN}(\cdot, \theta^\text{prior})$ and
    $m_\text{CNN}(\cdot, \theta^\text{post})$ are very similar in the face-averaging case
    (as can be seen in the zoom-in window), until they reach the point of decorrelation.
    In the volume-averaging case however, 
    $m_\text{CNN}(\cdot, \theta^\text{prior})$
    performs more poorly than 
    $m_\text{CNN}(\cdot, \theta^\text{post})$.
}
For $\mathcal{M}_{\text{DIF}}$, the high CNN errors from figure~\ref{fig:convergence} are
confirmed by the rapid growth of the total kinetic energy. Similar growth in
energy of unconstrained neural closure models after a period of seemingly good
overlap with the reference energy has been observed by Beck and Kurz
\cite{Beck2019,Kurz2021}, although in a different configuration. Training with
$L^\text{post}$ improves the stability, and the energy stays close to the
reference for a longer period. This was not sufficient to stabilize
$\mathcal{M}_{\text{DIF}}$, however, and the energy eventually starts increasing.

The growth in kinetic energy and resulting lack of stability for the
divergence-inconsistent model can be explained by the energy-conservation
properties of our spatial discretization. The convective terms are discretized
with a second order central scheme (in so-called divergence form), which can be
shown to be equivalent to a skew-symmetric, energy-conserving form
\textit{provided that the velocity field is divergence-free}~\cite{Verstappen2003}.
If the velocity field is not divergence-free, there is no
guarantee that the convective terms are energy-conserving, which can lead to
growth in kinetic energy and loss of stability. 

\revone{
    \R{energy_in_loss}
    It is interesting to point out that during the first time unit for
    $\mathcal{M}_\text{DCF}$, the CNN models manage to stay close to the target
    energy level without being explicitly trained to do so. It is also possible
    to add an energy mismatch term in the loss function, thus encouraging the
    CNN to produce a correct energy level~\cite{List2022}. Such training could
    potentially be used to further improve the stability of
    $\mathcal{M}_\text{DIF}$.
}

\subsubsection{Energy spectra}

\begin{figure}
    \centering
    \includegraphics[width=\textwidth]{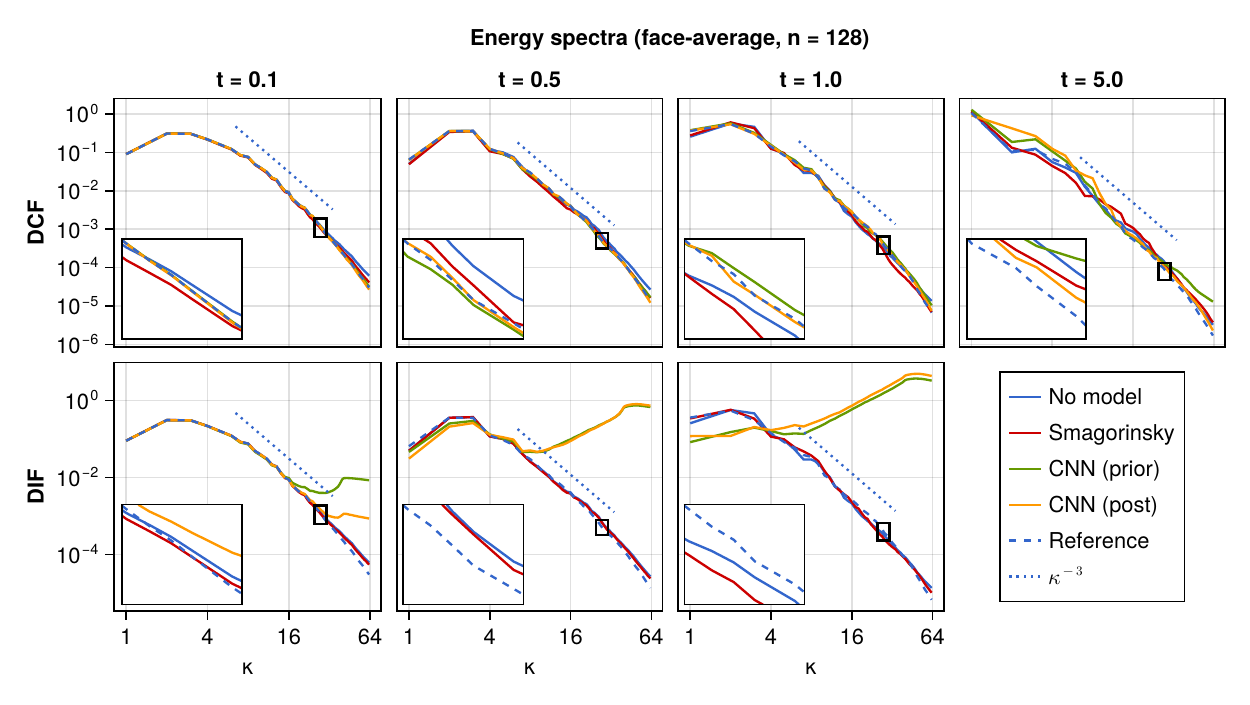}
    \includegraphics[width=\textwidth]{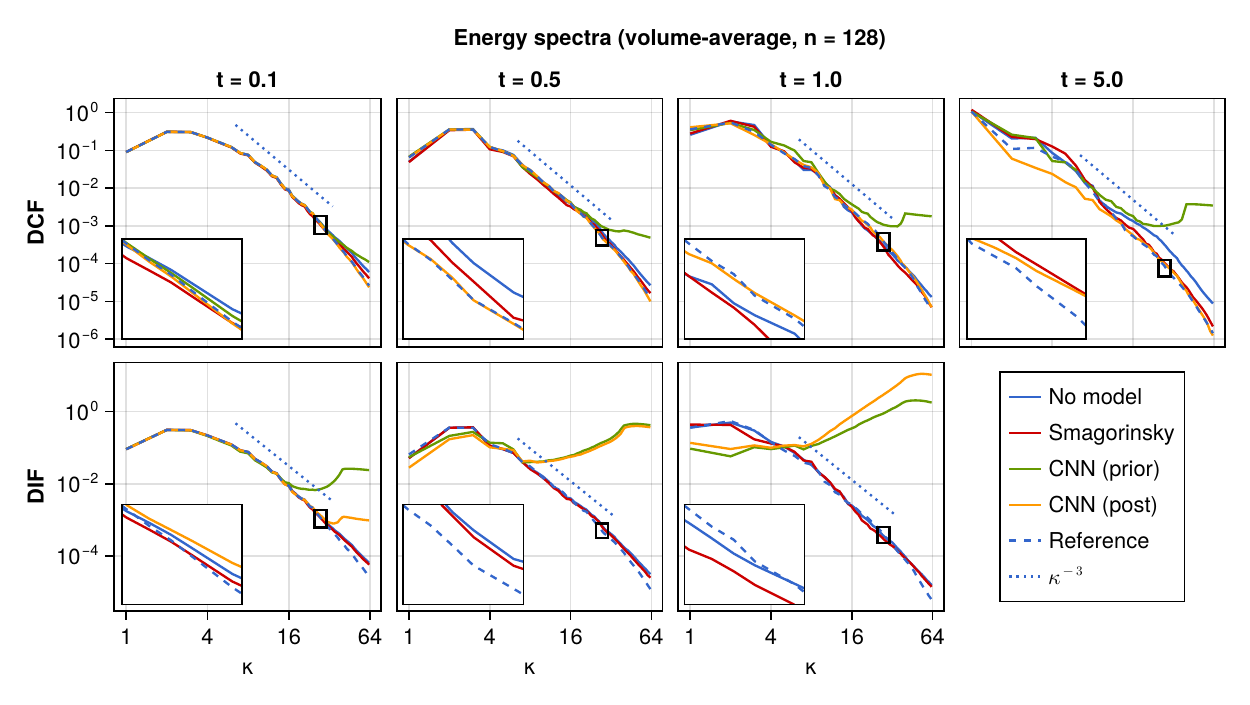}
    \caption{Energy spectra at final time for
        $\bar{n} = 128$.
        \textbf{Top:} Face-averaging filter.
        \textbf{Bottom:} Volume-averaging filter.
        \revboth{Since the $\mathcal{M}_\text{DIF}$ simulations become unstable,
        we do not show their spectrum at time $t = 5$.}
    }
    \label{fig:energy_spectra}
\end{figure}

The energy spectra at \revtwo{ different times} are shown in
figure~\ref{fig:energy_spectra}.
The no-closure model spectrum is generally close to
the reference spectrum, but contains too much energy in the high wavenumbers.
This is likely due to the oscillations discussed in the previous sections.
The Smagorinsky closure is correcting for this (as intended)
\revtwo{}.
For $\mathcal{M}_{\text{DCF}}$,
the a-posteriori trained CNN spectrum is visually very close to the reference
spectrum
\revtwo{, also at the final time, long after the LES trajectory has 
diverged from the filtered DNS trajectory. This shows that the trained model
can capture turbulence statistics correctly even though the system is chaotic.}
For $\mathcal{M}_{\text{DIF}}$, 
\revtwo{
    
    all the CNNs produce too much energy in the high wave numbers.
}
Training with $L^\text{post}$ seems to improve this somewhat.
    
\subsubsection{LES solution fields}

\begin{figure}
    \centering
    \includegraphics[width=\textwidth]{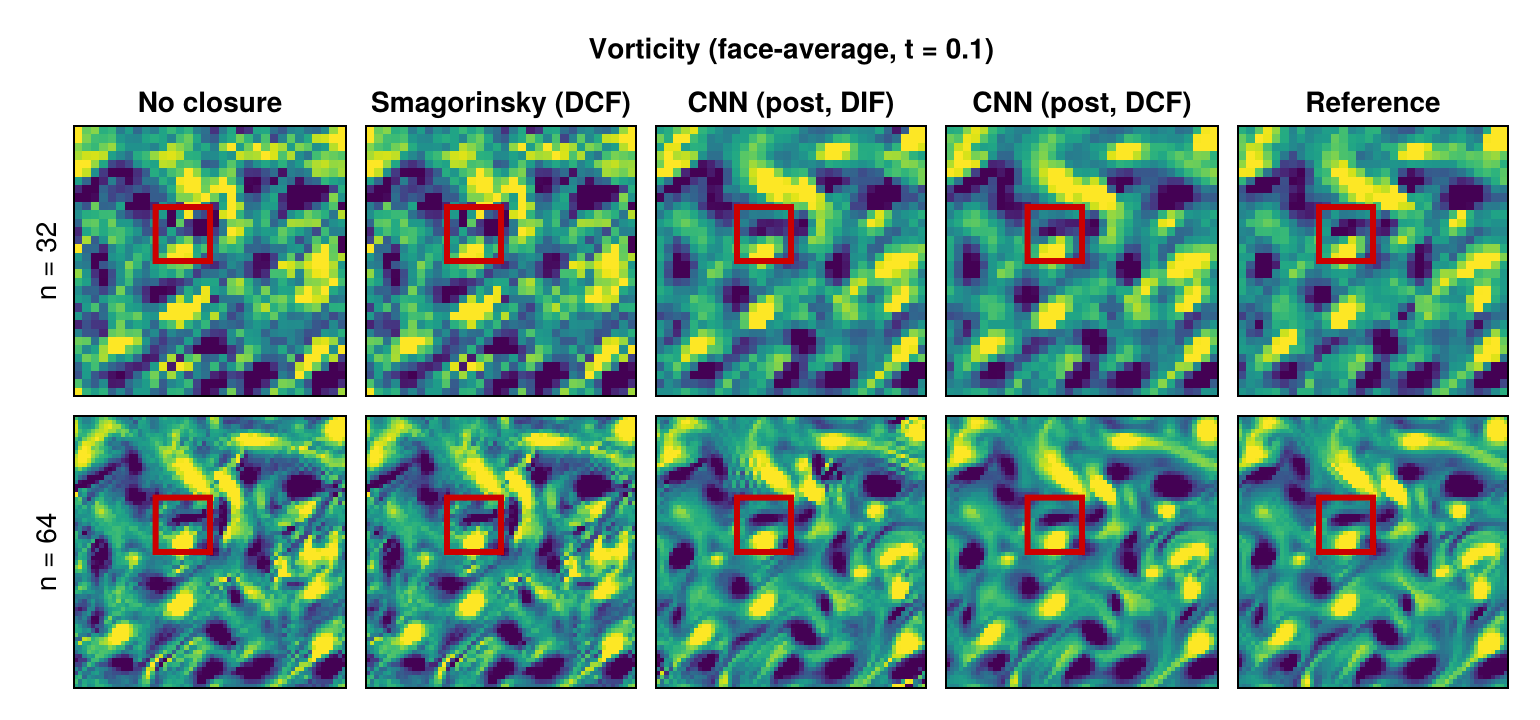}
    \includegraphics[width=\textwidth]{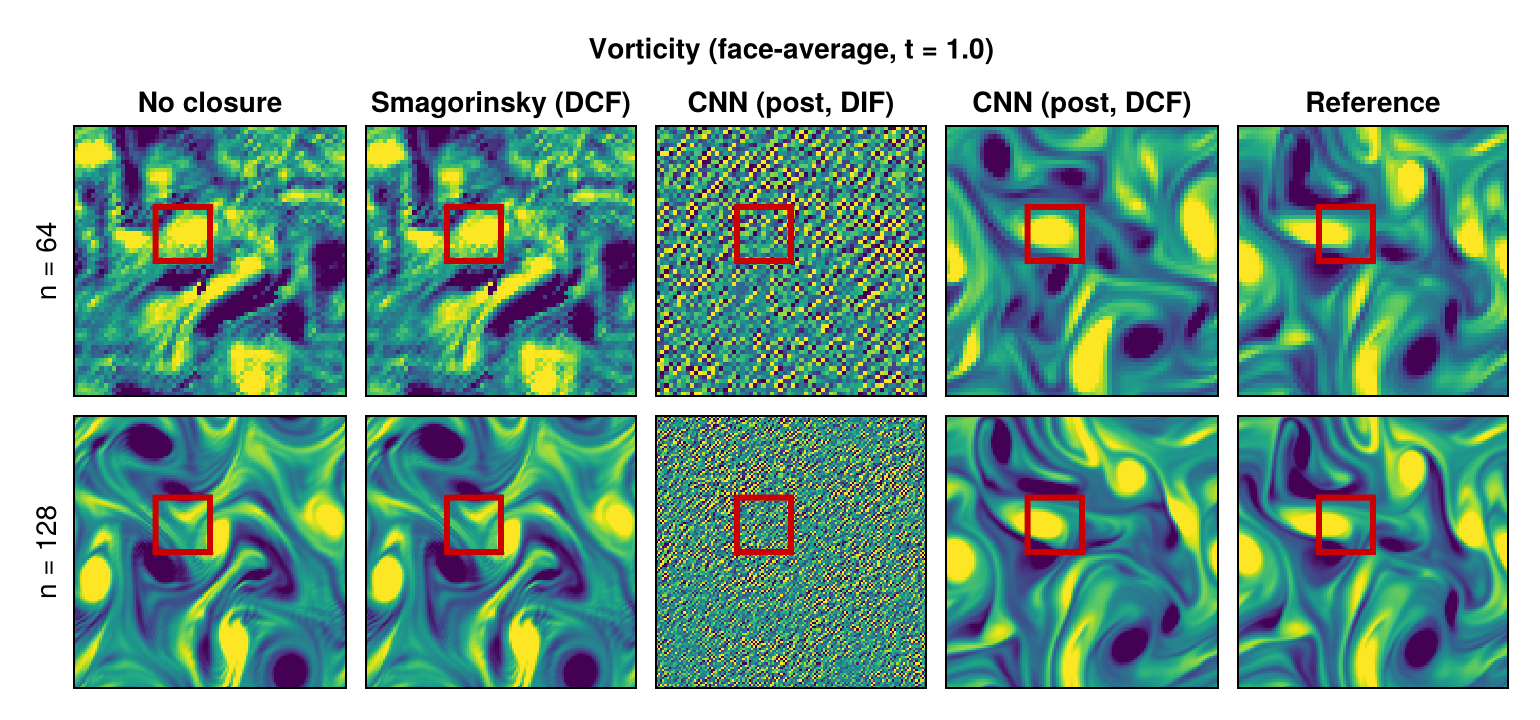}
    \caption{
        Vorticity of filtered DNS solution $\bar{u}$ (right) and LES solutions
        $\bar{v}$ computed using
        no closure (first column),
        Smagorinsky model (second column),
        $\mathcal{M}_{\text{DIF}}$ (third column),
        and $\mathcal{M}_{\text{DCF}}$ (fourth column)
        at \revboth{ different times}.
    }
    \label{fig:les_fields}
\end{figure}

The curl of the LES solutions (LES vorticity) at the final time is shown in
figure \ref{fig:les_fields}. The reference solution $\bar{u}$ has a smooth
vorticity field for all resolutions, since the low-pass filtering operation
removes high wavenumber components. The no-closure model solution, on the other
hand, shows sharp oscillations. Only the initial conditions are smooth, since
$\bar{v}(0) = \bar{u}(0)$. As the solution evolves in time, the central
difference scheme used in the right hand side $\bar{P} \bar{F}$ is too coarse
for the given Reynolds numbers, and produces well known oscillations
\cite{Shyy1985}. For
\revboth{ higher resolutions},
these oscillations go away, and $\bar{v}$ starts visually resembling $\bar{u}$.

Adding an $\mathcal{M}_{\text{DCF}}$ constrained CNN closure term
trained using $L^\text{post}_{\text{DCF}}$ seems to correct for the
oscillations of the no-closure model, and we recover the smooth fields with
recognizable features from the reference solution (compare $\bar{u}$ and
$\bar{v}$ inside the red square). We are effectively doing \emph{large eddy}
simulation, as large eddies are visually found in the right positions after
simulation. However, if we remove the projection and use the
$\mathcal{M}_{\text{DIF}}$ CNN closure model, the solution becomes
unstable, even after accounting for this by training with
$L^\text{post}_{\text{DIF}}$.
This confirms the observations from figure~\ref{fig:convergence}.

We stress again that the CNN closure model is accounting for the
total commutator error resulting from the discrete filtering procedure. This
commutator error includes the coarse grid discretization error, and also the
oscillations produced by the central difference scheme.

\subsubsection{Divergence}

\begin{figure}
    \centering
    \includegraphics[width=\textwidth]{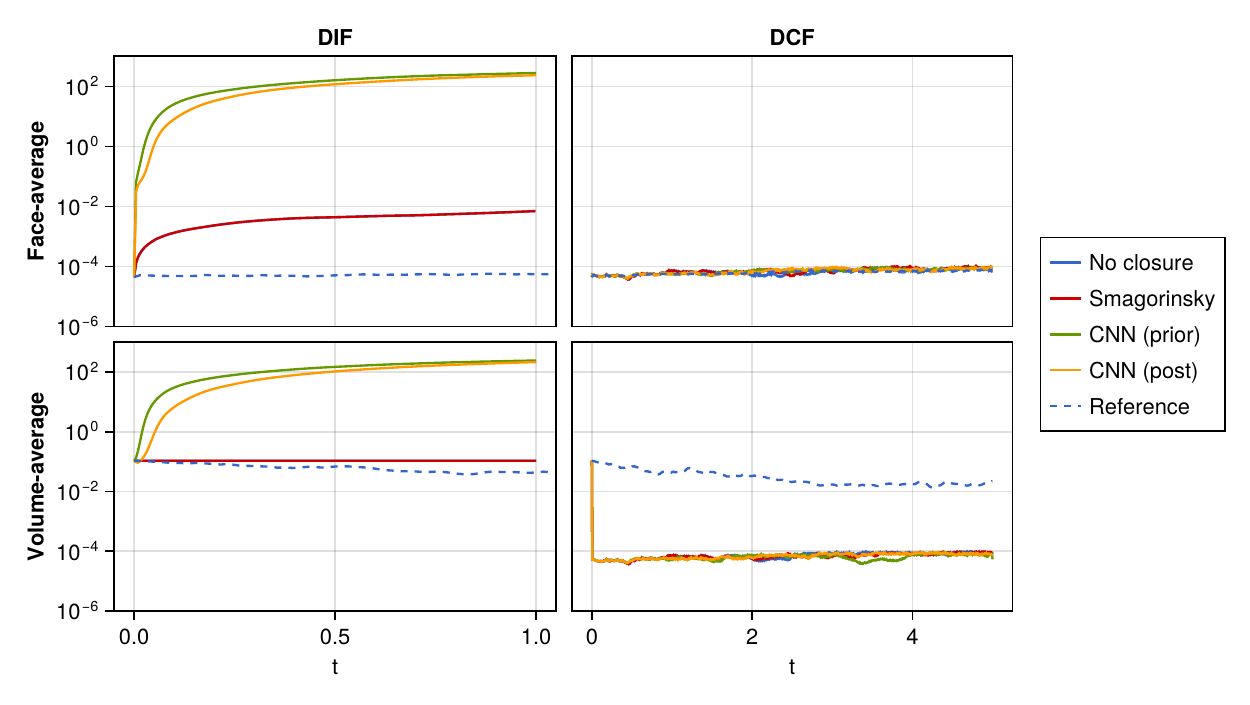}
    \caption{Divergence evolution for
        \revboth{ $\bar{n} = 256$}.
        \textbf{Left:} Unprojected closure model $\mathcal{M}_{\text{DIF}}$.
        \textbf{Right:} Constrained model $\mathcal{M}_{\text{DCF}}$.
        \textbf{Top:} Face-averaging filter. \textbf{Bottom:}
        Volume-averaging filter.
    }
    \label{fig:divergence_evolution}
\end{figure}

Figure~\ref{fig:divergence_evolution} shows the evolution of the average
divergence $\sqrt{\frac{1}{\bar{N}} \sum_I (\bar{D} \bar{v})_I^2}$ for the
two LES models
$\mathcal{M}_{\text{DIF}}$ and
$\mathcal{M}_{\text{DCF}}$.

For the face-averaging filter, the filtered DNS divergence (reference) is of the
order of
\revboth{ $10^{-4}$}
due to the single precision arithmetic. For $\mathcal{M}_{\text{DIF}}$,
all closure models produce a divergence increasing
in time. This is likely due to the instability observed in the previous
figures.
For $\mathcal{M}_{\text{DCF}}$, all closure models produce divergences at the same
order of magnitude as the reference, since the solution is projected at every
time step.

For the volume-averaging filter, the filtered DNS divergence (reference) is of
the order of \revboth{ $10^{-1}$}, since $\Phi^\text{VA}$ does not
preserve the divergence constraint. For $\mathcal{M}_{\text{DIF}}$, all closure models
produce an increasing divergence, just like for $\Phi^\text{FA}$.
For $\mathcal{M}_{\text{DCF}}$, all closure models produce divergence free solutions,
even though the reference solution is actually \emph{not} divergence-free.

\subsubsection{Computational cost}

\begin{table}
    \centering
    % \revboth{
    \begin{tabular}{c c c c c c}
    \toprule

    $\bar{u}$
    & Filter
        & $\bar{n}$
                  & $L^\text{prior}$
                               & $L^\text{post}_{\text{DIF}}$
                                          & $L^\text{post}_{\text{DCF}}$ \\

    \midrule

    \multirow{8}{*}{$1.5739\e{4}$} &
      FA & $\p32$ & $\p160.6$ & $4460.6$ & $4449.3$ \\
    & FA & $\p64$ & $\p244.0$ & $4509.1$ & $4584.7$ \\
    & FA &  $128$ & $\p497.4$ & $4901.2$ & $4853.9$ \\
    & FA &  $256$ & $ 1562.3$ & $5822.8$ & $5803.7$ \\
    & VA & $\p32$ & $\p157.8$ & $4490.8$ & $4569.4$ \\
    & VA & $\p64$ & $\p243.9$ & $4498.5$ & $4608.6$ \\
    & VA &  $128$ & $\p513.2$ & $4914.2$ & $4873.3$ \\
    & VA &  $256$ & $ 1583.6$ & $5793.1$ & $5849.2$ \\

    \bottomrule
\end{tabular}

% vim: conceallevel=0 textwidth=0

    % }
    \caption{
        Computational time (in seconds) for DNS and filtering of training,
        validation, and testing data ($\bar{u}$), a-priori training for 10000
        iterations ($L^\text{prior}$), and a-posteriori training for
        \revone{ 1000}
        iterations ($L^\text{post}$).
    }
    \label{tab:comptime}
\end{table}

The computational time for generating filtered DNS data and training the neural
networks is shown in table~\ref{tab:comptime}.
\revboth{
    
    One single DNS trajectory is used to compute filtered DNS trajectories 
    for all filters and coarse grid sizes. In total, we compute 8 DNS trajectories.
}
The two models $\mathcal{M}_{\text{DIF}}$ and $\mathcal{M}_{\text{DCF}}$
both produce similar timings, as exactly the same number of
operations are performed (but in a different order). The similar times for
\revboth{
    
    the three lowest resolutions
}
are likely due to the fact
that the same number of time steps are unrolled for all LES resolutions,
resulting in an equal number of GPU kernel calls. Since
\revboth{ $32^2$} and $128^2$ are both relatively small,
the kernels for the differential operators (convection, Poisson solves, etc.)
do not show any \revboth{ speed-up} for the lower resolution.
This does not seem to be the case for the $L^\text{prior}$ training however,
where most of the kernel calls are to the highly optimized convolutional
operators in CUDNN~\cite{Chetlur2014}.

The benefit of a-posteriori training is clear in situations where the a-priori
trained model is unstable (in our case: $\mathcal{M}_\text{DIF}$), but for a-priori
trained models that are stable and accurate
(such as $\mathcal{M}_\text{DCF}$ \revboth{for the face-averaging filter}),
one could ask whether the additional cost of a-posteriori training is worth it. We
think the additional accuracy is useful, since the learned weights can be reused
without having to retrain the model, as long as the same configuration is used
(grid size, Reynolds number etc.).

\section{Conclusion} \label{sec:conclusion}

The use of neural networks for LES closure models is a promising approach.
Neural networks are discrete by nature, and thus require consistent discrete
training data. To achieve model-data consistency, we propose the paradigm
``discretize, differentiate constraint, filter, and close'' (in that particular
order). This ensures full model-data consistency and circumvents pressure
problems. We do not need to define a pressure filter, only a velocity filter is
needed. Our framework allows for training the neural closure model using
a-priori loss functions, where the model is trained to predict the target
commutator error directly, or using a-posteriori loss functions, where the model
is trained to approximate the target filtered DNS-solution trajectory.

To ensure that the filtered DNS velocity stays divergence-free, we
\revone{
    \R{novel7}
    
    employed the divergence-consistent face-averaging filter used by
    Kochkov et al.~\cite{Kochkov2021}
}.
This allows for using a divergence-constrained LES model similar to those
commonly used in LES, but with the important difference that the training data
obtained through discrete DNS is fully consistent with the LES environment. This
resulted in our new divergence-consistent LES model, which was found to be
stable with both a-priori and a-posteriori training.
\revboth{Using the same formulation with a divergence-inconsistent filter (VA)
did however require a-posteriori training to achieve the same accuracy}.
A-posteriori training is more expensive than a-priori training, but still has
the advantage of increased accuracy.

\revone{\R{novel8}  The}
divergence-consistent filter stands out from commonly used
(volume-averaging) discrete filters, for which the filtered DNS solution is
generally not divergence-free. We showed that the resulting LES model can
produce instabilities. A-posteriori trained models were found to improve the
stability over a-priori trained models, but this was not sufficient to fully
stabilize the model in our experiment. With
\revone{ a divergence-consistent formulation},
such stability issues did not occur.
\revone{
    \R{noise}
    Another common approach to achieve stability is to add noise to the training
    data. However, this could destroy the divergence-free property of the data,
    and appropriate measures would need to be taken to avoid this.
    We tried adding different levels of noise to the training data, but
    the resulting validation error was found to be better without noise.
}

Another important property of our \revboth{discrete} approach is that the (coarse grid)
discretization error is included in the training data and learned by the neural
network. Turbulence simulations with DNS and LES rely on non-dissipative
discretization methods, such as the second order central difference
discretization we used in this work, but they produce oscillations and
instabilities on coarse grids. This could limit their use in LES, where using
coarse grids is one of the main goals. Various smoothing methods, such as
explicit LES or (overly) diffusive closure models are commonly used to address
this issue. We let the closure model learn to account for the oscillations. The
fact that the coarse grid discretization effects (and thus the oscillations) are
included in the training data allows the neural network to recognize and correct
for these oscillations, even when training with a-priori loss functions.

We realize that the choice of a divergence-consistent face-averaging filter
imposes some constraints on the filter choice. The weights have to be uniform
(top-hat filter like) to ensure that all the sub-filter velocities cancel out,
and the extension to other filters like Gaussians is an open problem. In
addition, we took the filter width to be equal to the grid spacing. The
face-average filter naturally extends to unstructured grids, as shown in
figure~\ref{fig:discretization}. However, this requires that the DNS grid perfectly
overlaps with the faces of the LES grid. This can be achieved by designing the
DNS and LES grids at the same time. Lastly,
\revone{\R{novel9}  the face-averaging}
filter is only divergence-consistent with respect to the
second-order central difference divergence operator on a staggered grid. We
intend to explore the use of filters that are divergence-consistent with respect
to higher-order discretization methods (such as the discontinuous Galerkin
element method) and non-staggered grids.
\revone{
    \R{collocated}
    The ``discretize first'' approach can also be applied on collocated grids.
    While is possible to define a face-averaging filter on a collocated
    \emph{Cartesian} grid, the divergence-consistent property of the filter
    would no longer hold. On an unstructured collocated grid, it might not make
    sense to define a face-averaging filter.
}

Future work will include the use of non-uniform grids with solid walls boundary
conditions. This modification will require a change in the neural network
architecture, but not in the face-averaging filtering procedure itself. The CNN
architecture assumes that the grid is uniform, but weighting the kernel stencil
with non-uniform volume sizes could be a way to design a discretization-informed
neural network. We intend to exploit the fact that discrete DNS boundary
conditions naturally extend to the LES model for the face-average filter. We
also intend to incorporate other constraints into the LES model, such as
conservation of quantities of interest, in particular kinetic energy
conservation, as studied in~\cite{Vangastelen2023}.

% Statements
% \section{Reproducibility statement} \label{sec:reproducibility}
\section*{Software and reproducibility statement}

The Julia scripts and source code used to generate all the results are
available at \\
\url{https://github.com/agdestein/IncompressibleNavierStokes.jl}. \\
It is released under the MIT license.
\revtwo{
    
    Each simulation was run on a single Nvidia H100 GPU
    on the Dutch national Snellius supercomputer.
}

To produce the results of this article, pseudo-random numbers were used for
\begin{itemize}
    \item data generation (DNS initial conditions for training, validation, and
        testing data);
    \item neural network parameter initialization;
    \item batch selection during stochastic gradient descent.
\end{itemize}
The scripts include the seeding numbers used to initialize the pseudo-random
number generators. Given a seeded pseudo-random generator, the code is
deterministic and the results are reproducible.

\section*{CRediT author statement}

% https://www.elsevier.com/researcher/author/policies-and-guidelines/credit-author-statement
% https://www.sciencedirect.com/journal/journal-of-computational-physics/publish/guide-for-authors#9060

\textbf{Syver Døving Agdestein}:
Conceptualization,
Methodology,
Software,
Validation,
Visualization,
Writing -- Original Draft

\textbf{Benjamin Sanderse}:
Formal analysis,
Funding acquisition,
Project administration,
Supervision,
Writing -- Review \& Editing

\section*{Declaration of Generative AI and AI-assisted technologies in the writing process}

% https://www.sciencedirect.com/journal/journal-of-computational-physics/publish/guide-for-authors#7300

During the preparation of this work the authors used GitHub Copilot in order to
propose wordings and mathematical typesetting. After using this tool/service,
the authors reviewed and edited the content as needed and take full
responsibility for the content of the publication.

\section*{Declaration of competing interest}

% https://www.sciencedirect.com/journal/journal-of-computational-physics/publish/guide-for-authors#7250

The authors declare that they have no known competing financial interests or
personal relationships that could have appeared to influence the work reported
in this paper.

\section*{Acknowledgements}

This work is supported by the projects ``Discretize first, reduce next'' (with
project number VI.Vidi.193.105) of the research programme NWO Talent Programme
Vidi and ``Discovering deep physics models with differentiable programming''
(with project number EINF-8705) of the Dutch Collaborating University Computing
Facilities (SURF), both financed by the Dutch Research Council (NWO). We thank
SURF (www.surf.nl) for the support in using the Dutch National Supercomputer
Snellius.

% https://servicedesk.surf.nl/wiki/display/WIKI/Snellius+FAQ#SnelliusFAQ-Miscellaneous

\appendix

\section{Finite volume discretization} \label{sec:discretization}

In this work, the integral form of the Navier-Stokes equations is considered,
which is used as starting point to develop a spatial discretization:
\begin{align}
    \frac{1}{|\mathcal{O}|} \int_{\partial \mathcal{O}}
    u \cdot n \, \mathrm{d} \Gamma & = 0, \label{eq:mass_integral} \\
    \frac{\mathrm{d} }{\mathrm{d} t} \frac{1}{|\mathcal{O}|}
    \int_\mathcal{O} u \, \mathrm{d} \Omega
    & = \frac{1}{|\mathcal{O}|} \int_{\partial \mathcal{O}}
    \left( - u u^\mathsf{T} - p I + \nu \nabla u \right) \cdot n \,
    \mathrm{d} \Gamma +
    \frac{1}{|\mathcal{O}|}\int_\mathcal{O} f \mathrm{d} \Omega,
    \label{eq:momentum_integral}
\end{align}
where $\mathcal{O} \subset \Omega$ is an arbitrary control volume with boundary
$\partial \mathcal{O}$, normal $n$, surface element $\mathrm{d} \Gamma$, and
volume size $|\mathcal{O}|$. We have divided by the control volume sizes in the
integral form, so that system
\eqref{eq:mass_integral}-\eqref{eq:momentum_integral} has the same dimensions as
the system
\eqref{eq:mass_continuous}-\eqref{eq:momentum_continuous}.

\subsection{Staggered grid configuration}

In this section we describe a finite volume discretization of equations
\eqref{eq:mass_integral}-\eqref{eq:momentum_integral}. Before doing so, we
introduce our notation, which is such that the mathematical description of the
discretization closely matches the software implementation.

The $d$ spatial dimensions are indexed by $\alpha \in \{1, \dots, d\}$.
The $\alpha$-th unit vector is denoted $2 h_\alpha = (2 h_{\alpha \beta})_{\beta
= 1}^d$, where the (half) Kronecker symbol $h_{\alpha \beta}$ is $1 / 2$ if
$\alpha = \beta$ and $0$ otherwise. The Cartesian index $I = (I_1, \dots,
I_d)$ is used to avoid repeating terms and equations $d$ times, where
$I_\alpha$ is a scalar index (typically one of $i$, $j$, and $k$ in common
notation). This notation is dimension-agnostic, since we can write $u_I$ instead
of $u_{i j}$ in 2D or $u_{i j k}$ in 3D. In our Julia implementation of the
solver we use the same Cartesian notation (\verb|u[I]| instead of \verb|u[i, j]|
or \verb|u[i, j, k]|).

\begin{figure}
    \centering
    \def\svgwidth{0.60\columnwidth}
    %% Creator: Inkscape 1.4 (e7c3feb100, 2024-10-09), www.inkscape.org
%% Accompanies image file 'figures_finitevolumes.pdf' (pdf, eps, ps)
\begingroup%
  \makeatletter%
  \providecommand\color[2][]{%
    \errmessage{(Inkscape) Color is used for the text in Inkscape, but the package 'color.sty' is not loaded}%
    \renewcommand\color[2][]{}%
  }%
  \providecommand\transparent[1]{%
    \errmessage{(Inkscape) Transparency is used (non-zero) for the text in Inkscape, but the package 'transparent.sty' is not loaded}%
    \renewcommand\transparent[1]{}%
  }%
  \providecommand\rotatebox[2]{#2}%
  \newcommand*\fsize{\dimexpr\f@size pt\relax}%
  \newcommand*\lineheight[1]{\fontsize{\fsize}{#1\fsize}\selectfont}%
  \ifx\svgwidth\undefined%
    \setlength{\unitlength}{553.25643224bp}%
    \ifx\svgscale\undefined%
      \relax%
    \else%
      \setlength{\unitlength}{\unitlength * \real{\svgscale}}%
    \fi%
  \else%
    \setlength{\unitlength}{\svgwidth}%
  \fi%
  \global\let\svgwidth\undefined%
  \global\let\svgscale\undefined%
  \makeatother%
  \begin{picture}(1,0.71221847)%
    \lineheight{1}%
    \setlength\tabcolsep{0pt}%
    \put(0,0){\includegraphics[width=\unitlength,page=1]{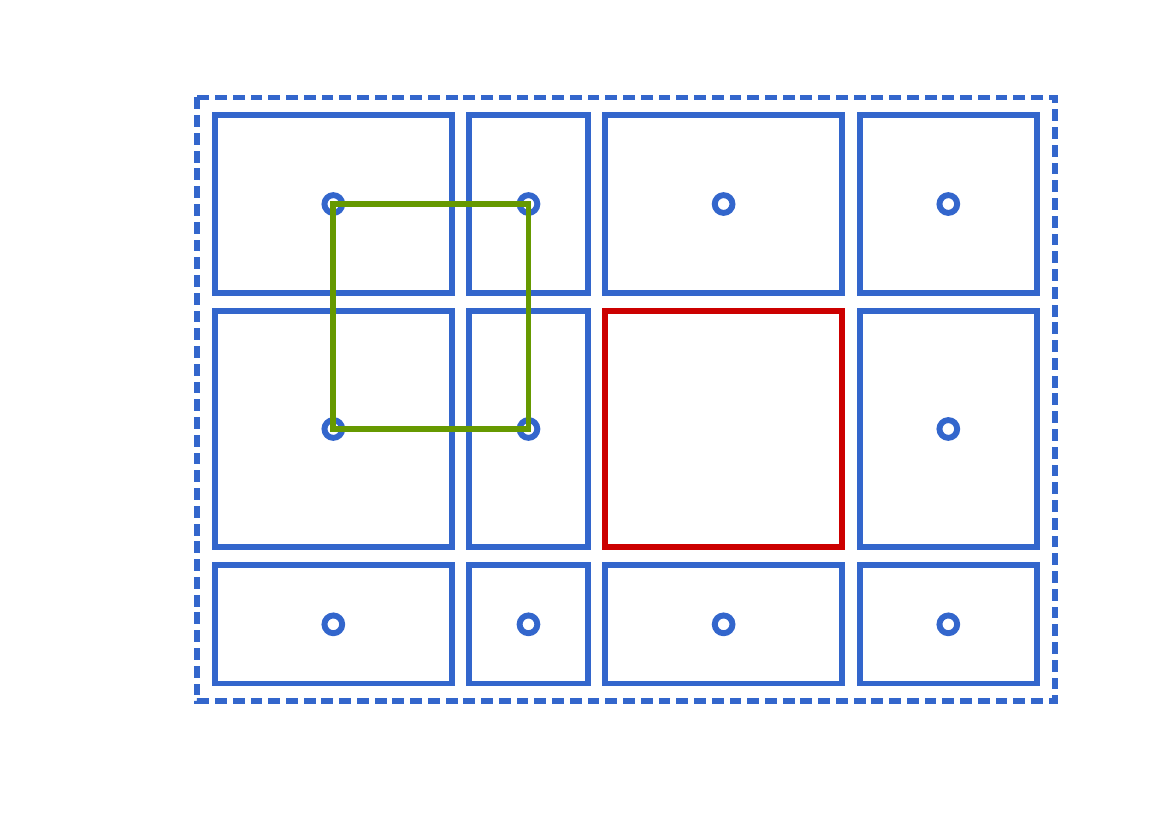}}%
    \put(0.59472317,0.64743233){\color[rgb]{0.8,0,0}\makebox(0,0)[lt]{\lineheight{1.25}\smash{\begin{tabular}[t]{l}$\Delta^1_{I_1}$\end{tabular}}}}%
    \put(0.92601781,0.44598741){\color[rgb]{0.4,0.6,0}\makebox(0,0)[lt]{\lineheight{1.25}\smash{\begin{tabular}[t]{l}$\Delta^2_{J_2}$\end{tabular}}}}%
    \put(0.30008766,0.35732246){\color[rgb]{0.4,0.6,0}\makebox(0,0)[lt]{\lineheight{1.25}\smash{\begin{tabular}[t]{l}$\Omega_J$\end{tabular}}}}%
    \put(0.17282437,0.6435522){\color[rgb]{0.2,0.4,0.8}\makebox(0,0)[lt]{\lineheight{1.25}\smash{\begin{tabular}[t]{l}$\Omega$\end{tabular}}}}%
    \put(0,0){\includegraphics[width=\unitlength,page=2]{figures_finitevolumes.pdf}}%
    \put(0.00600755,0.33330426){\color[rgb]{0.8,0,0}\makebox(0,0)[lt]{\lineheight{1.25}\smash{\begin{tabular}[t]{l}$I_2$\end{tabular}}}}%
    \put(0.38195858,0.01333494){\color[rgb]{0.4,0.6,0}\makebox(0,0)[lt]{\lineheight{1.25}\smash{\begin{tabular}[t]{l}$J_1$\end{tabular}}}}%
    \put(0.00732338,0.44236409){\color[rgb]{0.8,0,0}\makebox(0,0)[lt]{\lineheight{1.25}\smash{\begin{tabular}[t]{l}$I_2 + \frac{1}{2}$\end{tabular}}}}%
    \put(0.00732338,0.22364646){\color[rgb]{0.8,0,0}\makebox(0,0)[lt]{\lineheight{1.25}\smash{\begin{tabular}[t]{l}$I_2 - \frac{1}{2}$\end{tabular}}}}%
    \put(0.44456623,0.01333494){\color[rgb]{0.4,0.6,0}\makebox(0,0)[lt]{\lineheight{1.25}\smash{\begin{tabular}[t]{l}$J_1 + \frac{1}{2}$\end{tabular}}}}%
    \put(0.22602366,0.01333494){\color[rgb]{0.4,0.6,0}\makebox(0,0)[lt]{\lineheight{1.25}\smash{\begin{tabular}[t]{l}$J_1 - \frac{1}{2}$\end{tabular}}}}%
    \put(0.88680947,0.01166483){\color[rgb]{0.2,0.4,0.8}\makebox(0,0)[lt]{\lineheight{1.25}\smash{\begin{tabular}[t]{l}$x^1$\end{tabular}}}}%
    \put(0.05299244,0.62692028){\color[rgb]{0.2,0.4,0.8}\makebox(0,0)[lt]{\lineheight{1.25}\smash{\begin{tabular}[t]{l}$x^2$\end{tabular}}}}%
    \put(0,0){\includegraphics[width=\unitlength,page=3]{figures_finitevolumes.pdf}}%
    \put(0.58386744,0.36856025){\color[rgb]{0.8,0,0}\makebox(0,0)[lt]{\lineheight{1.25}\smash{\begin{tabular}[t]{l}$p_I$\end{tabular}}}}%
    \put(0.5335305,0.25416113){\color[rgb]{0.8,0,0}\makebox(0,0)[lt]{\lineheight{1.25}\smash{\begin{tabular}[t]{l}$\Omega_I$\end{tabular}}}}%
    \put(0.75402925,0.37263451){\color[rgb]{0.8,0,0}\makebox(0,0)[lt]{\lineheight{1.25}\smash{\begin{tabular}[t]{l}$u^1_{I + h_1}$\end{tabular}}}}%
    \put(0.64584411,0.47342431){\color[rgb]{0.8,0,0}\makebox(0,0)[lt]{\lineheight{1.25}\smash{\begin{tabular}[t]{l}$u^2_{I + h_2}$\end{tabular}}}}%
  \end{picture}%
\endgroup%

    \caption{
        Finite volume discretization on a staggered grid. Note that the grid can
        be non-uniform, as long as each volume in a given column has the same
        width and each volume in a given row has the same height. Here, $I$ and
        $J$ are two arbitrary Cartesian indices, with $I \in \mathcal{I}$ in a
        volume center and $J \in \mathcal{I} + h_1 + h_2$ in a volume corner for
        illustrative purposes.
    }
    \label{fig:finitevolumes}
\end{figure}

For the discretization scheme, we use a staggered Cartesian grid as proposed by
Harlow and Welch~\cite{Harlow1965}. Staggered grids have excellent conservation
properties~\cite{Lilly1965,Perot2011}, and in particular their exact
divergence-freeness is important for this work. Consider a rectangular domain
$\Omega = \prod_{\alpha = 1}^d [a_\alpha, b_\alpha]$, where $a_\alpha <
b_\alpha$ are the domain boundaries and $\prod$ is a Cartesian product. Let
$\Omega = \bigcup_{I \in \mathcal{I}} \Omega_I$ be a partitioning of $\Omega$,
where $\mathcal{I} = \prod_{\alpha = 1}^d \{ \frac{1}{2}, 2 - \frac{1}{2},
\dots, N_\alpha - \frac{1}{2} \}$ are volume center indices, $N = (N_1, \dots,
N_d) \in \mathbb{N}^d$ are the number of volumes in each dimension,
$\Omega_I = \prod_{\alpha = 1}^d \Delta^\alpha_{I_\alpha}$ is a finite
volume, $\Gamma^\alpha_I = \Omega_{I - h_\alpha} \cap \Omega_{I + h_\alpha} =
\prod_{\beta \neq \alpha} \Delta^\beta_{I_\beta}$ is a volume face,
$\Delta^\alpha_i = \left[ x^\alpha_{i - \frac{1}{2}}, x^\alpha_{i + \frac{1}{2}}
\right]$ is a volume edge, $x^\alpha_0, \dots, x^\alpha_{N_\alpha}$ are volume
boundary coordinates, and $x^\alpha_i = \frac{1}{2} \left(x^\alpha_{i -
\frac{1}{2}} + x^\alpha_{i + \frac{1}{2}}\right)$ for $i \in \{ 1 / 2, \dots,
N_\alpha - 1 / 2\}$ are volume center coordinates. We also define the operator
$\delta_\alpha$ which maps a discrete scalar field $\varphi = (\varphi_I)_I$ to
\begin{equation}
    (\delta_\alpha \varphi)_I = \frac{\varphi_{I + h_\alpha} - \varphi_{I -
    h_\alpha}}{| \Delta^\alpha_{I_\alpha} |}.
\end{equation}
It can be interpreted as a discrete equivalent of the continuous operator
$\frac{\partial}{\partial x^\alpha}$. All the above definitions are extended to
be valid in volume centers $I \in \mathcal{I}$, volume faces $I \in \mathcal{I}
+ h_\alpha$, or volume corners $I \in \mathcal{I} + \sum_{\alpha = 1}^d
h_\alpha$. The discretization is illustrated in figure~\ref{fig:finitevolumes}.

\subsection{Equations for unknowns}

We now define the unknown degrees of freedom. The average pressure in
$\Omega_I$, $I \in \mathcal{I}$ is approximated by the quantity $p_I(t)$. The
average $\alpha$-velocity on the face $\Gamma^\alpha_I$, $I \in \mathcal{I} +
h_\alpha$ is approximated by the quantity $u^\alpha_I(t)$. Note how the pressure
$p$ and the $d$ velocity fields $u^\alpha$ are each defined in their own
canonical positions $x_I$ and $x_{I + h_\alpha}$ for $I \in \mathcal{I}$. This
is illustrated for a given volume $I$ in figure~\ref{fig:finitevolumes}. In the
following, we derive equations for these unknowns.

Using the pressure control volume $\mathcal{O} = \Omega_I$ with $I \in
\mathcal{I}$ in the integral constraint \eqref{eq:mass_integral} and
approximating the face integrals with the mid-point quadrature rule
$\int_{\Gamma_I} u \, \mathrm{d} \Gamma \approx | \Gamma_I | u_I$ results in the
discrete divergence-free constraint
\begin{equation} \label{eq:mass_discrete}
    \sum_{\alpha = 1}^d
    (\delta_\alpha u^\alpha)_I = 0.
\end{equation}
Note how dividing by the volume size results in a discrete equation resembling
the continuous one (since $| \Omega_I | = | \Gamma^\alpha_I | |
\Delta^\alpha_{I_\alpha} |$).

Similarly, choosing an $\alpha$-velocity control volume $\mathcal{O} =
\Omega_{I}$ with $I \in \mathcal{I} + h_\alpha$ in equation
\eqref{eq:momentum_integral}, approximating the volume- and face integrals using
the mid-point quadrature rule, and replacing remaining spatial derivatives in
the diffusive term with a finite difference approximation gives the discrete
momentum equations
\begin{equation} \label{eq:momentum_discrete}
        \frac{\mathrm{d}}{\mathrm{d} t} u^\alpha_{I} =
        - \sum_{\beta = 1}^d
        (\delta_\beta (u^\alpha u^\beta))_{I}
        + \nu \sum_{\beta = 1}^d
        (\delta_\beta \delta_\beta u^\alpha)_{I}
        + f^\alpha(x_{I})
        - (\delta_\alpha p)_{I}.
\end{equation}
where we made the assumption that $f$ is constant in time for simplicity.
The outer discrete derivative in $(\delta_\beta \delta_\beta u^\alpha)_{I}$ is
required at the position $I$, which means that the inner derivative is evaluated
as $(\delta_\beta u^\alpha)_{I + h_\beta}$ and $(\delta_\beta u^\alpha)_{I -
h_\beta}$, thus requiring $u^\alpha_{I - 2 h_\beta}$, $u^\alpha_{I}$, and
$u^\alpha_{I + 2 h_\beta}$, which are all in their canonical positions. The two
velocity components in the convective term $u^\alpha u^\beta$ are required at
the positions $I - h_\beta$ and $I + h_\beta$, which are outside the canonical
positions. Their value at the required position is obtained using averaging with
weights $1 / 2$ for the $\alpha$-component and with linear interpolation for the
$\beta$-component. This preserves the skew-symmetry of the convection operator,
such that energy is conserved (in the convective term)~\cite{Verstappen2003}.

\section{Numerical experiment details} \label{sec:experiments}

\subsection{Energy spectra} \label{sec:energy}

For our discretization, we define the energy at a wavenumber $k$ as $\hat{E}_k =
\frac{1}{2} \| \hat{u}_k \|^2$, where $\hat{u}^\alpha =
\operatorname{DFT}(u^\alpha)$ is the discrete Fourier transform of $u$. Since $k
\in \mathbb{Z}^d$, it is not immediately clear how to compute a discrete
equivalent of the scalar energy spectrum as a function of $\| k \|$. We proceed
as follows. The energy at a scalar level $\kappa > 0$ is defined as the sum over
all energy components of the dyadic bin $\mathcal{K}_\kappa = \{ k \ | \ \kappa
/ a \leq \| k \| \leq \kappa a \}$ as
\begin{equation} \label{eq:energyspectrum}
    \hat{E}_\kappa = \sum_{k \in \mathcal{K}_\kappa} \hat{E}_k.
\end{equation}
The parameter $a > 1$ determines the width of the interval. Lumley argues to
use the golden ratio $a = (1 + \sqrt{5}) / 2 \approx 1.6$~\cite{Gatski1996}.
Note that there is no averaging factor in front of the sum
\eqref{eq:energyspectrum}, even though the number of wavenumbers in the set
increases with $\kappa$.

For homogeneous decaying isotropic turbulence, the spectrum should behave as
follows. In the inertial region, for large Reynolds numbers, the theoretical
decay of $\hat{E}$ should be $\hat{E}_\kappa = \mathcal{O}(\kappa^{-3})$ in 2D
and $\mathcal{O}(\kappa^{-5/3})$ in 3D~\cite{Pope2000}. For the lowest
wavenumbers, we should have $\hat{E}_\kappa = \mathcal{O}(\kappa^4)$ or
$\hat{E}_\kappa = \mathcal{O}(\kappa^2)$.

It is also common to use a linear bin such as $\mathcal{K}_\kappa = \{k \, | \,
\kappa - \frac{1}{2} \leq \| k \| < \kappa + \frac{1}{2} \}$
\cite{Orlandi2000,San2012,Maulik2017}. However, this leads to a different power
law scaling in the inertial range than the well known $\kappa^{-3}$ and
$\kappa^{-5/3}$.

\subsection{Initial conditions} \label{sec:ic}

To generate initial conditions in a periodic box, we consider a prescribed
energy spectrum $\hat{E}_k$. We want to create an initial velocity field $u$
with the following properties:
\begin{itemize}
    \item The Fourier transform of $u$, noted $\hat{u}$, should be such that
        $\frac{1}{2} \| \hat{u}_k \|^2 = \hat{E}_k$ for all $k$.
    \item $u$ should be divergence-free with respect to our
        discretization: $D u = 0$.
    \item $u$ should be parameterized by controllable random numbers, such
        that a wide variety of initial conditions can be generated.
\end{itemize}
These properties are achieved by sampling a velocity field in spectral space,
projecting (making it divergence-free), transforming to physical space, and
projecting again. In detail, let $a_k = \sqrt{2 \hat{E}_k} \mathrm{e}^{2 \pi
\mathrm{i} \tau_k}$, where  $\tau_k = \sum_{\alpha = 1}^d \xi_k^\alpha$ is phase
shift, $\xi_k^\alpha \sim \mathcal{U}[0, 1]$ is a random uniform number if
$k_\beta \geq 0$ for all $\beta$. If $k_\beta < 0$ for any $\beta$, we add the
symmetry constraint $\xi_k^\alpha = \operatorname{sign}(k_\alpha)
\xi_{|k|}^\alpha$ where $|k| = (|k_\beta|)_{\beta = 1}^d$. Then $\| \hat{u}_k \|
= | a_k |$, and $a_k$ has a random phase shift. We then multiply the scalar
$a_k$ with a random unit vector $e_k$ projected onto the divergence-free
spectral grid as follows: $\hat{u}^\alpha_k = a_k \hat{P}_k e^\alpha_k / \|
\hat{P}_k e^\alpha_k \|$, where $\hat{P}_k = I - \frac{k
k^\mathsf{T}}{k^\mathsf{T} k} \in \mathbb{C}^{d \times d}$ is a projector
for each $k$ ensuring that $2 \pi \mathrm{i} k^\mathsf{T} \hat{u}_k = 0$ for all
$k$ (which is the equivalent of $\nabla \cdot u = 0$ in spectral space)
\cite{Pope2000}. The normalization with respect to $\| \hat{P}_k e^\alpha_k \|$
ensures that no energy is lost in the projection step. In 2D, we choose a random
vector on the unit circle $e_k = (\cos(\theta_k), \sin(\theta_k))$ with
$\theta_k \sim \mathcal{U}[0, 2 \pi]$. In 3D, we choose a random vector on the
unit sphere $e_k = (\sin(\theta_k) \cos(\phi_k), \sin(\theta_k) \sin(\phi_k),
\cos(\theta_k))$ with $\theta_k \sim \mathcal{U}[0, \pi]$ and $\phi_k \sim
\mathcal{U}[0, 2 \pi]$. Finally, we obtain the velocity field $u$ by taking the
inverse discrete Fourier transform, and also projecting it again since
divergence-freeness on the ``spectral grid'' and on the staggered grid are
slightly different. This gives the random initial field
\begin{equation}
    u = P \operatorname{DFT}^{-1}(\hat{u}).
\end{equation}
Note that the second projection may result in a slight loss of energy, but since
$u$ is already divergence-free on the spectral grid, our experience is that the loss is non-significant.

\subsection{CNN architecture} \label{sec:architecture}

The CNN architecture is shown in table \ref{tab:cnn}. We use periodic padding.
For the last convolutional layer, we use no activation and no bias, in order not
to limit the expressiveness of $m$. For the inner layers, we use bias and
the $\tanh$ activation function.

\begin{table}
    \centering
    \begin{tabular}{l c c c c c}
    \toprule

    Layer &
    Radius &
    Channels &
    Activation &
    Bias &
    Parameters \\
    
    \midrule

    Interpolate$_{u \to p}$   &     &               &               &     &           \\
    $\operatorname{Conv}$     & $2$ & $\p 2 \to 24$ & $\tanh$       & Yes & $\p 1224$ \\
    $\operatorname{Conv}$     & $2$ & $24 \to 24$   & $\tanh$       & Yes &   $14424$ \\
    $\operatorname{Conv}$     & $2$ & $24 \to 24$   & $\tanh$       & Yes &   $14424$ \\
    $\operatorname{Conv}$     & $2$ & $24 \to 24$   & $\tanh$       & Yes &   $14424$ \\
    $\operatorname{Conv}$     & $2$ & $24 \to \p 2$ & $x \mapsto x$ &  No & $\p 1200$ \\
    Interpolate$_{p \to u}$   &     &               &               &     &           \\ \midrule
                              &     &               &               &     &   $45696$
\end{tabular}

% vim: conceallevel=0 textwidth=0

    \caption{CNN architecture, with $\epsilon = 1 / 100$. The total radius is
        $8$, which means that the component $m(\bar{u}, \theta)_I$ depends on
        the components $\bar{u}_{I + J}$ for $J \in [-8, 8]^2$. In comparison,
        the diffusion operator has a radius of $1$.}
    \label{tab:cnn}
\end{table}

The choice of channel sizes is chosen solely to have sufficient expressive
capacity in the closure model. The kernel radius on the other hand, is chosen to
be small ($r = 2$, diameter $5$) to ensure that the closure model uses local
information. In addition, the resulting small stencils that are learned can 
possibly be interpreted as discrete differential operators of finite difference
type, separated by simple non-linearities. In this way, the CNN can be thought
of as a generalized Taylor series expansion of the commutator error in terms of
the filtered velocity field, similar to certain continuous filter expansions
\cite{Sagaut2005}. We do not investigate this further in this study, but it
could be a direction for future research.

The same CNN architecture $m$ is used for all grids and filters. We choose a
simple architecture since the goal of the study is not to get the most accurate
closure model, but rather to compare different filters, LES formulations, and
loss functions for the same closure architecture. For the LES formulation, we
only consider the two models $\mathcal{M}_{\text{DIF}}$ and $\mathcal{M}_{\text{DCF}}$.

\subsection{Data generation} \label{sec:datasets}

\revboth{
    
    To create data, we run $8$ DNS simulations of size $N = 4096^2$.
    We use adaptive time-stepping.
}
For every random initial flow field $u(0)$, we let the DNS
run for a burn-in time \revboth{ $t_\text{burn} = 0.5$}
to initialize the flow beyond the artificial initial spectrum
\eqref{eq:energyprofile}.
We then start saving $\bar{u}$ and $c$ every time step until $t_\text{end}
\revboth{= 5}$.
\revboth{
    Every $50$ time steps we compute $\bar{u}$ and $c(u)$ on four coarse
    grids of size $\bar{n} \in \{32, 64, 128, 256\}$. The first $6$ trajectories
    are used for training, and the remaining $2$ for validation and testing.
}

\subsection{Training} \label{sec:training}

\begin{table}
    \centering
    % \revone{
    \begin{tabular}{c c c c c}
    \toprule
    $\bar{n}$
         & $\theta^\text{FA}_{\text{DIF}}$
                 & $\theta^\text{VA}_{\text{DIF}}$
                         & $\theta^\text{FA}_{\text{DCF}}$
                                 & $\theta^\text{VA}_{\text{DCF}}$ \\
    \midrule
    \p32 & 0.129 & 0.131 & 0.144 & 0.142 \\
    \p64 & 0.114 & 0.116 & 0.143 & 0.141 \\
     128 & 0.061 & 0.064 & 0.096 & 0.093 \\
     256 & 0.000 & 0.000 & 0.023 & 0.014 \\
    \bottomrule
\end{tabular}

    % }
    \caption{\revone{Optimized Smagorinsky coefficients for each coarse resolution, filter type, and LES model.}}
    \label{tab:smagcoeffs}
\end{table}

Both the Smagorinsky model and the CNN are parameterized and require training.
Since the Smagorinsky parameter is a scalar, we perform a grid search to find
the optimal parameter for each of the \revone{grid sizes},
filter types and projection orders.
The relative a-posteriori error for the training set is evaluated
\revone{}.
We choose the value of
\revone{
    
    $\theta \in \{ 0, 1/1000, 2/1000, \dots, 300/1000 \}$
}
that gives the lowest training error.
The resulting Smagorinsky constants are
\revone{
    \R{smagcoeffs}
    
    shown in table~\ref{tab:smagcoeffs}.
    Note that the Smagorinsky constants are theoretically grid-independent,
    since information about the grid is incorporated using the filter width that
    enters the expression separately. We still optimize the coefficient for each
    grid size in order to achieve a fair comparison with the CNN.
    It is important to stress that the Smagorinsky closure term (divergence of
    weigthed strain tensor in equation \eqref{eq:smagorinsky_closure}) is
    \emph{not} divergence-free, just like the other right hand side terms like convection and diffusion. As a result, the optimal Smagorinsky coefficients for $\mathcal{M}_\text{DIF}$
    at the higher LES resolutions (where the CNN is unstable) go to zero, thus creating a divergence-free right hand side.
}

For the CNN, the initial model parameters $\theta_0$ are sampled from a uniform
distribution. They are improved by minimizing the stochastic loss function using
the ADAM optimizer~\cite{Kingma2017}. The gradients are obtained using reverse
mode automatic differentiation.

\begin{figure}
    \centering
    \includegraphics[width=1\textwidth]{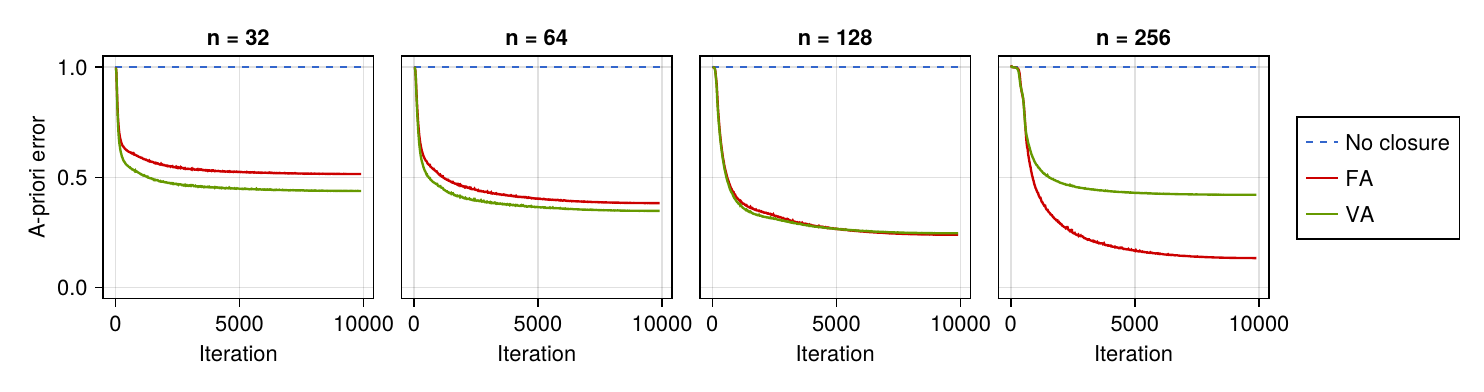}
    \caption{\revone{
        Relative a-priori error on validation set during a-priori training for
        $10^4$ iterations. From left to right: $\bar{n} = 32, 64, 128, 256$.
    }}
    \label{fig:priortraining}
\end{figure}

\begin{figure}
    \centering
    \includegraphics[width=1\textwidth]{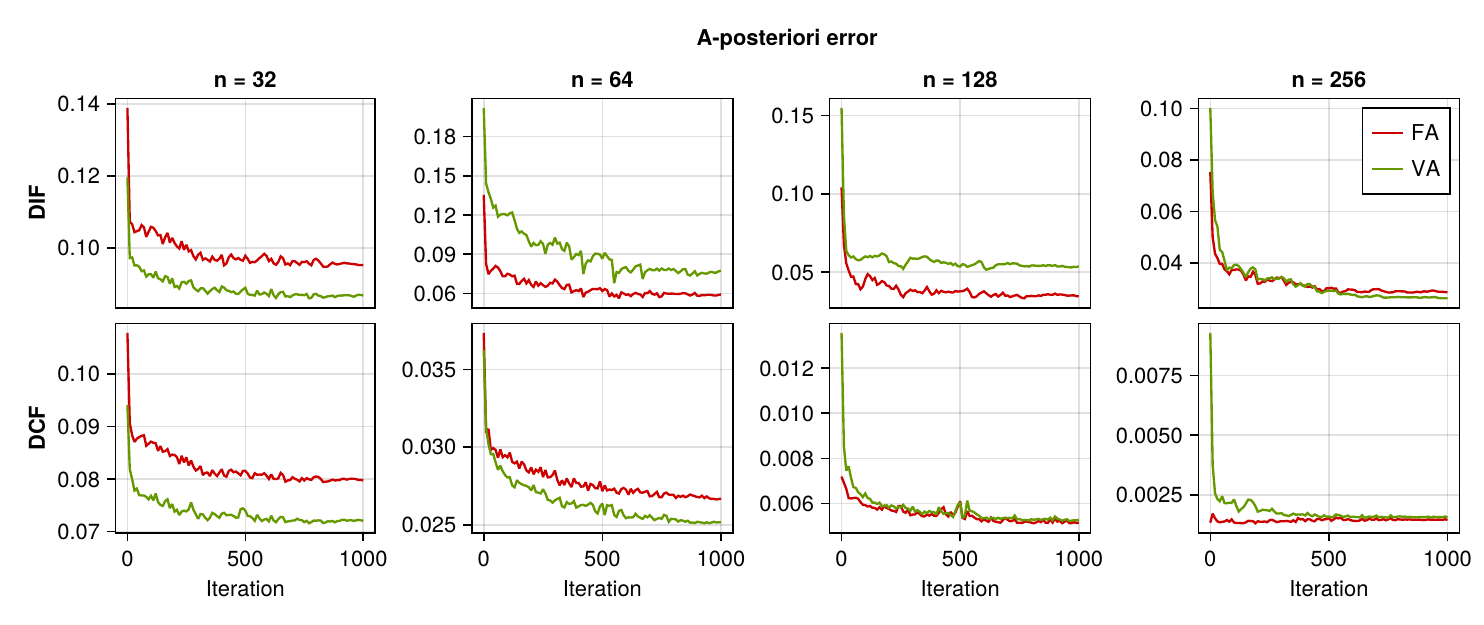}
    \caption{\revone{
        Relative a-posteriori error on validation set during a-posteriori
        training for $1000$ iterations. From left to right:
        $\bar{n} = 32, 64, 128, 256$.
    }}
    \label{fig:posttraining}
\end{figure}

We start by training using the a-priori loss function \eqref{eq:Lprior}.
We learn one set of parameters $\theta^\text{prior}$ for
each of the training grids and filter types. Each time, the model is trained for
\revone{
    
    $50$ epochs using the Adam optimizer~\cite{Kingma2017}
    with default weight decay and momentum parameters. 
    Each epoch consists of iterating through the 198 training batches, each
    containing 64 $(\bar{u}, c)$ snapshot pairs, resulting in roughly $10^4$
    iterations of stochastic gradient descent.
    The learning rate is set before each epoch using a cosine annealing
    scheduler~\cite{Loshchilov2017} with initial learning rate $10^{-3}$ and
    final learning rate $10^{-6}$.
    Every $20$ iterations the a-priori error is evaluated on the
    validation dataset.
    The validation error is shown in figure~\ref{fig:priortraining}. Note
    that the commutator error is different for each LES resolution and filter type,
    and so the relative a-priori error is not directly comparable between
    setups.}
The parameters giving the lowest validation error are retained after training.

Since the a-priori loss $L^\text{prior}$ does not take into account the effect
of the LES model we use, we
\revone{
    \R{finetuning}
    
    fine-tune the a-priori trained CNN parameters by training 
}
using the a-posteriori loss function $L^\text{post}_\mathcal{M}$.
This time, the LES model
$\mathcal{M} \in \{\mathcal{M}_{\text{DIF}}, \mathcal{M}_{\text{DCF}}\}$
is part of the loss function definition. We use
\revone{ $n_\text{unroll} = 50$}, and
train for 1000 iterations\revone{, as shown in figure~\ref{fig:posttraining}}.
\revone{
    Since the weights are already converged from the a-priori training,
    we use a smaller learning rate with the scheduler, starting at $10^{-4}$ and
    ending at $10^{-6}$.}
The best parameters from the
a-priori training ($\theta^\text{prior}$) are used as initial parameters for the
a-posteriori training. The Adam optimizer is reinitialized without the history
terms from the a-priori training session. The parameters
$\theta^\text{post}_\mathcal{M}$ giving the lowest a-posteriori error on the
validation dataset during training are retained after training.

\revone{
    For FA with $\bar{n} = 256$, the parameters found with a-priori training are
    already close to the optimum, and subsequent a-posteriori training only
    leads to minor improvements.
    For VA with $\bar{n} = 256$, the same a-priori training validation error is
    much higher than for FA. A-posteriori training therefore improves the validation
    error significantly, and eventually reaches the same error as for FA. This
    is likely due to the fact that for VA, the commutator error used in a-priori
    training is inconsistent with $\mathcal{M}_\text{DCF}$. Training with
    $L^\text{post}_\text{DCF}$ corrects for this inconsistency. For FA, the
    computed commutator errors are fully consistent with
    $\mathcal{M}_\text{DCF}$, and a-posteriori training is not needed.
}

\section{Divergence-free filter} \label{sec:divfree}

In this appendix we give the proof that the face-averaging filter is
divergence-free. The proof is a natural consequence of the continuous
divergence-theorem.

\subsection{Simple example} \label{sec:divfree_simple}

Consider a uniform 2D-grid with $3 \times 3$ fine volumes in each coarse volume.
The divergence-free constraint for a volume $\Omega_{i, j}$ of the fine grid reads:
\begin{equation}
    (D u)_{i, j} = (\delta_1 u^1)_{i, j} + (\delta_2 u^2)_{i, j} = 0,
\end{equation}
where
\begin{equation}
    (\delta_1 u^1)_{i, j} =
    \frac{u^1_{i + \frac{1}{2}, j} - u^1_{i - \frac{1}{2}, j}}{\Delta^1}
    \quad \text{and} \quad
    (\delta_2 u^2)_{i, j} =
    \frac{u^2_{i, j + \frac{1}{2}} - u^2_{i, j - \frac{1}{2}}}{\Delta^2}.
\end{equation}
Let $(a, b)$ and $(i, j)$ be coarse-grid and fine-grid indices such that
$\Omega_{i, j}$ is in the center of $\bar{\Omega}_{a, b}$. The four
face-averaged velocities at the boundary of $\bar{\Omega}_{a, b}$ are
\begin{equation}
\begin{split}
    \bar{u}^1_{a + \frac{1}{2}, b}
    & = \tfrac{1}{3} u^1_{i + \frac{3}{2}, j - 1}
    + \tfrac{1}{3} u^1_{i + \frac{3}{2}, j}
    + \tfrac{1}{3} u^1_{i + \frac{3}{2}, j + 1} \\
    \bar{u}^1_{a - \frac{1}{2}, b} & = \tfrac{1}{3} u^1_{i - \frac{3}{2}, j - 1}
    + \tfrac{1}{3} u^1_{i - \frac{3}{2}, j}
    + \tfrac{1}{3} u^1_{i - \frac{3}{2}, j + 1} \\
    \bar{u}^2_{a, b + \frac{1}{2}} & = \tfrac{1}{3} u^2_{i - 1, j + \frac{3}{2}}
    + \tfrac{1}{3} u^2_{i, j + \frac{3}{2}}
    + \tfrac{1}{3} u^2_{i + 1, j + \frac{3}{2}} \\
    \bar{u}^2_{a, b - \frac{1}{2}} & = \tfrac{1}{3} u^2_{i - 1, j - \frac{3}{2}}
    + \tfrac{1}{3} u^2_{i, j - \frac{3}{2}}
    + \tfrac{1}{3} u^2_{i + 1, j - \frac{3}{2}}.
\end{split}
\end{equation}
The coarse-grid divergence reads
\begin{equation}
    (\bar{D} \bar{u})_{a, b} =
    (\bar{\delta}_1 \bar{u}^1)_{a, b} +
    (\bar{\delta}_2 \bar{u}^2)_{a, b}.
\end{equation}
For the first term, we get
\begin{equation}
\begin{split}
    (\bar{\delta}_1 \bar{u}^1)_{a, b}
    & = \frac{\bar{u}^1_{a + \frac{1}{2}, b}
    - \bar{u}^1_{a - \frac{1}{2}, b}}{\bar{\Delta}^1} \\
    & = \tfrac{1}{3} \frac{u^1_{i + \frac{3}{2}, j - 1}
    - u^1_{i - \frac{3}{2}, j - 1}}{3 \Delta^1}
    + \tfrac{1}{3} \frac{u^1_{i + \frac{3}{2}, j}
    - u^1_{i - \frac{3}{2}, j}}{3 \Delta^1}
    + \tfrac{1}{3} \frac{u^1_{i + \frac{3}{2}, j + 1}
    - u^1_{i - \frac{3}{2}, j + 1}}{3 \Delta^1} \\
    & = \frac{1}{9 \Delta^1} \left( u^1_{i + \frac{3}{2}, j - 1}
    - (u^1_{i + \frac{1}{2}, j - 1} - u^1_{i + \frac{1}{2}, j - 1})
    - (u^1_{i - \frac{1}{2}, j - 1} - u^1_{i - \frac{1}{2}, j - 1})
    - u^1_{i - \frac{3}{2}, j - 1} \right) \\
    & + \frac{1}{9 \Delta^1} \left( u^1_{i + \frac{3}{2}, j}
    - (u^1_{i + \frac{1}{2}, j} - u^1_{i + \frac{1}{2}, j})
    - (u^1_{i - \frac{1}{2}, j} - u^1_{i - \frac{1}{2}, j})
    - u^1_{i - \frac{3}{2}, j} \right) \\
    & + \frac{1}{9 \Delta^1} \left( u^1_{i + \frac{3}{2}, j + 1}
    - (u^1_{i + \frac{1}{2}, j + 1} - u^1_{i + \frac{1}{2}, j + 1})
    - (u^1_{i - \frac{1}{2}, j + 1} - u^1_{i - \frac{1}{2}, j + 1})
    - u^1_{i - \frac{3}{2}, j + 1} \right) \\
    & = \frac{1}{9} \Big( (\delta_1 u^1)_{i + 1, j - 1}
    + (\delta_1 u^1)_{i, j - 1}
    + (\delta_1 u^1)_{i - 1, j - 1} \\
    & + (\delta_1 u^1)_{i + 1, j}
    + (\delta_1 u^1)_{i, j}
    + (\delta_1 u^1)_{i - 1, j} \\
    & + (\delta_1 u^1)_{i + 1, j + 1}
    + (\delta_1 u^1)_{i, j + 1}
    + (\delta_1 u^1)_{i - 1, j + 1} \Big).
\end{split}
\end{equation}
With a similar derivation for $\bar{\delta}_2 \bar{u}^2$, we get the following
expression for the coarse grid divergence:
\begin{equation}
\begin{split}
    (\bar{D} \bar{u})_{a, b} & = \frac{1}{9} \Big(
    (D u)_{i - 1, j - 1} + (D u)_{i, j - 1} + (D u)_{i + 1, j - 1} \\
    & + (D u)_{i - 1, j} + (D u)_{i, j} + (D u)_{i + 1, j} \\
    & + (D u)_{i - 1, j + 1} + (D u)_{i, j + 1} + (D u)_{i + 1, j + 1} \Big) \\
    & = 0.
\end{split}
\end{equation}
The filtered velocity field is indeed divergence-free.

It is also interesting to point out that $\bar{D} \bar{u}$ can be seen as a
volume average of $D u$, i.e.\ $\bar{D} \bar{u} = \Psi D u$, or $\bar{D} \Phi =
\Psi D$, for a certain pressure filter $\Psi$ built with uniform $3 \times
3$-stencils of weights $1 / 9$. In other words: the face-averaging velocity
filter goes hand in hand with a volume-averaging pressure filter.

\subsection{Proof for the general case}

A similar proof can be shown for a general non-uniform grid in 2D or 3D.
Consider a coarse grid index $J$. The fine grid volumes $\Omega_I$ contained
inside $\bar{\Omega}_J$ are indexed by $I \in \mathcal{K}_J = \{ I \ | \
\Omega_I \subset \bar{\Omega}_J\}$. The face-averaging filter is defined by
\begin{equation}
    \bar{u}^\alpha_J = \sum_{I \in \mathcal{F}^\alpha_J} \rho^\alpha_{J, I}
    u^\alpha_I,
\end{equation}
where $\mathcal{F}^\alpha_J = \{I \ | \ \Gamma^\alpha_I \in
\bar{\Gamma}^\alpha_J \}$ contains the face-indices and $\rho^\alpha_{J, I}$ are
weights to be determined. We assume that $\rho^\alpha_{J, I}$ is independent of
$I_\alpha$ and $J_\alpha$. The fine grid divergence is given by
\begin{equation}
    (D u)_I =
    \sum_{\alpha = 1}^d (\delta_\alpha u^\alpha)_I = 
    \sum_{\alpha = 1}^d \frac{u^\alpha_{I + h_\alpha} - u^\alpha_{I -
    h_\alpha}}{| \Delta^\alpha_{I_\alpha} |}
    = 0.
\end{equation}
The coarse grid divergence is given by
\begin{equation} \label{eq:divfreefilter}
    \begin{split}
        (\bar{D} \bar{u})_J
        & = \sum_{\alpha = 1}^d (\bar{\delta}_\alpha \bar{u}^\alpha)_J \\
        & = \sum_{\alpha = 1}^d \frac{\bar{u}^\alpha_{J + h_\alpha} - \bar{u}^\alpha_{J -
        h_\alpha}}{| \bar{\Delta}^\alpha_{J_\alpha} |}
        \quad \revboth{\text{by definition of finite difference operator}} \\
        & = \sum_{\alpha = 1}^d 
        \frac{1}{| \bar{\Delta}^\alpha_{J_\alpha} |}
        \left(
            \sum_{I \in \mathcal{F}^\alpha_{J + h_\alpha}}
            \rho^\alpha_{J, I}
            u^\alpha_I
            - \sum_{I \in \mathcal{F}^\alpha_{J - h_\alpha}}
            \rho^\alpha_{J, I}
            u^\alpha_I
        \right)
        \quad \text{by definition of $\bar{u}^\alpha_J$} \\
        & =
        \sum_{\alpha = 1}^d 
        \frac{1}{| \bar{\Delta}^\alpha_{J_\alpha} |}
        \sum_{I \in \mathcal{K}_J}
        \left(
            \rho^\alpha_{J, I + h_\alpha} u^\alpha_{J, I + h_\alpha}
            - \rho^\alpha_{J, I - h_\alpha} u^\alpha_{J, I - h_\alpha}
        \right)
        \quad \text{telescoping sum over $I_\alpha$} \\
        & =
        \sum_{\alpha = 1}^d 
        \frac{1}{| \bar{\Delta}^\alpha_{J_\alpha} |}
        \sum_{I \in \mathcal{K}_J}
        \rho^\alpha_{J, I} 
        \left( u^\alpha_{I + h_\alpha} - u^\alpha_{I - h_\alpha} \right)
        \quad \text{since $\rho^\alpha_{J, I}$ is independent of $I_\alpha$} \\
        & =
        \sum_{I \in \mathcal{K}_J}
        \frac{| \Omega_I |}{| \bar{\Omega}_J |}
        \sum_{\alpha = 1}^d 
        \frac{| \bar{\Gamma}^\alpha_J |}{| \Gamma^\alpha_{I} |}
        \rho^\alpha_{J, I}
        \frac{u^\alpha_{I + h_\alpha} - u^\alpha_{I - h_\alpha}}{| \Delta^\alpha_{I_\alpha} |}
        \quad \text{rewrite terms with} \ | \Omega_I | = | \Gamma^\alpha_I | | \Delta^\alpha_I | 
        \ \forall \alpha \\
        & = 
        \sum_{I \in \mathcal{K}_J}
        \frac{| \Omega_I |}{| \bar{\Omega}_J |}
        (D u)_I
        \quad \text{if we choose} 
        \ \rho^\alpha_{J, I} =  | \Gamma^\alpha_{I} | / | \bar{\Gamma}^\alpha_{J} | \\
        & = 0
        \quad \revboth{\text{since} \ (D u)_I = 0.}
    \end{split}
\end{equation}
The chosen $\rho^\alpha_{J, I}$ is indeed independent of $I_\alpha$ and
$J_\alpha$. We also get the property $\sum_{I \in \mathcal{F}^\alpha_J}
\rho^\alpha_{J, I} = 1$, so constant fine grid velocities are preserved upon
filtering. In other words, choosing $\bar{u}^\alpha_J$ as a weighted average of
the DNS-velocities passing through the coarse volume face
$\bar{\Gamma}^\alpha_J$ gives a divergence-free $\bar{u}$. Note that the
divergence constraint only holds for the face-size weights chosen above. Using
arbitrary weights such as Gaussian weights would not work.

We also observe that for the general case, just like for the simple case in
\ref{sec:divfree_simple}, we can write $\bar{D} \Phi = \Psi D$ for a certain
pressure filter $\Psi$. It is defined by
\begin{equation}
    (\Psi p)_J =
    \sum_{I \in \mathcal{K}_J}
    \frac{| \Omega_I |}{| \bar{\Omega}_J |}
    p_I
\end{equation}
for all fields $p$ defined in the pressure points. If we consider the $3 \times
3$ compression from \ref{sec:divfree_simple}, with $J = (a, b)$, we
effectively get $\mathcal{K}_{a, b} = \{i - 1, i, i + 1\} \times \{j - 1, j, j +
1\}$ and $| \Omega_I | / | \bar{\Omega}_J | = 1 / 9$ for all $I \in
\mathcal{K}_J$.

\section{Discretize-then-filter (without differentiating the constraint)} \label{sec:discretize_filter}

In this appendix we show a problem that arises when discretely filtering the differential-algebraic system
\eqref{eq:mass}-\eqref{eq:momentum}, instead of filtering the ``pressure-free'' equation \eqref{eq:dns}. 
Since the discrete DNS system \eqref{eq:mass}-\eqref{eq:momentum} includes a divergence term and pressure term, we need to define a pressure-filter $\Psi$ (in addition to the velocity filter $\Phi$), such that $\bar{p} = \Psi p$. This results in the following set of equations:
\begin{align}
    \Psi D u & = 0, \\
    \frac{\mathrm{d} \bar{u}}{\mathrm{d} t} & = \Phi F(u) - \Phi G p.
\end{align}
These can be rewritten as follows:
\begin{align}
    \bar{D} \bar{u} & = c_{D}(u,\bar{u}),  \\
    \frac{\mathrm{d} \bar{u}}{\mathrm{d} t} & = \bar{F}(\bar{u}) + \tilde{c}(u,\bar{u}) - \bar{G} \bar{p} + c_{P}(p,\bar{p}) .
\end{align}
Here $c_{D}(u,\bar{u}) := (\bar{D} \Phi - \Psi D) u$ represents the commutator error between the discrete divergence operator and filtering, $\tilde{c} (u,\bar{u}) := (\Phi F(u) - \bar{F}(\bar{u}))$ represents the commutator error arising from filtering $F(u)$ (note that it is different from the one used in equation \eqref{eq:commutator_error}), and the commutator error for the pressure $c_{P} (p,\bar{p}) = (\bar{G} \Psi - \Phi G) p$. 
In case of a face-averaging filter, one has $c_{D} = 0$, which is the constraint that needs to be enforced by the filtered pressure, when moving to the LES equations. However, in the above formulation an additional commutator error for the pressure appears, which is unwanted and it is unclear how it should be modelled. In the discretize-differentiate-filter approach, this issue with the pressure is circumvented.

% \revtwo{
\section{Continuous filters and transfer functions} \label{sec:continuous_filtering}

The discrete volume-averaging filter $\Phi^\text{VA}$ is an
approximation to the continuous top-hat volume-averaging filter $g$ defined by
\begin{equation}
    \begin{split}
        (g * \varphi)(x)
        & = \int_{\mathbb{R}^d} g(y) \varphi(x - y) \, \mathrm{d} \Omega(y) \\
        & = \frac{1}{\bar{\Delta}^d}
        \int_{\left[ -\frac{\bar{\Delta}}{2}, \frac{\bar{\Delta}}{2} \right]^d}
        \varphi(x - y) \, \mathrm{d} \Omega(y),
    \end{split}
\end{equation}
where $\varphi$ is a scalar field,
\begin{equation*}
    g(x) = \frac{1}{\bar{\Delta}^d}
        \prod_{\alpha = 1}^d
        I \left(x^\alpha \in
        \left[ -\frac{\bar{\Delta}}{2}, \frac{\bar{\Delta}}{2} \right]
        \right)
\end{equation*}
is the convolutional filter kernel, and $I(a \in A)$ is an indicator function equal
to $1$ if $a \in A$ and $0$ otherwise.
The continuous filter $g$ can also be interpreted in terms of its transfer
function $G$, defined by
\begin{equation}
    G_k = \prod_{\alpha = 1}^d \frac{\sin(\pi k_\alpha \bar{\Delta} /
    2)}{\pi k_\alpha \bar{\Delta} / 2}
\end{equation}
at a given wavenumber vector $k = (k_1, \dots, k_d)$.
Then the Fourier transform of $g * \varphi$ is given by
$\widehat{(g * \varphi)}_k = G_k \hat{\varphi}_k$.

Similarly, we can define the continuous top-hat face-averaging filter (with face normal
to the direction $\alpha$) as
\begin{equation}
    (g^\alpha * \varphi)(x) = \int_{\mathbb{R}^d} g^\alpha(y) \varphi(x - y) \,
    \mathrm{d} \Omega(y),
\end{equation}
where the kernel
\begin{equation}
    g^\alpha(x) = \frac{1}{\bar{\Delta}^{d - 1}} \delta(x^\alpha)
    \prod_{\beta \neq \alpha}
    I \left(x^\beta \in
    \left[ -\frac{\bar{\Delta}}{2}, \frac{\bar{\Delta}}{2} \right]
    \right)
\end{equation}
is the same as for the volume-averaging filter, but with one interval indicator
replaced by a Dirac delta function $\delta$.
The transfer function of the face-averaging filter is given by
\begin{equation}
    G^\alpha_k = \prod_{\beta \neq \alpha} \frac{\sin(\pi k_\beta \bar{\Delta} / 2)}{\pi k_\beta \bar{\Delta} / 2}.
\end{equation}
In other words, the $\alpha$-face-averaging filter is performing top-hat averaging
in every direction except for direction $\alpha$, where it leaves signals intact.

\begin{figure}
    \centering
    \includegraphics[width=0.9\textwidth]{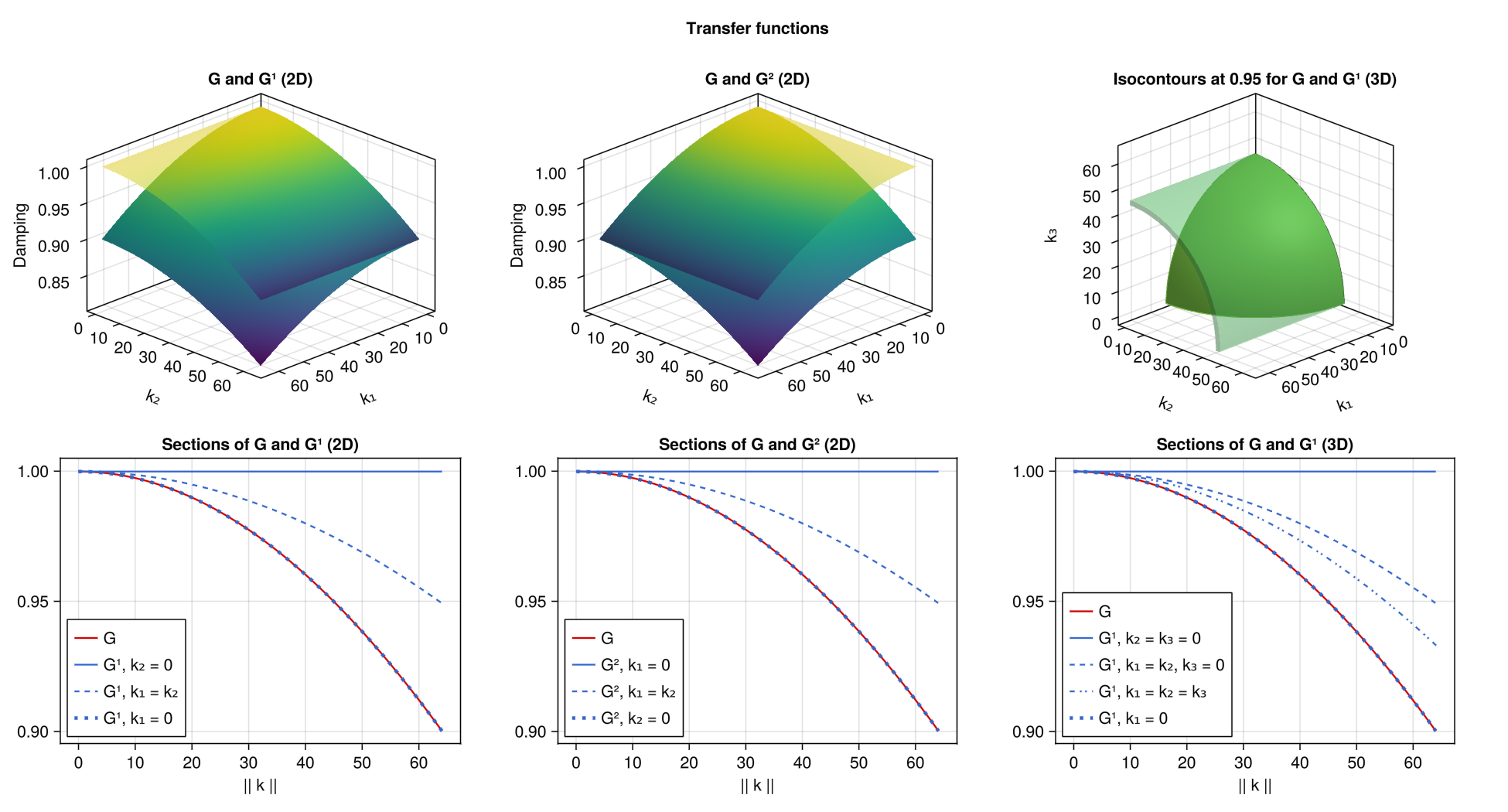} 
    \caption{\revtwo{
        Transfer functions in 2D and 3D. The two top left plots show the 2D
        transfer functions of the volume-averaging and face-averaging filters.
        The upper surface is $G^\alpha$, the lower one is $G$.
        The top-right plot shows isocontours of the 3D transfer functions $G$
        (behind) and $G^1$ (transparent, in front) at a given damping level
        $0.95$. The bottom row shows sections of the transfer functions in
        different wavenumber directions $k / \| k \|$ as a function of the
        wavenumber magnitude $\| k \|$. $G$ does not depend on the direction.
    }}
    \label{fig:transferfunctions}
\end{figure}

The transfer functions $G$, $G^1$, and $G^2$ are shown in
figure~\ref{fig:transferfunctions}.
The filter width is set to $\bar{\Delta} = 1 / 128$, corresponding to a cut-off
wavenumber level at $\| k \|_\infty = 64$. We also show the transfer
functions $G$ and $G^1$ for the 3D versions of the filters.
In 3D, the isocontour of $G$ at a given damping level takes the shape of a
sphere, as $G_k$ only depends on $\| k \|$. The isocontour of $G^\alpha$ takes
the shape of a cylinder with axis along the $\alpha$-direction, as the value of
$G^\alpha_k$ does not depend on $k_\alpha$.

We can clearly see that for \emph{all} wavenumbers $k$,
the volume-averaging filter is damping more than the face-averaging filter
($\forall k, G_k \leq G^\alpha_k$).
For $k_\alpha = 0$, we get $G_k = G^\alpha_k$ (VA and FA coincide),
and for $k_\alpha > 0$, we get the strict inequality $G_k < G^\alpha_k$.
For all $k$ such $k_\beta = 0$ for $\beta \neq \alpha$, we get
$G^\alpha_k = 1$ (FA does not filter in the normal direction).

% }

% \revboth{
\section{Discrete and continuous commutator errors for the 2D Taylor-Green Vortex}
\label{sec:commutator}

In this section, we compare the discrete and continuous convective commutator
errors. For this purpose we employ an analytical solution to the Navier-Stokes
equations on a periodic domain $\Omega = [0, 2 \pi]^2$, being the Taylor-Green
vortex~\cite{Green1937}:
\begin{equation} \label{eq:taylorgreen}
    \begin{split}
        u(x, y, t) & = - \sin(x) \cos(y) \mathrm{e}^{-2 \nu t}, \\
        v(x, y, t) & = \phantom{+} \cos(x) \sin(y) \mathrm{e}^{-2 \nu t}, \\
        p(x, y, t) & = \frac{1}{4} \left( \cos(2 x) + \cos(2 y) \right)
        \mathrm{e}^{-4 \nu t}.
    \end{split}
\end{equation}
In the remainder of this appendix, $u$, $v$, and $p$ will refer to this solution
only, and the results that follow are specific to this solution.
We will also drop the time dependence and consider $u(x, y) = u(x, y, 0)$ and similarly for $v$ and $p$.

\subsection{Top-hat filter}

First, we note that
\begin{equation}
    % \begin{split}
    %     \int_{x - \Delta / 2}^{x + \Delta / 2} \sin(x') \, \mathrm{d} x'
    %     & = -\cos(x + \Delta / 2) + \cos(x - \Delta / 2) \\
    %     & = - (\cos(x) \cos(\Delta / 2) - \sin(x) \sin(\Delta / 2) \\
    %     & \phantom{=} + (\cos(x) \cos(\Delta / 2) + \sin(x) \sin(\Delta / 2)) \\
    %     & = 2 \sin(x) \sin(\Delta / 2), \\
    %     \int_{x - \Delta / 2}^{x + \Delta / 2} \cos(x') \, \mathrm{d} x'
    %     & = \sin(x + \Delta/2) - \sin(x - \Delta/2) \\
    %     & = (\sin(x) \cos(\Delta / 2) + \cos(x) \sin(\Delta / 2)) \\
    %     & - (\sin(x) \cos(\Delta / 2) - \cos(x) \sin(\Delta / 2)) \\
    %     & = 2 \cos(x) \sin(\Delta / 2).
    % \end{split}
    \int_{x - \Delta / 2}^{x + \Delta / 2} \sin(x') \, \mathrm{d} x'
    = 2 \sin(x) \sin(\Delta / 2), \quad
    \int_{x - \Delta / 2}^{x + \Delta / 2} \cos(x') \, \mathrm{d} x'
    = 2 \cos(x) \sin(\Delta / 2).
\end{equation}

Applying a volume-averaging top-hat filter
$\varphi \mapsto \bar{\varphi}$
of width $\Delta$
with
\begin{equation}
    \bar{\varphi}(x, y) = \frac{1}{\Delta^2} \int_{x - \Delta/2}^{x + \Delta/2} \int_{y - \Delta/2}^{y + \Delta/2}
    \varphi(x', y') \, \mathrm{d} x' \mathrm{d} y'
\end{equation}
gives the following filtered solution:
\begin{equation}
    % \begin{split}
    %     \bar{u}(x, y)
    %     & = \frac{1}{\Delta^2} \int_{x - \Delta/2}^{x + \Delta/2} \int_{y - \Delta/2}^{y + \Delta/2}
    %     u(x', y', t) \, \mathrm{d} x' \mathrm{d} y', \\
    %     & = - \frac{1}{\Delta^2}
    %     \int_{x - \Delta/2}^{x + \Delta/2} \sin(x') \, \mathrm{d} x'
    %     \int_{y - \Delta/2}^{y + \Delta/2} \cos(y') \, \mathrm{d} y' \\
    %     & = - \frac{4}{\Delta^2} \sin^2(\Delta / 2) \sin(x) \cos(y) \\
    %     & = \sinc^2(\Delta / 2) u(x, y), \\
    % \end{split}
    \bar{u}(x, y)
    = - \frac{1}{\Delta^2}
    \int_{x - \Delta/2}^{x + \Delta/2} \sin(x') \, \mathrm{d} x'
    \int_{y - \Delta/2}^{y + \Delta/2} \cos(y') \, \mathrm{d} y'
    = \sinc^2(\Delta / 2) u(x, y),
\end{equation}
where $\sinc(x) = \sin(x) / x$ is the transfer function of the 1D top-hat
filter. A similar derivation for $\bar{v}$ and $\bar{p}$ gives
\begin{equation}
    \bar{u} = \sinc^2(\Delta / 2) u, \quad
    \bar{v} = \sinc^2(\Delta / 2) v, \quad
    \bar{p} = \sinc(\Delta) p.
\end{equation}

\subsection{Continuous convective commutator error}

For the continuous commutator error in the Navier-Stokes equations,
we need the nonlinearities
$\bar{u} \bar{u}$, $\bar{u} \bar{v}$, $\bar{v} \bar{u}$, and $\bar{v} \bar{v}$.
They are given by
\begin{equation}
    \begin{split}
        \bar{u} \bar{u}
        % & = \sinc^4(\Delta / 2) \sin^2(x) \cos^2(y) \mathrm{e}^{-4 \nu t}, \\
        % & = \frac{1}{4} \sinc^4(\Delta / 2) (1 - \cos(2 x)) (1 + \cos(2 y)) \mathrm{e}^{-4 \nu t}, \\
        & = \frac{1}{4} \sinc^4(\Delta / 2) (1 - \cos(2 x) + \cos(2 y) - \cos(2 x) \cos(2 y)) , \\
        \bar{v} \bar{v}
        & = \frac{1}{4} \sinc^4(\Delta / 2) (1 + \cos(2 x) - \cos(2 y) - \cos(2 x) \cos(2 y)) , \\
        \bar{u} \bar{v}
        & = - \frac{1}{4} \sinc^4(\Delta / 2) \sin(2 x) \sin(2 y) . \\
    \end{split}
\end{equation}
Next, we are interested in the quantities
$\overline{u u}$, $\overline{u v}$, $\overline{v u}$, and $\overline{v v}$.
Integrating $u u$, $v v$, and $u v$ over $x$ and $y$ separately gives
\begin{equation}
    \begin{split}
        \overline{u u}
        % & = \frac{1}{4} (1 - \sinc(\Delta) \cos(2 x)) (1 + \sinc(\Delta) \cos(2 y)) -4 \nu t} \\
        & = \frac{1}{4} (1 - \sinc(\Delta)(\cos(2 x) - \cos(2 y)) - \sinc^2(\Delta) \cos(2x) \cos(2 y)) , \\
        \overline{v v}
        % & = \frac{1}{4} (1 + \sinc(\Delta) \cos(2 x)) (1 - \sinc(\Delta) \cos(2 y)) -4 \nu t} \\
        & = \frac{1}{4} (1 + \sinc(\Delta)(\cos(2 x) - \cos(2 y)) - \sinc^2(\Delta) \cos(2x) \cos(2 y)) , \\
        \overline{u v}
        & = - \frac{1}{4} \sinc^2(\Delta) \sin(2 x) \sin(2 y) . \\
    \end{split}
\end{equation}
Finally, we can assemble the sub-filter stress tensor
$\tau$:
\begin{equation}
    \begin{split}
        % tau_xx
        \tau_{x x}
        &= \overline{u u} - \bar{u} \bar{u} \\
        % & = \frac{1}{4}
        % \left( 1 - \sinc(\Delta) \cos(2 x) + \sinc(\Delta) \cos(2 y) - \sinc^2(\Delta) \cos(2 x) \cos(2 y) \right)
        % -4 \nu t}, \\
        % & - \frac{1}{4} \sinc^4(\Delta / 2)
        % \left( 1 - \cos(2 x) + \cos(2 y) - \cos(2 x) \cos(2 y) \right)
        % -4 \nu t}, \\
        & = \frac{1}{4} \Big(
        1 - \sinc^4(\Delta / 2)
        + (\sinc^4(\Delta / 2) - \sinc(\Delta)) (\cos(2 x) - \cos(2 y)) \\
        & + (\sinc^4(\Delta / 2) - \sinc^2(\Delta)) \cos(2 x) \cos(2 y) \Big)
        , \\
        % % tau_yy
        % \tau_{y y}
        % & = \overline{v v} - \bar{v} \bar{v} \\
        % & = \frac{1}{4} \Big(
        % 1 - \sinc^4(\Delta / 2)
        % - (\sinc^4(\Delta / 2) - \sinc(\Delta)) (\cos(2 x) - \cos(2 y)) \\
        % & + (\sinc^4(\Delta / 2) - \sinc^2(\Delta)) \cos(2 x) \cos(2 y) \Big)
        % , \\
        % tau_xy
        \tau_{x y}
        & = \overline{u v} - \bar{u} \bar{v} \\
        & =
        \frac{1}{4} (\sinc^4(\Delta / 2) - \sinc^2(\Delta))
        \sin(2 x) \sin(2 y), \\
    \end{split}
\end{equation}
and similarly for $\tau_{y y}$.
% The terms in $\nabla \cdot \tau$ then follow as
% \begin{equation}
%     \begin{split}
%         \frac{\partial \tau_{x x}}{\partial x}
%         & =
%         - \frac{1}{2} (\sinc^4(\Delta / 2) - \sinc^2(\Delta)) \sin(2 x) \cos(2 y)
%         -\frac{1}{2} (\sinc^4(\Delta / 2) - \sinc(\Delta)) \sin(2 x)
%         , \\
%         \frac{\partial \tau_{x y}}{\partial y}
%         & = \frac{1}{2} (\sinc^4(\Delta / 2) - \sinc^2(\Delta)) \sin(2 x) \cos(2 y) , \\
%         \frac{\partial \tau_{y x}}{\partial x}
%         & = \frac{1}{2} (\sinc^4(\Delta / 2) - \sinc^2(\Delta)) \cos(2 x) \sin(2 y) , \\
%         \frac{\partial \tau_{y y}}{\partial y}
%         & =
%         - \frac{1}{2} (\sinc^4(\Delta / 2) - \sinc^2(\Delta)) \cos(2 x) \sin(2 y)
%         - \frac{1}{2} (\sinc^4(\Delta / 2) - \sinc(\Delta)) \sin(2 y)
%         , \\
%     \end{split}
% \end{equation}
Finally, the sub-filter force terms
$c_x =
\frac{\partial \tau_{x x}}{\partial x} +
\frac{\partial \tau_{x y}}{\partial y}$
and
$c_y =
\frac{\partial \tau_{y x}}{\partial x} +
\frac{\partial \tau_{y y}}{\partial y}$
are given by
\begin{empheq}[box=\boxstyle]{equation} \label{eq:commutator_continuous}
    \begin{split}
        c_x & =
        % \frac{\partial \tau_{x x}}{\partial x} +
        % \frac{\partial \tau_{x y}}{\partial y} =
        -\frac{1}{2} (\sinc^4(\Delta / 2) - \sinc(\Delta)) \sin(2 x), \\
        % = -(\sinc^4(\Delta / 2) - \sinc(\Delta))
        % \left( \frac{\partial}{\partial x} (u u) + \frac{\partial}{\partial y} (u v) \right), \\
        c_y & =
        % \frac{\partial \tau_{y x}}{\partial x} +
        % \frac{\partial \tau_{y y}}{\partial y} =
        -\frac{1}{2} (\sinc^4(\Delta / 2) - \sinc(\Delta)) \sin(2 y). \\
        % = -(\sinc^4(\Delta / 2) - \sinc(\Delta))
        % \left( \frac{\partial}{\partial x} (v u) + \frac{\partial}{\partial y} (v v) \right), \\
    \end{split}
\end{empheq}
Since the filter width is small compared to the domain size ($2 \pi$), we can assume that
high powers of $\Delta$ may be negligible. Using the Taylor series expansions
\begin{equation} \label{eq:taylor_sinc}
    \begin{split}
        \sinc(x) & = 1 - \frac{1}{6} x^2 + \frac{1}{120} x^4 + \mathcal{O}(x^6), \\
        % \sinc^2(x) & = 1 - \frac{1}{3} x^2 + \frac{2}{45} x^4 + \mathcal{O}(x^6), \\
        \sinc^4(x) & = 1 - \frac{2}{3} x^2 + \frac{1}{5} x^4 + \mathcal{O}(x^6), \\
    \end{split}
\end{equation}
we get the coefficient
\begin{equation}
    \begin{split}
        \sinc^4(\Delta / 2) - \sinc(\Delta)
        % & = \frac{1}{80} \Delta^4 - \frac{1}{120} \Delta^4 + \mathcal{O}(\Delta^6), \\
        & = \frac{1}{240} \Delta^4 + \mathcal{O}(\Delta^6). \\
        % \sinc^4(\Delta / 2) - \sinc^2(\Delta)
        % % & = \frac{1}{6} \Delta^2 + \frac{1}{80} \Delta^4 - \frac{2}{45} \Delta^4 + \mathcal{O}(\Delta^6), \\
        % & = \frac{1}{6} \Delta^2 - \frac{23}{720} \Delta^4 + \mathcal{O}(\Delta^6).
    \end{split}
\end{equation}
If we drop the terms in $\mathcal{O}(\Delta^6)$ and higher, we get the following simplified expressions for the continuous commutator error:
\begin{equation}
    \begin{split}
        c_x & = -\frac{\Delta^4}{480} \sin(2 x)  + \mathcal{O}(\Delta^6), \\
        c_y & = -\frac{\Delta^4}{480} \sin(2 y)  + \mathcal{O}(\Delta^6). \\
    \end{split}
\end{equation}

\subsection{Discrete convective commutator error}

Next, we consider the discrete commutator error computed
from a discrete DNS. Assume for simplicity that the DNS field is the
exact solution sampled at the DNS grid points. Consider a grid point $(x, y)$.
Let $d$ denote the DNS grid spacing.

\subsubsection{Discrete filter}

Assume that $\Delta / 2 = n d$ for some integer $n$. A simple discretization of the
top-hat volume-averaging filter at a coarse grid point $(x, y)$ is given
by
\begin{equation}
    \begin{split}
        (\Phi u)(x, y)
        & =
        \frac{1}{(2 n + 1)^2}
        \sum_{i = -n}^{n}
        \sum_{j = -n}^{n}
        u(x + i d, y + j d) \\
        & =
        - \frac{1}{(2 n + 1)^2}
        \sum_{i = -n}^{n}
        \sin(x + i d)
        \sum_{j = -n}^{n}
        \cos(y + j d), \\
    \end{split}
\end{equation}
and similarly for $\Phi v$ and $\Phi p$.
For the $x$-part, we get
\begin{equation}
    \sum_{i = -n}^{n} \sin(x + i d)
    = \sum_{i = -n}^{n} (\sin(x) \cos(i d) + \cos(x) \sin(i d))
    = \left(\sum_{i = -n}^{n} \cos(i d) \right) \sin(x).
\end{equation}
% Similarly, for the $y$-part, we get
% \begin{equation}
%     \begin{split}
%         \sum_{j = -n}^{n} \cos(y + j d)
%         & = \sum_{j = -n}^{n} (\cos(y) \cos(j d) - \sin(y) \sin(j d)) \\
%         & = \left(\sum_{j = -n}^{n} \cos(j d) \right) \cos(y). \\
%     \end{split}
% \end{equation}
A similar derivation for the $y$-part, $v$ and $p$ give
\begin{equation}
    \Phi u = G_{n, d}^2 u, \quad
    \Phi v = G_{n, d}^2 v, \quad
    \Phi p = G_{n, 2d} p,
\end{equation}
where $G_{n, d} = \frac{1}{2 n + 1} \sum_{i = -n}^n \cos(i d)$
is the transfer function of the discrete top-hat filter in 1D.

% \subsubsection{Diffusion}
% 
% The finite difference diffusion operator
% $\operatorname{diff}_{x x, d} \approx \frac{\partial^2 }{\partial x^2}$
% is defined as
% 
% \begin{equation}
%     \begin{split}
%         \operatorname{diff}_{x x, d}(u)(x, y)
%         & = \frac{1}{d^2} (u(x + d, y) - 2 u(x, y) + u(x - d, y)) \\
%         & = - \frac{1}{d^2} (\sin(x + d) - 2 \sin(x) + \sin(x - d)) \cos(y) \\
%         & = + \frac{2}{d^2} (1 - \cos(d)) \sin(x) \cos(y) \\
%         & = - \frac{2}{d^2} (1 - \cos(d)) u(x, y), \\
%     \end{split}
% \end{equation}
% and similarly for $\operatorname{diff}_{y y, d}(u)(x, y)$. We thus get
% \begin{equation}
%     \operatorname{diff}_d(u)(x, y) =
%     \operatorname{diff}_{x x, d}(u)(x, y) +
%     \operatorname{diff}_{y y, d}(u)(x, y) = - \frac{4}{d^2} (1 - \cos(d)) u(x, y).
% \end{equation}
% For comparison, we have
% \begin{equation}
%     \nabla^2 u(x, y) = \frac{\partial^2 u}{\partial x^2}(x, y) +
%     \frac{\partial^2 u}{\partial y^2}(x, y) = -2 u(x, y),
% \end{equation}
% since 
% \begin{equation}
%     \frac{2}{d^2} (1 - \cos(d)) = \frac{2}{d^2} \left(1 - 1 + \frac{d^2}{2} -
%     \frac{d^4}{24} + \mathcal{O}(d^6) \right) = 1 - \frac{d^2}{12} + \mathcal{O}(d^4).
% \end{equation}
% 
% Since both the discrete filter and the discrete diffusion operator scale the
% velocity field by constants, we can assert that
% $\Phi_x \operatorname{diff}(u) = \operatorname{diff}(\Phi_x u)$, and similarly
% for $v$.

\subsubsection{Convective term}

We now compute the discrete convective term
$\operatorname{conv}_{x, d}(u, v)$ which contributes to the equation for $u$.
This term is required at a $u$-velocity grid point $(x, y)$, and we thus need
the quantity $u u$ at $(x + d / 2, y)$ and $(x - d / 2, y)$, and the quantity
$u v$ at $(x, y + d / 2)$ and $(x, y - d / 2)$. Knowing that the native points
of $v$ are $(x \pm d / 2, y \pm d / 2)$ when $(x, y)$ is the $u$-point,
we proceed through interpolation.

$A_x(u)$ is obtained by interpolating $u$ left and right:
\begin{equation}
    \begin{split}
        % A_x(u)
        A_x(u)(x + d / 2, y)
        & = \frac{1}{2} (u(x, y) + u(x + d, y)) \\
        % & = - \frac{1}{2} (\sin(x) + \sin(x + d)) \cos(y) , \\
        & = - \frac{1}{2} (\sin(x) \cos(y) (1 + \cos(d)) + \cos(x) \cos(y) \sin(d)) , \\
    \end{split}
\end{equation}
$A_x(v)$ is obtained by interpolating $v$ left and right:
\begin{equation}
    \begin{split}
        % A_x(v)
        A_x(v)(x, y + d / 2)
        & = \frac{1}{2} (v(x - d / 2, y + d / 2) + v(x + d / 2, y + d / 2)) \\
        % & = \frac{1}{2} (\cos(x - d / 2) + \cos(x + d / 2)) \sin(y + d / 2) , \\
        % & = \cos(x) \cos(d / 2) (\sin(y) \cos(d / 2) + \cos(y) \sin(d / 2)) , \\
        & = \frac{1}{2} \cos(x) (\sin(y) (1 + \cos(d)) + \cos(y) \sin(d)) , \\
    \end{split}
\end{equation}
$A_y(u)$ is obtained by interpolating $u$ up and down:
\begin{equation}
    \begin{split}
        % A_y(u)
        A_y(u)(x, y + d / 2)
        & = \frac{1}{2} (u(x, y) + u(x, y + d)) \\
        % & = - \frac{1}{2} \sin(x) (\cos(y) + \cos(y + d)) , \\
        & = - \frac{1}{2} \sin(x) (\cos(y) (1 + \cos(d)) - \sin(y) \sin(d)) , \\
    \end{split}
\end{equation}
Now we can assemble the quadratic terms.
\begin{equation}
    \begin{split}
        % A_x(u) A_x(u)
        (A_x(u) A_x(u))(x + d / 2, y)
        % & = \frac{1}{4}
        % (
        % \sin^2(x) \cos^2(y) (1 + \cos(d))^2 \\
        % & + \sin(2 x) \cos^2(y) (1 + \cos(d)) \sin(d) \\
        % & + \cos(x)^2 \cos^2(y) \sin(d)^2
        % ) \\
        & = \frac{1}{16}
        (
        (1 - \cos(2 x)) (1 + \cos(2 y)) (1 + \cos(d))^2 \\
        & + 2 \sin(2 x) (1 + \cos(2 y)) \sin(d) (1 + \cos(d)) \\
        & + (1 + \cos(2 x)) (1 + \cos(2 y)) \sin(d)^2
        ), \\
    \end{split}
\end{equation}
\begin{equation}
    \begin{split}
        % A_y(u) A_x(v)
        (A_y(u) A_x(v))(x, y + d / 2)
        % & = - \frac{1}{4}  \\
        % & \cos(x) (\sin(y) (1 + \cos(d)) + \cos(y) \sin(d)) \\
        % & \sin(x) (\cos(y) (1 + \cos(d)) - \sin(y) \sin(d)) \\
        % & = - \frac{1}{4}  \sin(x) \cos(x) ( \\
        % & + \sin(y) \cos(y) (1 + \cos(d))^2 \\
        % & - \sin^2(y) \sin(d) (1 + \cos(d)) \\
        % & + \cos^2(y) \sin(d) (1 + \cos(d)) \\
        % & - \sin(y) \cos(y) \sin^2(d) \\
        % & ) \\
        % & = - \frac{1}{16} \sin(2 x) ( \\
        % & \sin(2 y) (1 + 2 \cos(d) + \cos(2 d)) + \\
        % & 2 \cos(2 y) \sin(d) (1 + \cos(d)) \\
        % & ) ,
        & = - \frac{1}{16} \sin(2 x) \sin(2 y) (1 + 2 \cos(d) + \cos(2 d)) \\
        & \phantom{=} - \frac{1}{8} \sin(2 x) \cos(2 y) \sin(d) (1 + \cos(d)),
    \end{split}
\end{equation}

The discrete convective term at a point $(x, y)$ is defined as
\begin{equation}
    \begin{split}
        % Cxx
        C_{x x, d}
        & = \frac{1}{d} \left( (A(u) A(u))(x + d / 2, y) - (A(u) A(u))(x - d / 2, y) \right) \\
        & = \frac{1}{4 d} \sin(2 x) (1 + \cos(2 y)) \sin(d) (1 + \cos(d)) , \\
    \end{split}
\end{equation}
\begin{equation}
    \begin{split}
        % Cxy
        C_{x y, d}
        & = \frac{1}{d} \left( (A_y(u) A_x(v))(x, y + d / 2) - (A_y(u) A_x(v))(x, y - d / 2) \right) \\
        & = - \frac{1}{4 d} \sin(2 x) \cos(2 y) \sin(d) (1 + \cos(d)), \\
    \end{split}
\end{equation}
\begin{equation}
    % Cxx + Cxy
    \operatorname{conv}_{x, d}(u, v) =
    C_{x x, d} + C_{x y, d} =
    % \frac{1}{4 d} \sin(2 x) \sin(d) (1 + \cos(d)) =
    \frac{1}{4} \sin(2 x) (\sinc(d) + \sinc(2 d)) ,
\end{equation}
A similar derivation for the $y$-component gives
\begin{equation}
    % Cyx + Cyy
    \operatorname{conv}_{y, d}(u, v) =
    C_{y x, d} + C_{y y, d} =
    \frac{1}{4} \sin(2 y) (\sinc(d) + \sinc(2 d)) .
\end{equation}

Using the Taylor series expansion \eqref{eq:taylor_sinc} for the expression
\begin{equation}
    \sinc(d) + \sinc(2 d) =
    % 1 - \frac{1}{6} d^2 + \frac{1}{120} d^4 + \mathcal{O}(d^6) +
    % 1 - \frac{4}{6} d^2 + \frac{16}{120} d^4 + \mathcal{O}(d^6) =
    2 - \frac{5}{6} d^2 + \frac{17}{120} d^4 + \mathcal{O}(d^6),
\end{equation}
and remembering that the continuous convective $x$-term is
$\frac{1}{2} \sin(2 x) $, we get
\begin{equation}
    \operatorname{conv}_{x, d}(u, v) =
    \frac{\partial }{\partial x} (u u) +
    \frac{\partial }{\partial y} (u v) +
    \mathcal{O}(d^2),
\end{equation}
and similarly for the $y$-component. This is just the confirmation that the computation of the derivatives is second order accurate.

Next, we compute the filtered convective force
$\Phi \operatorname{conv}_x(u, v)$. We get 
\begin{equation}
    \begin{split}
        \Phi \operatorname{conv}_x(u, v)
        & =
        \frac{1}{(2 n + 1)^2}
        \sum_{i = -n}^{n}
        \sum_{j = -n}^{n}
        \frac{1}{4} \sin(2 (x + i d)) (\sinc(d) + \sinc(2 d))
         \\
        % & =
        % \frac{1}{4} 
        % (\sinc(d) + \sinc(2 d))
        % \frac{1}{2 n + 1}
        % \sum_{i = -n}^{n}
        % (\sin(2 x) \cos(2 i d) +
        % \cos(2 x) \sin(2 i d)) \\
        % & =
        % \frac{1}{4} 
        % \sin(2 x) 
        % (\sinc(d) + \sinc(2 d))
        % \frac{1}{2 n + 1}
        % \sum_{i = -n}^{n}
        % \cos(2 i d) \\
        & =
        G_{n, 2 d} \operatorname{conv}_{x, d}(u, v),
    \end{split}
\end{equation}
where the transfer function is now applied on $2 d$ instead of $d$.
Noting that $\Phi u = G_{n, d}^2 u$ and $\Phi v = G_{n, d}^2 v$ everywhere on
the domain, we get
\begin{equation}
    \begin{split}
        A_{x, \Delta}(\Phi u) A_{x, \Delta}(\Phi u) & = G_{n, d}^4 A_{x, \Delta}(u) A_{x, \Delta}(u), \\
        A_{y, \Delta}(\Phi u) A_{x, \Delta}(\Phi v) & = G_{n, d}^4 A_{x, \Delta}(u) A_{y, \Delta}(v), \\
        \operatorname{conv}_{x, \Delta}(\Phi u, \Phi v) & = G_{n, d}^4 \operatorname{conv}_{x, \Delta}(u, v)
        = \frac{1}{4} G_{n, d}^4 (\sinc(\Delta) + \sinc(2 \Delta)) \sin(2 x) ,
    \end{split}
\end{equation}

Using a similar derivation for the $y$-component,
we can now compute the \emph{discrete} convective commutator error
\begin{empheq}[box=\boxstyle]{equation}\label{eq:commutator_discrete}
    \begin{split}
        \Phi \operatorname{conv}_{x, d}(u, v) -
        \operatorname{conv}_{x, \Delta}(\Phi u, \Phi v)
        & = -E \frac{1}{4} \sin(2 x), \\
        \Phi \operatorname{conv}_{y, d}(u, v) -
        \operatorname{conv}_{y, \Delta}(\Phi u, \Phi v)
        & = -E \frac{1}{4} \sin(2 y),
    \end{split}
\end{empheq}
where
\begin{equation}
    \begin{split}
        E
        & =
        G_{n, d}^4 (\sinc(\Delta) + \sinc(2 \Delta)) -
        G_{n, 2 d} (\sinc(d) + \sinc(2 d))
    \end{split}
\end{equation}
is the commutator coefficient.

% \subsection{Comparison}

It is important to note that at any given point $(x, y)$, the value of the
discrete and continuous convective commutator errors
\eqref{eq:commutator_discrete}
and
\eqref{eq:commutator_continuous}
are different.
The dependence on $(x, y)$ is the same, but
the coefficient in front is different. If a closure model is imposed using the
traditional route of ``filtering first'', but trained using commutator error
target data computed using DNS (``discretizing first''), there is an
inconsistency between the model and the learning environment.

\section{Additional numerical experiments} \label{sec:other}

We here include some additional numerical experiments
to confirm that our methods work in a variety of settings.

\subsection{LES of decaying turbulence (2D)} \label{sec:decaying}

We repeat the experiment from \ref{sec:analysis_les}, but without the body
force. This results in a decaying turbulence setup, where energy dissipates over
time. 

\begin{figure}
    \centering
    \includegraphics[width=0.49\textwidth]{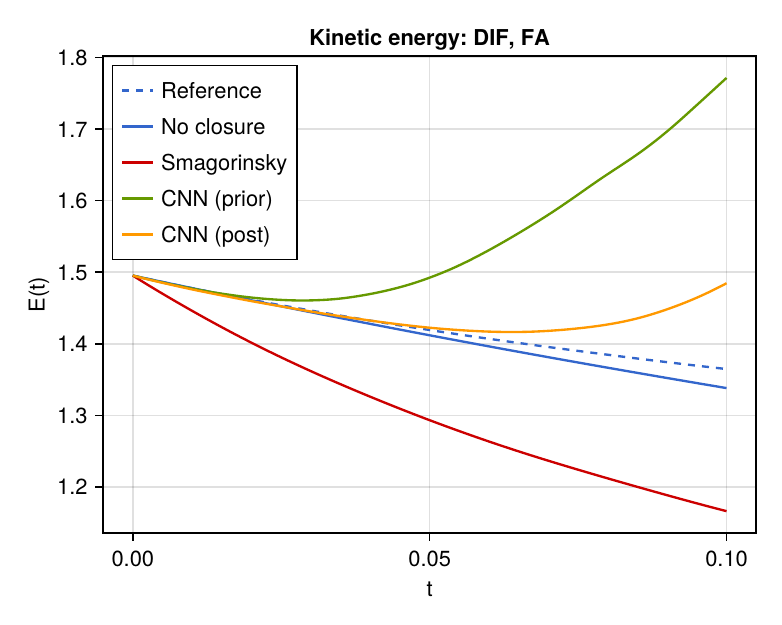}
    \includegraphics[width=0.49\textwidth]{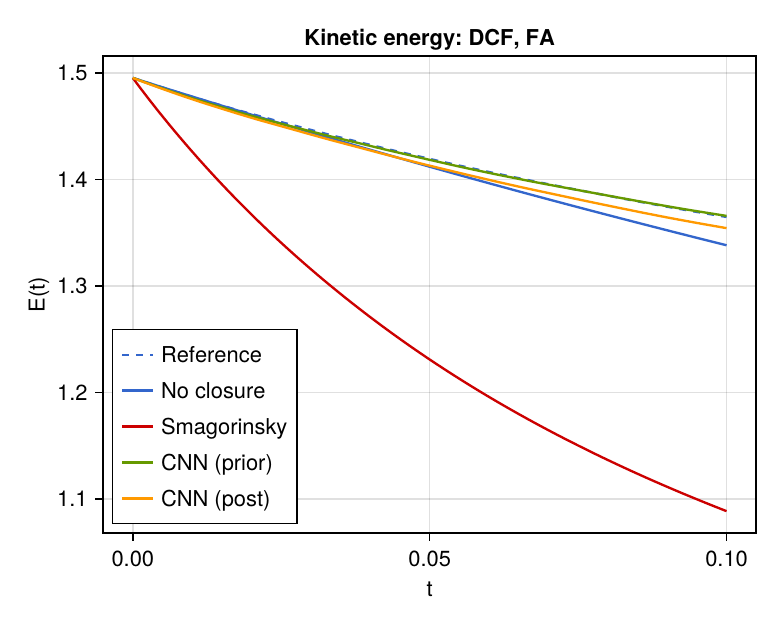}
    \includegraphics[width=0.49\textwidth]{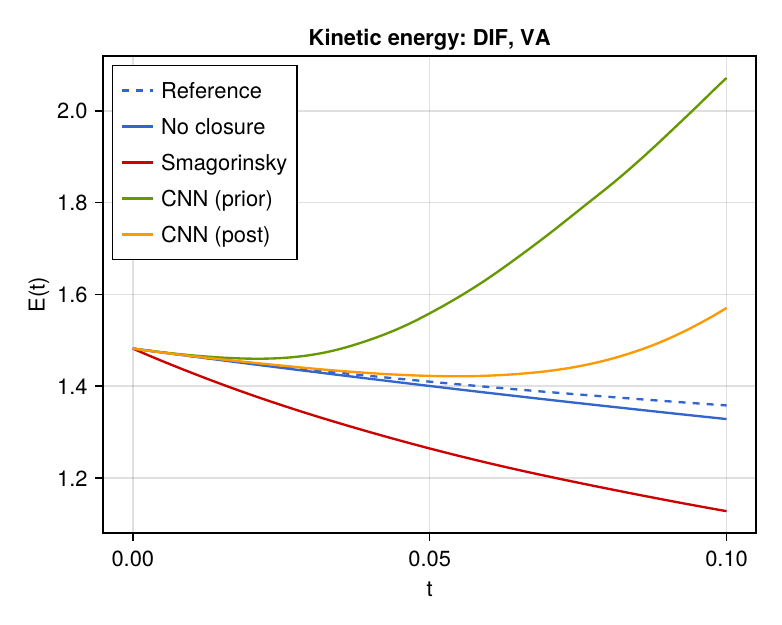}
    \includegraphics[width=0.49\textwidth]{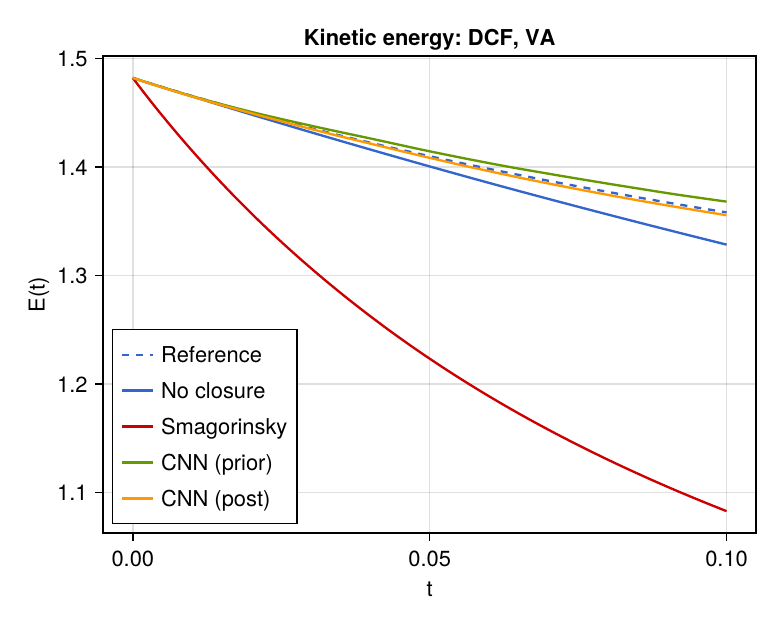}
    \caption{Total kinetic energy evolution for
        the 2D decaying turbulence case at
        $\bar{n} = 128$.
        \textbf{Left:} Unprojected closure model $\mathcal{M}_{\text{DIF}}$.
        \textbf{Right:} Constrained model $\mathcal{M}_{\text{DCF}}$.
        \textbf{Top:} Face-averaging filter.
        \textbf{Bottom:} Volume-averaging filter.
    }
    \label{fig:energy_evolution_decaying}
\end{figure}

Figure \ref{fig:energy_evolution_decaying} shows the total kinetic energy
evolution for the decaying turbulence test case. Since our discretization is
energy-conserving, the DNS energy cannot increase in the absence of a body
force, and it can only decrease due to dissipation. If $D_2$ denotes the
discrete diffusion operator, it can be shown that the kinetic energy
$E = \frac{1}{2} \| u \|_\Omega^2 = \frac{1}{2} \langle u, u \rangle_\Omega$
satisfies
\begin{equation}
    \frac{\mathrm{d} E}{\mathrm{d} t} = \langle u, D_2 u \rangle_\Omega.
\end{equation}
For $\mathcal{M}_\text{DCF}$, all models produce decaying energy profiles, but
the CNN energy is the most accurate. The no-closure model is slightly too
dissipative, and the Smagorinsky model is even more dissipative.
For $\mathcal{M}_\text{DIF}$, the no-closure and Smagorinsky models have similar
profiles as for $\mathcal{M}_\text{DCF}$, but the CNN models become unstable.
Training with the a-posteriori loss function corrects for this instability,
leading to correct energy levels for a longer time, but the instability is still
present.

\subsection{LES of forced turbulence (3D)} \label{sec:les3D}

\R{les3D}
To show that the formulation also works in 3D,
we consider a decaying turbulence test case in
a periodic box $\Omega = [0, 1]^3$ with $N = 1024^3$ finite volumes.
We use single precision floating point numbers,
Reynolds number $\mathrm{Re} = 2000$,
simulation time $t_\text{end} = 3$,
burn-in time $t_\text{burn} = 0.5$, and
$10$ random initial conditions.
Unlike for the extensive analysis in section~$\ref{sec:analysis_les}$,
we only consider $2$ LES resolutions $\bar{n} \in \{ 32, 64 \}$,
one filter (face-averaging),
and one LES model ($\mathcal{M}_\text{DCF}$).
We use the same body force as in the 2D case, namely $f$ with $f^\alpha(x) = 5
\delta_{\alpha = 1} \sin(8 \pi x^2)$.
We use the same CNN architecture as in the 2D case, but with a kernel size
$5 \times 5 \times 5$ instead of $5 \times 5$,
resulting in $234\,096$ learnable parameters instead of $45\,696$.
For the Smagorinsky model, we use the common choice of $\theta = 0.17$.

\begin{figure}
    \centering
    \includegraphics[width=0.9\textwidth]{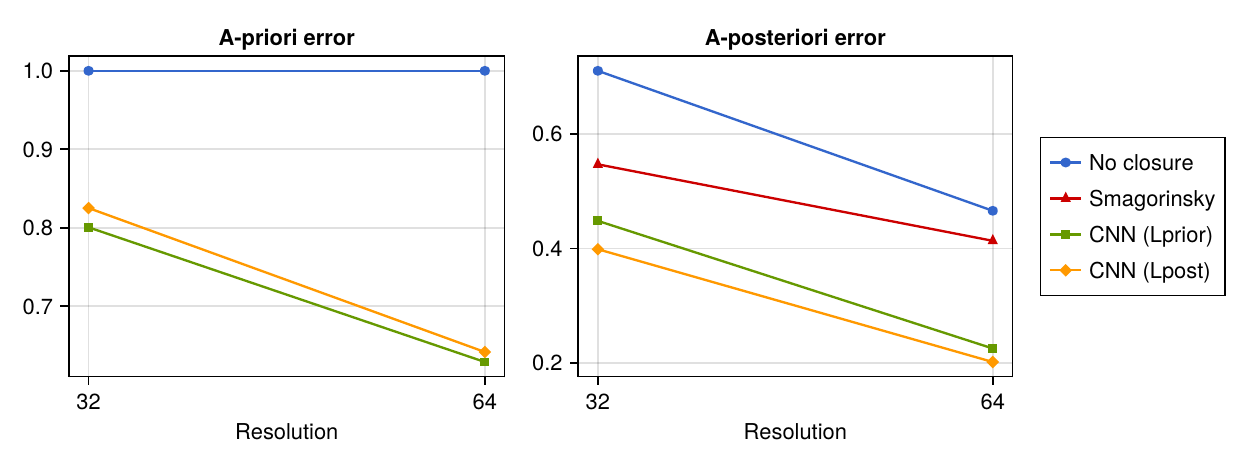}
    \caption{\revboth{
        A-priori and a-posteriori errors at time $t = 0.29$ for the 3D case.
    }}
    \label{fig:errors3D}
\end{figure}

Figure~\ref{fig:errors3D} shows the a-priori and a-posteriori errors for the 3D
case at time $t = 0.29$. As in the 2D case (figures \ref{fig:convergence_prior}
and \ref{fig:convergence}), the a-priori trained CNN produces lower a-priori
errors, while the a-posteriori trained CNN produces lower a-posteriori errors.
For both training methods, the CNN performs better than the no-closure model and
Smagorinsky.

\begin{figure}
    \centering
    \includegraphics[width=0.6\textwidth]{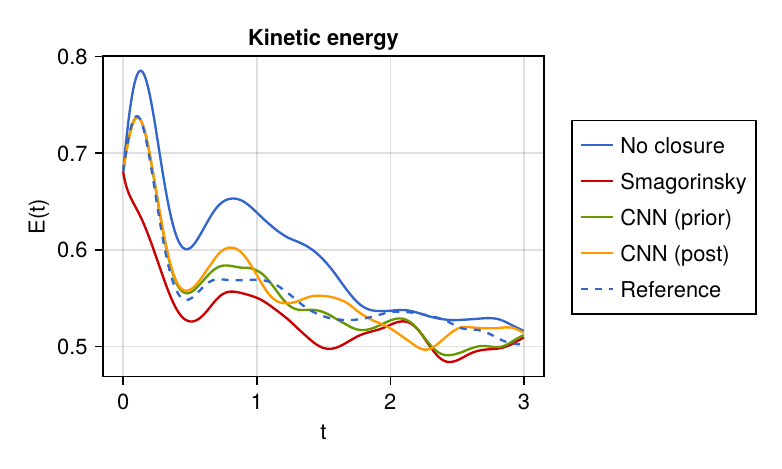}
    \caption{Energy evolution for the 3D case with $\bar{n} = 64$.}
    \label{fig:energy_evolution3D}
\end{figure}

Figure~\ref{fig:energy_evolution3D} shows the energy evolution in the 3D case.
As in the 2D case (figure~\ref{fig:energy_evolution}),
the CNN energy stays close to the target energy much longer
than the no-closure and Smagorinsky energies. Additionally, the Smagorinsky
energy seems to be more dissipative, while for the 2D case this was not visible
due to the energy injection from the body force.

\begin{figure}
    \centering
    \includegraphics[width=0.9\textwidth]{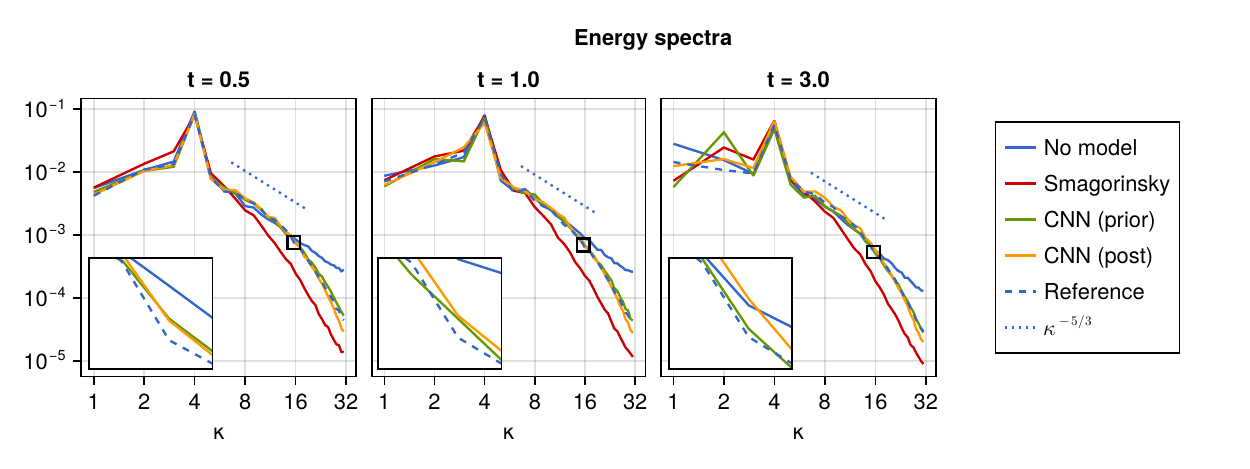}
    \caption{Energy spectra for the 3D case with $\bar{n} = 64$.}
    \label{fig:spectra3D}
\end{figure}

Figure~\ref{fig:spectra3D} shows the energy spectra in the 3D case.
A clear peak is visible at the body force injection wavenumber at $\kappa = 4$.
The Smagorinsky model is clearly too dissipative in
the high wavenumbers. The CNN is outperforming both of the two other models, and
the spectrum stays close long after the LES solution has become decorrelated
from the reference solution, confirming that while the closure model is not able
to predict the exact chaotic solution for long times, it is able to track the
turbulent statistics.

% }

% \bibliographystyle{plain}
% \bibliography{preamble_references}

\end{document}